\documentclass[a4paper,10pt]{amsart}

\usepackage[utf8]{inputenc}
\usepackage{amssymb,amsmath}

\usepackage[all]{xy}

   \usepackage[pagebackref]{hyperref}
  \hypersetup{
colorlinks=true,
linkcolor=cyan,
citecolor=cyan
}
 
\usepackage{mathtools}

\usepackage{color}

\theoremstyle{definition}

\numberwithin{equation}{section}

\newlength{\szer}
\newcommand{\Teiss}[2]{%
\settowidth{\szer}{$\displaystyle\frac{#1}{#2}$}%
\setlength{\szer}{0.5\szer}%
\left\{\hspace{\szer}%
\raisebox{0.14ex}{\makebox[0pt]{$\displaystyle\frac{#1}{\phantom{#2}}$}}%
\raisebox{-0.14ex}{\makebox[0pt]{$\displaystyle\frac{\phantom{#1}}{#2}$}}%
\hspace{\szer}\right\}%
}

\newcommand{\gras}[1]{{\mathbb #1}}

\newcommand{\Z}{\gras{Z}} 
 
\newcommand{\R}{\gras{R}} 
\newcommand{\C}{\gras{C}}

\def\elem(#1,#2){  \{ \frac{#1}  {\overline {\ #2\ }} \} }

\begin{document}

\title[An interview with Bernard Teissier]{An interview with Bernard Teissier}
\author[Patrick Popescu-Pampu]{Patrick Popescu-Pampu} 
\thanks{{\em Acknowledgments.} I acknowledge the support of the CDP C2EMPI, 
       together with the French State under the France-2030 programme, the University of Lille, 
       the Initiative of Excellence of the University of Lille, the European Metropolis of Lille for 
       their funding and support of the R-CDP-24-004-C2EMPI  project.  This work was funded 
       in part by l'Agence Nationale de la Recherche (ANR), project 
       ANR-22-CE40-0014 SINTROP.  I am grateful to Pedro Daniel Gonz\'alez P\'erez, Evelia Rosa Garc\'{\i}a Barroso, Antoni Rangachev, Paul-Emmanuel Timotei and Kroum Tzanev for their comments on previous versions of this text.}
\date{15 March 2026}



\begin{abstract}
    This text is a reworked version of a recorded interview with Bernard Teissier conducted in his house in Paris, on 28 and 29 September 2024. 
\end{abstract}

\vspace{1cm}

\maketitle

\vspace{-1cm}

\tableofcontents



\medskip
\section{Introduction}  \label{sect:intro}
\medskip

While he was my thesis supervisor in the period 1998--2001, I had weekly meetings with Bernard Teissier. They could last a whole afternoon. We discussed not only my ongoing mathematical work in singularity theory, but also many more aspects of mathematics and even extra-mathematical topics. I found those discussions fascinating. They contributed decisively in convincing me to spend my life as a researcher in singularity theory. Since then, I love discussing with him whenever possible. 

 \begin{figure}[h!] 
  \centering 
  \includegraphics[scale=0.25]{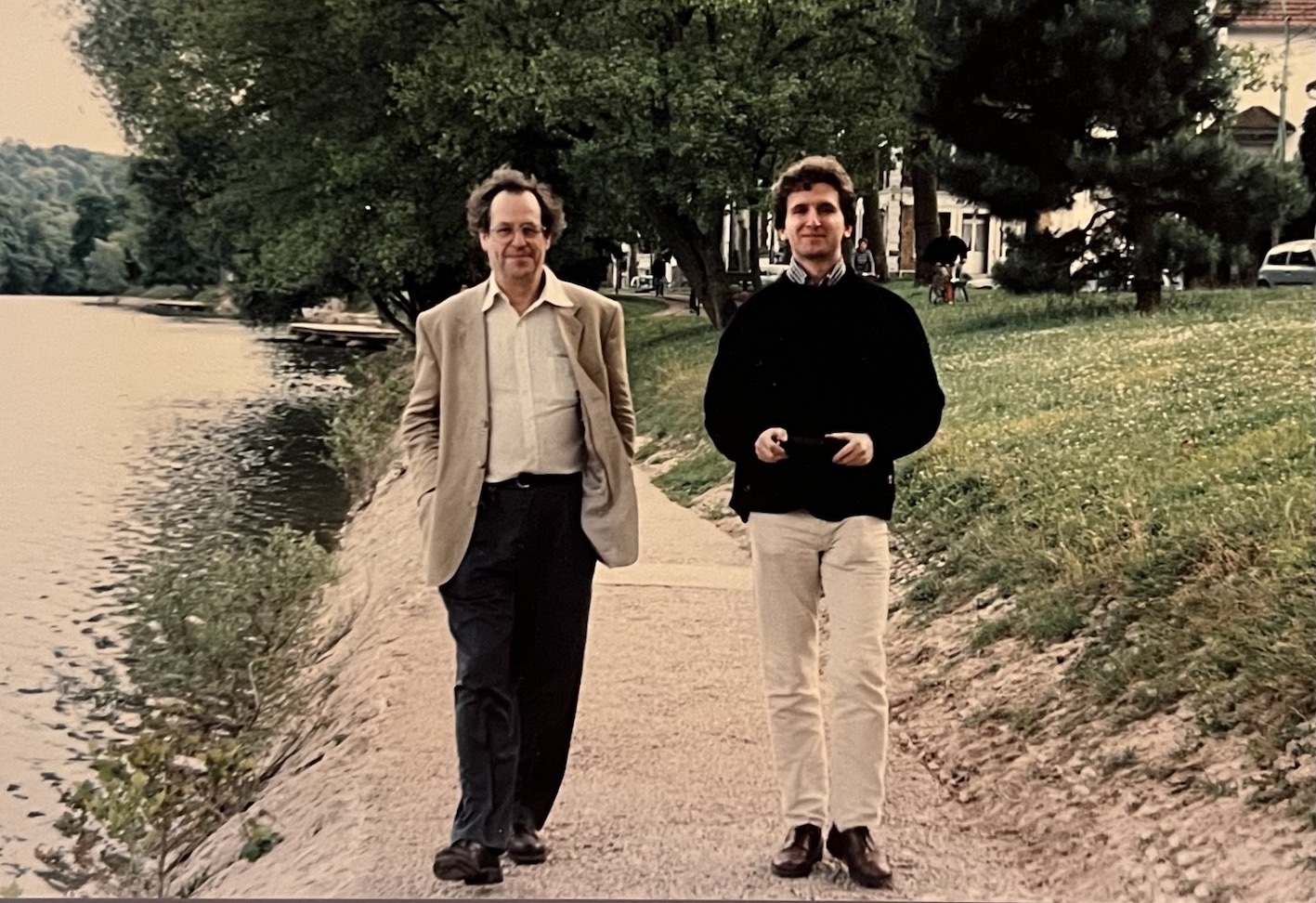} 
  \caption{Bernard and the interviewer discussing along the Marne river in June 2001}
  \label{fig:June2001}
  \end{figure} 

I conducted this interview during two days in September 2024, in order to communicate an image of his personality and of his mathematical universe to people who did not have the chance to have long discussions with Bernard\footnote{There is another interview of Bernard Teissier, conducted by Peter Galison, see \cite{G 07}.}. The text is structured into sections whose titles should help the reader orient among the various topics we discussed. We added footnotes in order to bring more information.

\begin{figure}[h!] 
  \centering 
  \includegraphics[scale=0.85]{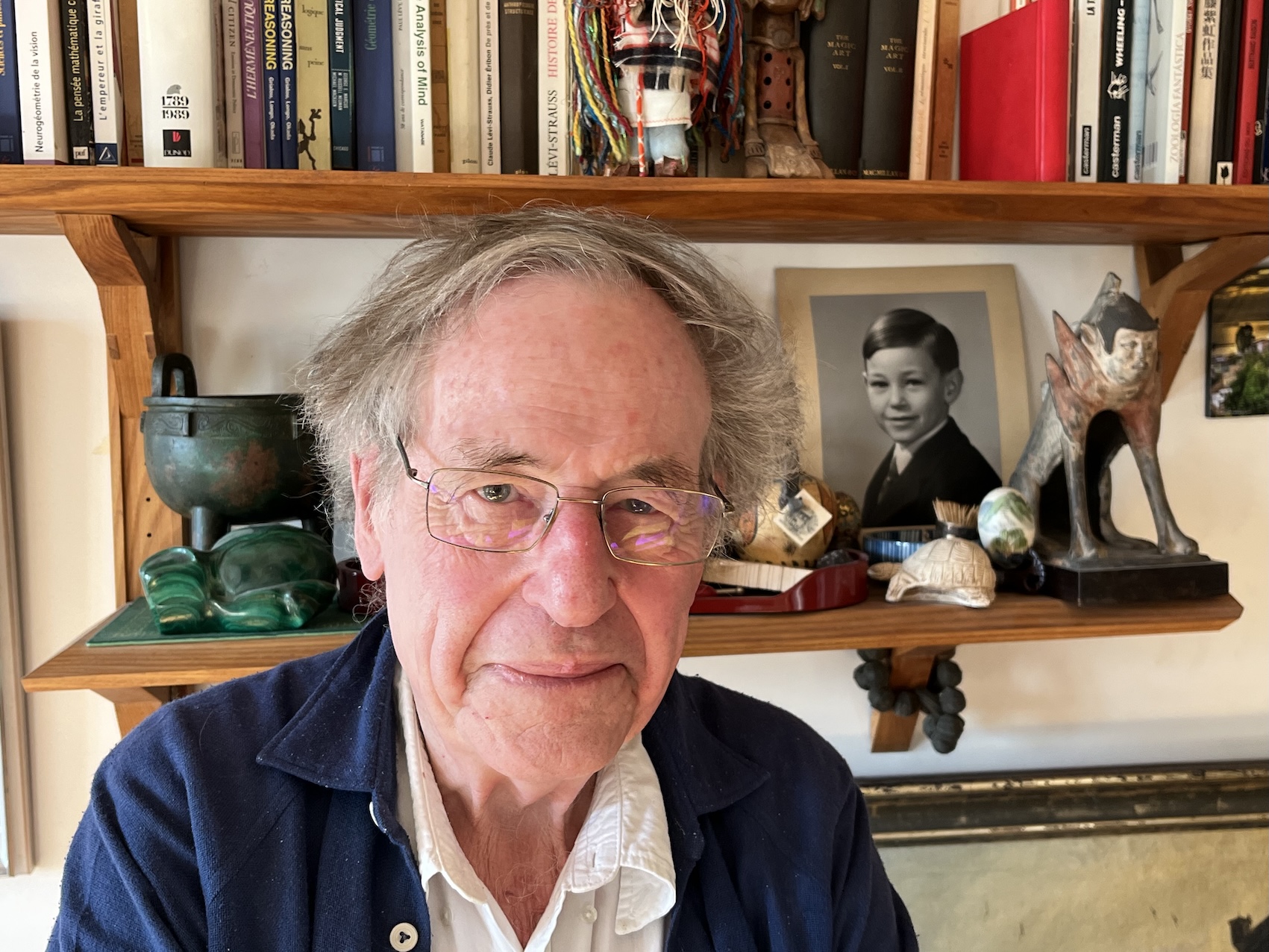} 
  \caption{Bernard near a photo of him as a child, during the interview}
   \label{fig:Sept2024}
  \end{figure}

\medskip
\section{Parents and studies} 
\label{sect:famX}
\medskip

\medskip
PPP: Dear Bernard, in this interview I will ask you questions about your life, your research interests, your works and your interactions with other persons. Let us start with your parents.  What could you tell us about them?

\medskip
BT: My mother belonged to a generation of women who had to fight to get higher education. She went to the university and became a journalist, writer and translator. She wrote for example the subtitles of many of the Marx Brothers films, but she also wrote historical and detective novels. The first under her maiden name of Jeanine Delpech, the second under various pseudonyms which you can find on her wikipedia page. She was a very capable woman. My father had no formal education, he had a complicated youth, because his father died very young. He was partly educated in the United States. He worked as a businessman, but that was not really what he liked. He was much more inclined to live in the literary world, where my mother lived. So he worked for various publishers, among other things. I admire my father because he had a gift for imagining the future in many practical ways and I've always admired that capacity in him and in other people.  I was fortunate to have very kind and understanding parents. 

\begin{figure}[h!] 
  \centering 
  \includegraphics[scale=0.5]{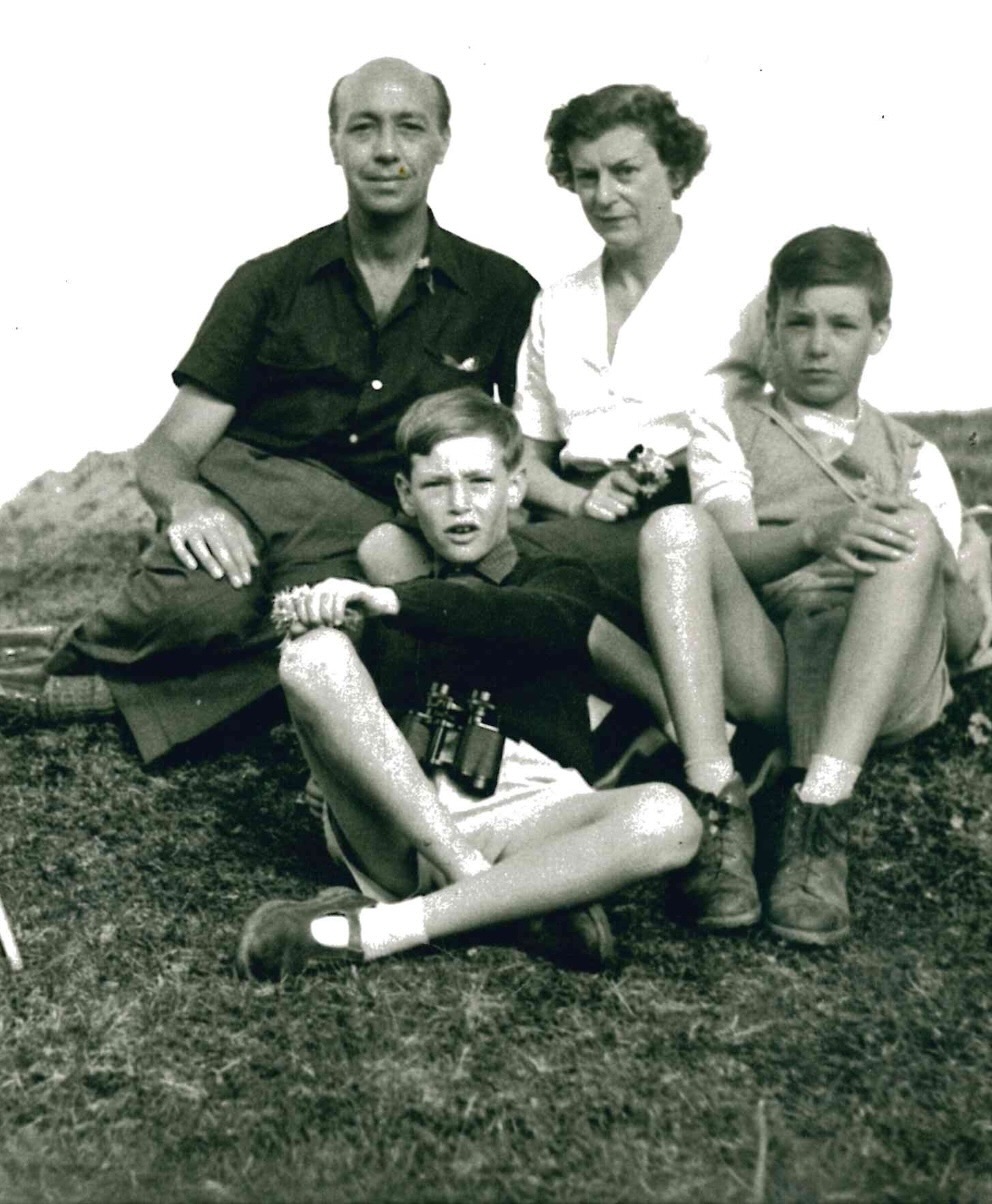} 
  \caption{Bernard (in the front) with his parents and his brother Jean-Renaud in the Alps, around 1955}
   \label{fig:Bernard-parents-frere}
  \end{figure}

\medskip
PPP: When did you begin to like mathematics, and for which reason? 

\medskip
BT: I think fairly late during my studies, perhaps one year before the baccalaur\'eat. The reason was that I had a very good teacher, Monsieur Spitzm\"uller. In fact, I think that for three years in a row I had the same mathematics teacher. I don't remember his teaching in detail, but  what I remember is that he wanted the class to laugh. He always tried to make us laugh, about many things. So, for me, studying mathematics was a pleasure. Even if I was not particularly good in maths at that time, I have a memory of being happy in class. Independently perhaps of what I understood.

\medskip
PPP: And where was it?

\medskip
BT: That was in the Coll\`ege Stanislas, because my maternal grandmother was a catholic bigot, and she induced my parents to put my brother and me in this catholic school. 

\medskip
PPP: Are you fond of other teachers too? 

\medskip
BT: In my secondary studies... yes, there was one, a Latin teacher, when I was much younger. I was maybe fourteen. He was a priest.  He was a very kind person. I remember his name as Giraudon. 

\medskip
PPP: Do you still read Latin? 

\medskip
BT: Oh yes, I can still read Latin. 

\medskip
PPP: You studied it for a long time? 

\medskip
BT: I studied it for six years, and Greek for five years. I was in a class where you had Latin and Science, but you could also have Greek. I like Greek much more than Latin. 

\medskip
PPP: Do you understand why? 

\medskip
BT: Greek is much more like a wild garden as opposed to a french-style garden. Greek is more difficult, in a way, to translate, because there can be many subtleties in each sentence. Latin is more straightforward. So, I think I liked Greek more because there was a frequent  feeling that you might be wrong, but it was a pleasure to discover at the end of the translation that you were hundred per cent right.  In Latin you did not really have that kind of surprise. Translating from ancient Greek teaches you not only their culture but also, because you have to constantly criticise your attempts at finding meaning until you finally get it right, it helps you to become a mathematician.

\medskip
PPP: Since then, did you keep reading texts in Latin or Greek? 

\medskip
BT: Not really, no. I found a mistake in a translation of Plato, in the Bud\'e edition of {\em The Republic}, but that's all. 

\medskip
PPP: What were your ideas about a future profession at the end of high-school? 

\medskip
BT: At the end of high-school I wanted to become a Hellenist. I did not really know what it meant, but I was ready to study Greek for the rest of my life. But once, our maths teacher gave us a problem about envelopes, which was very challenging. I spent a lot of time on it, and finally I found a solution.  I think that this was what started me really thinking about mathematics as a potential profession. 

\medskip
PPP: This was what made you enter into Classes Pr\'eparatoires Scientifiques, instead of Litt\'eraires? 

\medskip
BT: At that time I was really hesitating between Literary and Scientific Classes Pr\'eparatoires. What decided me was something rather trivial. I thought, probably mistakenly, that if I started in the scientific line,  I could always go back to the literary line if it did not work out for some reason, and that the opposite was impossible. 

\medskip
PPP: What would you like to tell us about that passage through Classes Pr\'eparatoires? 

\medskip
BT: I was in a Classe Pr\'eparatoire which was not at all in the top tier.

\medskip
PPP: It was still at Stanislas?

\medskip
BT: Yes. But it was managed by the Lyc\'ee Saint Louis. It was supposed to be an annex of it. It was in the same location, but it was completely different from the Stanislas system. It was not very demanding. Unlike many other people, I have a good memory of my passage through the Classes Pr\'eparatoires because, well, you worked a little, you thought a little, and things went on. It was not very competitive and I had some friends with whom I worked. 

I did not enter the \'Ecole Polytechnique at first try. There was some problem that we were not prepared for. That was one of the features of being in such not so good Classe Pr\'eparatoire. I was what they called {\em grand admissible} at \'Ecole Polytechnique, which means that I had good grades at the written exam and I could skip the first oral exam, but then at the main oral exam they asked me a question and I had absolutely no idea what it was about. So I failed, miserably. 

\medskip
PPP: You were interested only in \'Ecole Polytechnique?

\medskip
BT: No, no, I was also interested in \'Ecole Normale Sup\'erieure. But I did not pass. In my second year, I presented only two exams, \'Ecole Polytechnique and \'Ecole Normale. And I entered at \'Ecole Polytechnique.

\medskip
\section{First research interests} 
\label{sect:studX}
\medskip

\medskip
PPP: It was during Classes Pr\'eparatoires that you decided that you wanted to be a researcher in mathematics? 

\medskip
BT: Not really, no. I didn't know anything about it. I was just going through the system. In Classe Pr\'eparatoire, with two friends, we had started to translate Landau's book on number theory. We had an English edition, we started a French translation. I can't remember why, but I was interested in number theory. Then, at  \'Ecole Polytechnique, I studied a little more, and I got more interested in number theory. So, eventually, I was a candidate for the CNRS. At that time, the CNRS was recruiting young people who just came out of \'Ecole Polytechnique or \'Ecole Normale, and my research program was to apply probabilistic methods to additive number theory. 

\medskip
PPP: Was somebody influential in this decision? 

\medskip
BT: It was a little book by Mark Kac, about applications of probabilistic methods to number theory. I found it very beautiful. At that time there was practically no serious algebra course, everything was analysis. I didn't have the chance to get Laurent Schwartz as a teacher, because he taught only in alternate years and I was in the wrong year. So, no course really made a very deep impression on me at that time.

\medskip
PPP: You followed courses only at  \'Ecole Polytechnique, not at University?

\medskip
BT: Right. In my second year at  \'Ecole Polytechnique,  I went to the Number Theory Seminar of Pisot, Delange and Poitou, which took place at the Institut Henri Poincar\'e. I even gave a talk there, about Viggo Brun's work on twin primes\footnote{It appeared as \cite{T 66}.}. That was my only contact with the academic world, at that time.

\medskip
PPP: It was therefore at  \'Ecole Polytechnique that you imagined doing research in mathematics.

\medskip
BT: Yes. I liked to do maths. Most other students were going to do various things which did not attract me at all. So, it was more or less a choice by default, mostly because I could not imagine myself working under a boss. I was very, very fortunate in that Schwartz had obtained from the government and created two important things. One was called ``botte recherche'', it was the possibility for students at \'Ecole Normale or  \'Ecole  Polytechnique to go out without paying the fee that you have to pay if you leave the public service. In principle, if you go into any kind of non-public service when you leave \'Ecole Polytechnique or \'Ecole Normale, theoretically you have to reimburse what the School has paid for you. Schwartz had obtained that those students who wanted to go into research could delay the reimbursement for six years and if in that period they got a PhD, then there was no reimbursement at all to be made. At that time the members of CNRS had a special status and were not really civil servants. I think that this was a very important thing for many students. By chance, I was in the first few years when this system was put in place. 

The other important thing that Laurent Schwartz did, was to create a Mathematics Center of the \'Ecole  Polytechnique. My contact with research was a bit original, because it came first through the physicists. There was a Center for Theoretical Physics at \'Ecole Polytechnique, which was quite active. And then it happened that one day in the library at the old \'Ecole Polytechnique\footnote{Since its creation after the French Revolution,  \'Ecole Polytechnique was in the Quartier Latin in Paris, not far from the Panth\'eon. It moved to Palaiseau in the mid-seventies.}, I was looking at the books and I found a copy of the EGA\footnote{That is, Dieudonn\'e and Grothendieck's treatise ``{\em \'El\'ements de g\'eom\'etrie alg\'ebrique}''.}. I opened it, and at the first sentence I said ``I don't understand a word''. Then, for some reason I also borrowed a thick volume which had the notes of a course which Chevalley had given at the Institut Henri Poincar\'e. The notes had been taken by a physicist, Dimitri Fotiadi.  It was all functors, once again, things I did not understand a word of. So after a week or two, I put back the book. Some weeks later, I got a call to go to the theoretical physics lab. And there I met with this professor and researcher, Dimitri Fotiadi. He told me ``I see that you borrowed that from the library. Can you tell me why?'' Very embarrassing question. I could not really say that I borrowed it to read it. I borrowed it because it looked like a mysterious object, like a stone with strange pictures on it, or a piece of wood with an interesting shape you put on your desk or on your chimney. Something which could have a meaning but you have no idea what it is. It was that kind of interest. Then, he told me of the existence of the Centre de Math\'ematiques, saying that maybe I would be interested. I think I owe him a lot, because from that time on, it became clear to me that  I could choose the botte recherche and then I could go to this Centre de Math\'ematiques. And eventually this is what I did. So it was more, how to say, a lucky encounter.

\medskip
\section{Discovering Singularity Theory} \label{sect:firstmeetsing}
\medskip

\medskip
PPP: This encounter with Dimitri Fotiadi was therefore crucial for you. 

\medskip
BT: Yes, retrospectively it was crucial. And then it continued. I had a period of military service, during which I did research. At that time I was still thinking about additive number theory. But there was nobody in France doing that. So, it was a kind of dead end. I followed the course on algebraic number theory, but I couldn't think it was something I really wanted to do. 

Then, thanks again to Fotiadi, I met Heisuke Hironaka. Because Fotiadi and other physicists  like Fr\'ed\'eric Pham and Jean Lascoux were interested in resolution of singularities. So, they had been to a conference, in Battelle, maybe in 1966, where they had met Hironaka. And that year Hironaka was in France, at IHES, and Fotiadi defended his thesis. L\^e D\~{u}ng Tr\'ang and I went to this thesis defense, and specially to the party afterwards. By that time, we had more or less a project to rent a house somewhere in Scandinavia, and to spend a month there working together. None of us really knew what he wanted to do, but we were thinking we would read some books or some papers and we would try to understand things together. Then we met Hironaka, who looked like a young japanese, and we started talking to him, asking him about his speciality. He said ``I do a little commutative algebra''...

\medskip
PPP: What year was it, more or less? 

\medskip
BT: This was in late 1967. Then we found that he was a very congenial person to talk to, it was very easy to discuss with him. He told us that he had a seminar in Bures. So,  in the spring of 1968 we started going to his seminar. We told him at some point that we were planning to go to Scandinavia, and to work together for a few weeks. We invited him to come, because we thought he was a nice person. He accepted immediately. When we told that to Fotiadi and Lascoux, they became very excited, because they knew who Hironaka was. So, the whole thing took a completely different colour, because they found money. Finally we organised this meeting in Finland. It was partly funded by us, but mostly funded by the Center for Theoretical Physics. Nothing official, I don't know how it was managed. This kind of thing is totally impossible today.

\medskip
PPP: So, it was not an official conference, it was a meeting between friends. 

\medskip
BT: Yes, it was very much so. 

\medskip
PPP: For one week? 

\medskip
    BT: Oh no, it lasted almost a month, during the summer. And Hironaka gave us lectures every day, about resolution of singularities, in the complex analytic case.  Monique Lejeune-Jalabert, L\^e D\~{u}ng Tr\'ang and Audun Holme were there. Deshouillers was there, although he was already doing number theory. Pham was there, as well as Fotiadi and Jean Lascoux of course. Maybe we were ten mathematicians and physicists.  

\medskip
PPP: It was during that stay that you felt attracted by singularity theory? 

\medskip
BT: Yes. I had been attracted before, because of two lectures of Zariski in Paris, which got me excited.

\medskip
PPP: Do you remember the theme of those lectures?

\medskip
BT: Oh yes. They were about Zariski saturation\footnote{Zariski's main papers on this theme were \cite{Z 71}, \cite{Z 71bis} and \cite{Z 75}.}. 

\medskip
PPP: And it was them which triggered your work with Pham, about saturation\footnote{This work was published only in 2020, as \cite{PT 69}.}?

\medskip
BT: Oh yes, absolutely. This was at the time when I was not really knowing what I wanted to do, because I was not really attracted to anything besides additive number theory, in which there was nobody to talk to. Then I listened to two lectures of Zariski on saturation. I did not understand anything, but I was really enthusiastic. I remember that I started explaining to people who were with me what it was about, but I am sure I did not understand anything. But the style of Zariski, and the mixture of algebra and geometry... I thought suddenly, this is what I like, even if I don't understand it. 

After that, Samuel played a role. This must have been in the spring of 67, I guess, and Samuel convened a meeting to organise his seminar of the next year. I went there with L\^e. You know, we were more or less drifting around, going to hear things and trying to find our ways. We were very fortunate that we could  survive in the system just following our instinct. So anyway, we went to this meeting to organise the seminar.  There were various things which I did not understand, but at the end I suggested to L\^e to ask Samuel if he intended to talk about saturation in his seminar. I remember, this was in the Hermite amphith\'eatre at the Institut Henri Poincar\'e.  We went down to see Samuel, respectfully, because Samuel was Samuel. We asked ``Do you intend to talk about saturation in your seminar?'',  and he said ``Yes, good idea, you will do it''. So, we had to prepare  lectures about Zariski saturation. 

\medskip
PPP: You and L\^e. 

\medskip
BT: Yes, L\^e and me. I gave the first lecture, and L\^e gave the last one. And in the first paper, the first line started by ``Let $\mathcal{O}$ be a local domain of Krull dimension one''. I did not understand any of the words. So, that's really when I started to work on singularities, to understand what ``local domain'' meant...

\medskip
PPP: Do you remember what books you looked at? Zariski and Samuel, for instance? Or did you ask questions to specialists in order to discover this field? 

\medskip
BT: I went to the United States in 1966, to find out if I wanted to study there or not. And I decided not to do so. But in New York in a bookshop I found a used copy of volume I of Zariski-Samuel. I opened it, it looked like a very nice book, also pretty with its blue hard cover and the price was nothing, so I bought it. Maybe I started reading it at that time, and maybe that's how I learnt first some commutative algebra. Probably I also asked people who were around. I was of course greatly helped by very good people like Pham. I think that for commutative algebra, Fotiadi and Lascoux were not so much experts, but I learnt a lot from Pham...

\medskip
\section{The first article} \label{sect:firstart}
\medskip

\medskip
PPP: Your project with Pham\footnote{Namely, the project leading to the article \cite{PT 69}.} was to give a less commutative-algebraic view-point on saturation? 

\medskip
BT: No. It happened like this. I was at a meeting in Nice. This must have been late in 68, I don't remember exactly.  I gave a talk about saturation and Pham was there. He offered to drive me back to Paris, but he said ``I have to make a stop in CERN'', near Geneva. This was not the direct route. In the end, at 3 AM we were at CERN, and all night we talked about saturation.  In fact, he was explaining to me my talk, more or less. That's how our paper was born. 

\medskip
PPP: And if I am not wrong, that was the first paper in which were combined Singularity Theory and Lipschitz Geometry. 

\medskip
BT: Yes. 

\medskip
PPP: How did you connect these two things, do you remember?

\medskip
BT: I think there were a number of factors. I am not sure but I think we were aware, thanks to Hironaka, that integral dependence relations corresponded to inequalities of absolute values. 

\medskip
PPP: You had learnt this during your month in Finland? 

\medskip
BT: Yes, I suppose so. Pham was much more advanced than I was. Probably it was him who made the connection. At the time of my talk I had a problem on my mind: I was very frustrated because I had learnt about Puiseux exponents of a plane branch, and there was no intrinsic definition. You had to choose a coordinate system, and I thought that this is not nice, it is not as it should be. This was certainly one of the things we discussed during our drive in the night, what came out eventually was an intrinsic definition of the Puiseux exponents in terms of the blowing up of the cartesian product $\overline{X} \times \overline{X}$ of two copies of the normalisation of the branch $X$ along its relative cartesian product above $X$. I was very happy with that. I had the impression that something which is hidden in the mist suddenly becomes very clear. At the time, for me that was the main contribution of the article\footnote{See \cite[Section 10.6]{PT 69}.}, not its Lipschitz part.\par\noindent It strikes me that this approach of studying a germ of a singular space without presenting it \textit{a priori} as a ramified covering of a non-singular space is still with me today in my work on valuations.

\medskip
PPP: This was your first paper, before starting the PhD. 

\medskip
BT: Oh yes. 

\medskip
PPP: But why didn't you publish it then? It appeared only in 2020,  in the first textbook on Lipschitz Geometry of Singularities\footnote{This is the book \cite{NP 00}.}. 

\medskip
BT: You know,  at that time things did not work out like now. Many of the research centers communicated by preprints. There was a whole system of exchanging preprints from lab to lab. I remember very clearly that at the Centre de Math\'ematiques we thought that we write something, we make a preprint, we send it out, and that's it. We had a dedicated printing shop for that, which we shared with the Centre de Physique Th\'eorique.

\medskip
PPP: Very nice, much less complicated than now. 

\medskip
BT: Yes, you did not have to ask for referees and so on. I suppose that the system of evaluation was appropriate for the much smaller number of researchers. I mean, I was quite surprised several years later when I was in the United States, that  in the course of a conversation a rather eminent American mathematician said ``oh yes, but you have this Lipschitz interpretation''. Somehow, the preprint had travelled to his desk, or its content to his ears. I think that the things which were interesting circulated. The circulation system was much more demanding in some sense, because you had physically to send a copy to somebody and say ``Maybe you will find this interesting''. So, if things were not interesting to enough people, they just died out.

\medskip
\section{Early mathematician friends} \label{sect:earlyfriends}
\medskip

\medskip
PPP: You spoke about L\^e D\~{u}ng Tr\'ang and Monique Lejeune-Jalabert. They were very important for you. Can you remember how you met them and how developed your friendship? Would you mention somebody else of your generation who became important for you in the same years? 

\medskip
BT: L\^e had been admitted to \'Ecole Polytechnique, but then for some personal reasons he resigned. In fact, I tried to do the same thing, but I was a French citizen and he was not. When I went to my commanding officer saying ``I would like to leave'', he said ``Do you realise that this means you will spend the next three years in a regiment?''.  I said that I had not realised. Anyway. I think that L\^e was encouraged by Schwartz  to come to the Centre de Math\'ematiques.  At the time he was an assistant of Chevalley. 

Monique came from \'Ecole Normale Sup\'erieure de Jeunes Filles. This was like the \'Ecole Normale Sup\'erieure for boys,  but only for girls and with smaller promotions. Many prominent women mathematicians in France are alumni of that school, especially after Pierre Samuel became the Director of studies for Mathematics there, which was when Monique, Mich\`ele Vergne and Claire Delaroche, for example, were students,  the same years as my wife Maryvonne. The two schools were merged in 1985 if I remember correctly. In my time there were no girls at the \'Ecole Polytechnique, so Monique was the first woman mathematician to enter the Centre de Math\'ematiques. She was interested in singularities from early on. She was definitely interested in Hironaka's work, probably as a consequence of having followed a course of Jean Giraud at Orsay. 

\medskip
PPP: What was the topic of the course? 

\medskip
BT: It was about singularities, which was very original at the time. And so, she also came to Finland. 

\medskip
PPP: Therefore,  for you Zariski and Hironaka were important for discovering Singularity Theory and for her Giraud. What about L\^e? 

\medskip
BT: For L\^e it is more difficult to say. I don't know exactly,  but probably his interest arose because he gave these talks on saturation, and certainly because he followed Hironaka's seminar in Bures in the Spring of 1968.

\medskip
PPP: Jean-Jacques Risler, another great friend of yours, was already there? 

\medskip
BT: He arrived a little bit later. He had finished, or was finishing, his thesis. I think he was not in Finland. He wanted to do real algebraic geometry and he was a student of Cartan. I have a conjecture about how he came to join us. At that time, real algebraic geometry was not popular at all in France. It was like singularities, but even worse. I suspect that Cartan told Schwartz ``I have this student, he is good but I don't know what to do with him. Why don't you take him in your Centre de Math\'ematiques?'' Because we were all a little bit, you know, outcast, we were not part of the system, we were not at Bures, we were not at the \'Ecole Normale, we were not in the subjects which were popular at the time. I think that Risler came to join us because of that. Of course, we got along extremely well; we were a little bit in the margins, but just enough to be really comfortable and free.

\medskip
\section{The first stay at Harvard} \label{sect:Harvard1}
\medskip

\begin{figure}[h!] 
  \centering 
  \includegraphics[scale=1.8]{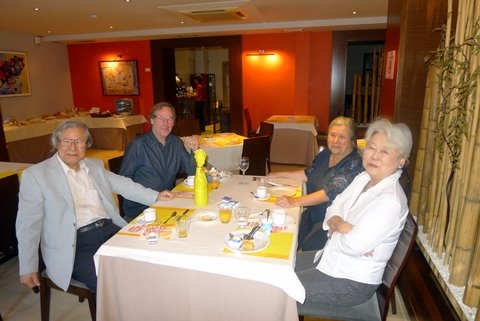} 
  \caption{Bernard and Maryvonne with Heisuke and his wife Wakako Hironaka at Tordesillas (Spain) in 2011}
   \label{fig:Tordesillas11}
  \end{figure}

\begin{figure}[h!] 
  \centering 
  \includegraphics[scale=0.25]{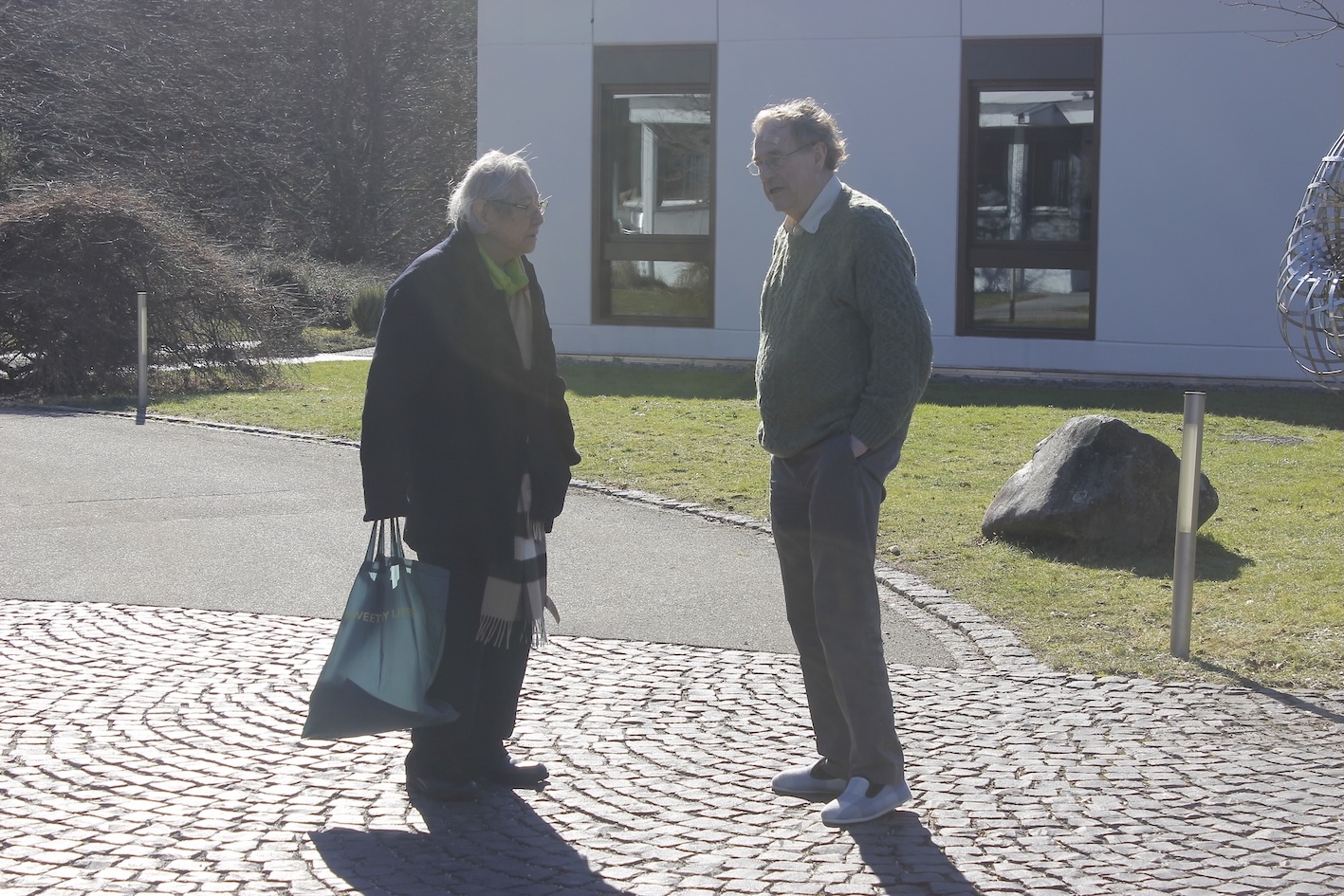} 
  \caption{Bernard with Heisuke Hironaka at Oberwolfach in 2019}
   \label{fig:Oberwolfach19}
  \end{figure}

  \medskip
PPP: Let us pass then to the United States. How did it happen that you started your PhD with Hironaka? 

\medskip
BT: What happened was that after the Finland meeting, Hironaka offered to Monique and me to come to Harvard to write down with him a book on complex analytic desingularization. Finally this book appeared only a few years ago\footnote{It is the book \cite{AHV 18}, which appeared in 2018.}. 

\medskip
PPP: Its authors are not Hironaka, Monique and you, but Hironaka and two  Spanish mathematicians, Jos\'e Manuel Aroca and Jos\'e Luis Vicente. 

\medskip
BT: But Monique and I were the first to start it. Hironaka offered us those positions at Harvard to do that writing, so we went. 

\medskip
PPP: It was not with the intention to do a PhD with him? 

\medskip
BT: No, not really. We were not graduate students but research fellows or something like that. We even had a small salary to complement our CNRS salaries.

\medskip
PPP: You did not even have the idea to make a PhD at that moment? 

\medskip
BT: I knew I had this six years limit with the ``botte recherche''. But I was not really aware what it meant to do a PhD. I knew that sometimes, in the near future, I will have to do that. But it was a remote and abstract idea. The purpose of this trip to Harvard was to write this book, not to get a PhD. We worked and we discussed with Hironaka, and of course Zariski, who was there, and many interesting people, like Mumford and Bott. It was really a great place to be. We did not finish the book, but we wrote a lot of pages. We returned to France, and it is at that time that Jean Giraud had an extremely positive influence, because we came back and we had no thesis advisor...

\medskip
PPP: After how much time did you come back? 

\medskip
BT: After one year. We went to Harvard from the fall of 1970 to June of 1971. Back from Harvard, the end of the six years was not very far. My guess is that Giraud and probably also Samuel, who I think followed Monique's work, convinced whoever needed to be convinced that we had sufficient material to write a thesis. 

\medskip
PPP: Because during this work on resolution of complex analytic spaces, you and Monique had developed a new theory. 

\medskip
BT: Yes. Starting of course from ideas of Hironaka, we had developed some pretty good, I think, technical tools, which we called the {\em theory of installations}.  It is the theory of a general kind of Newton polygons. We had also worked out the stratification of a complex space by the Hilbert-Samuel function of the local algebras, several aspects of integral dependence on ideals, rediscovered Rees valuations in the complex analytic framework and introduced fractional powers of ideals, etc. After the year 1971-72 was essentially devoted to writing our thesis, we explained part of this in our 1973-74 seminar\footnote{It was published much later, with updates, as \cite{LeT 08}.}.  And also, while I was at Harvard, I became very interested in a complex version of Thom's theory of versal deformations. 

\medskip
PPP: Had you started going to Thom's seminar before your stay in Harvard? 

\medskip
BT: I was aware of some of Thom's ideas, but I don't remember if I went to Bures to his seminar. He came to the \'Ecole Polytechnique to give a lecture, which I  liked very much. 

At Harvard, I had developed a question about the Milnor number, which obsessed me. During the week I worked at Hironaka's book, but during the week-end I worked on the geometry of discriminants and things like that. In the end, this gave rise to a part of my paper at Carg\`ese\footnote{That is, the article \cite{T 73}. It is discussed at length in Section \ref{sect:detres}.}, which ultimately was portion of my thesis. 

\medskip
PPP: During your stay at Harvard, did you meet Milnor? 

\medskip
BT: No. 

\medskip
PPP: You never interacted mathematically with him? 

\medskip
BT: No, never. 

\medskip
PPP: And with Whitney? 

\medskip
BT: No, unfortunately. I would have liked to meet him, but he was already quite old. I wrote to him when I understood something about Whitney conditions. He answered a very kind letter, saying ``I am no longer doing geometry, I am now doing didactics''. That was the extent of my interaction with him. 

\medskip
PPP: With whom did you interact seriously mathematically at Harvard? 

\medskip
BT: Mostly with Zariski and Hironaka of course. And also with Mather. Zariski organised a seminar and Mather attended it. I gave lectures on saturation, I mean Lipschitz saturation this time, and I remember Mather was very kind and said nice things to me.  We had very cordial relations with many of the members of the Mathematics Department like Mumford, Mazur and Bott, because they understood we all tried to do the same things, and they were benevolent. 

\medskip
PPP: Can you briefly describe how Zariski did mathematics? How did he choose his problems, how was it like speaking with him? Did he like telling you a lot of mathematical stories, explaining things? 

\medskip
BT: It's difficult to explain. He did not tell a lot of mathematical stories, but he made a lot of observations. For example, in a talk at his seminar in 1970-71 I presented a proof which I privately thought was quite nice. After everybody else had left he told me kindly ``That was a nice sketch of a proof''. So I went back to my proof and dissected it. I think that at that time he was still mostly interested in resolution of singularities. In fact, the theory of saturation is a subprogram of his program on resolution. He thought one should get out of the usual way to really understand what is going on. He had a vision of a method of resolution which was, I think, a variation of Jung's method, but which was much more sophisticated, using the geometry of discriminants in a much more subtle way than what comes out of the usual Jung method. More importantly, it was a precise inductive definition. Part of our interaction is described in the book of Carol Parikh\footnote{It is the biography \cite{P 08} of Zariski.}. What is not described there is how kind he was, inviting me to join his students for lunch in the Faculty Club, or in his office when he started to have hearing problems. In the 1970's and early 1980's I was visiting Harvard often and Maryvonne and I visited Oscar and Yole rather often in their Brookline apartment. One of our twins, I do not remember which one, walked for the first time in that apartment.

\medskip
PPP: Your interest in discriminants comes from Zariski? 

\medskip
BT: Not entirely. I think it comes from the observation that what Thom and Zariski were doing was all about bifurcation. Thom had his viewpoint, Zariski had his viewpoint, they are different, but the phenomenon was basically the same. You have a family of objects, varying analytically, and the shape changes when the parameters cross certain boundaries, and that is the phenomenon you have to understand.  I think that I was very  keenly aware  that these two giants were sharing the same class of problems, so to say. Zariski's approach was algebraic but one could see that it condensed a lot of geometry.

\medskip
PPP: Do you remember some precise question you wanted to attack? 

\medskip
BT: I had an obsession when I was at Harvard, what I called ``$\sum \mu_i = \mu$''. The fact that if you have a family and a singularity bifurcates, but the sum of the Milnor numbers is constant, then you can have only one singular point. This was later proved by L\^e\footnote{In the article \cite {L 73}.} using the irreducibility of the monodromy, then by several other people\footnote{By Lazzeri in \cite{La 73} and Gabrielov in \cite{Ga 74}.}. I wanted to get an algebraic proof, which still does not exist. 

\medskip
PPP: It is you who formulated this conjecture? 

\medskip
BT: Yes. 

\medskip
PPP: Could one say that it is because of this conjecture that you started working on discriminants? 

\medskip
BT: Yes. Absolutely. I really don't know how I arrived at this conjecture. It just seemed to me that it had to be that way. I tried to prove it, I tried very hard, and I failed. 

\medskip
PPP: Do you remember how you attacked the problem? 

\medskip
BT: I was trying to get a feeling for what the geometry of discriminants is.  This is a very specific problem. You can look at it in different ways. It revolves around the fact that the Milnor number is the multiplicity of the discriminant of the miniversal deformation and I was aware of the importance of equimultiplicity because of the work for Hironaka's resolution. That is why I developed this definition of the discriminant using Fitting ideals, and the decomposition formula which says that when you move away from the origin, the discriminant splits into a union whose components correspond to the various singularities in the fiber\footnote{It is explained in \cite[Chap. III, Section 2 and Section 3]{T 73}.}. That's one of the things I understood at that time.  I was trying to find rules which these discriminants have to obey, in the hope that the combination of those rules ultimately would give me the answer. I tried many things. Why should the geometry of the discriminant be the same when the multiplicity of the discriminant is constant? It comes to problems about understanding the equimultiplicity condition on the discriminant. It turns out to be a subtle problem, and Pham later gave an example showing that the geometry of the discriminant can indeed vary along the equimultiplicity stratum\footnote{See \cite[Theorem 2.9]{B 74} and \cite[Part II, Section 8]{T 74}.}. Still later, I showed that along the $\mu^{(*)}$-constant stratum, which is the same as the equimultiplicity stratum of the discriminant in the case of plane curves, it cannot vary ``too much'': some Newton polygon of the discriminant has to remain constant\footnote{See \cite[Th\'eor\`eme 6']{T 77}.}.

\medskip
\section{Research style} \label{sect:doingres}
\medskip

\medskip
PPP: So, your understanding ramified. How did you keep in memory everything you discovered? Did you have a notebook in which you marked the main ideas and proofs? 

\medskip
BT: It varies very much. Zariski had these notebooks full of examples. I was very much impressed by the fact that he could find many examples just by looking in his notebooks. I worked like that only sporadically. Sometimes, if I had a rather clear idea of what I wanted to do, then I wrote down things in a notebook. Half of the time, I just wrote things on pieces of paper, and tried to keep them in order. But I don't have a, how to say, research method, like a lab-book in a biology lab, where every experiment is described with all its conditions. You know, madness is repeating the same experiment thinking it is going to give different results. But this statement is true only if the experiment is independent of the experimenter. In the way I do research there is a lot of that. You try something, and it does not work out, and you try it with a slightly different frame of mind. Sometimes ideas are born like that. You may not have to change your computations, but you have to change the way you understand them. And sometimes you suddenly say ``Ah yes, I should have seen that from the beginning''. Of course, after a while you forget the computations and try something else.

\medskip
PPP: Do you perceive computations as an experimental part of your research? 

\medskip
BT: Yes. I do not use computers for research myself, but my collaborators do, and I make computations by hand, of course.

\begin{figure}[h!] 
  \centering 
  \includegraphics[scale=0.5]{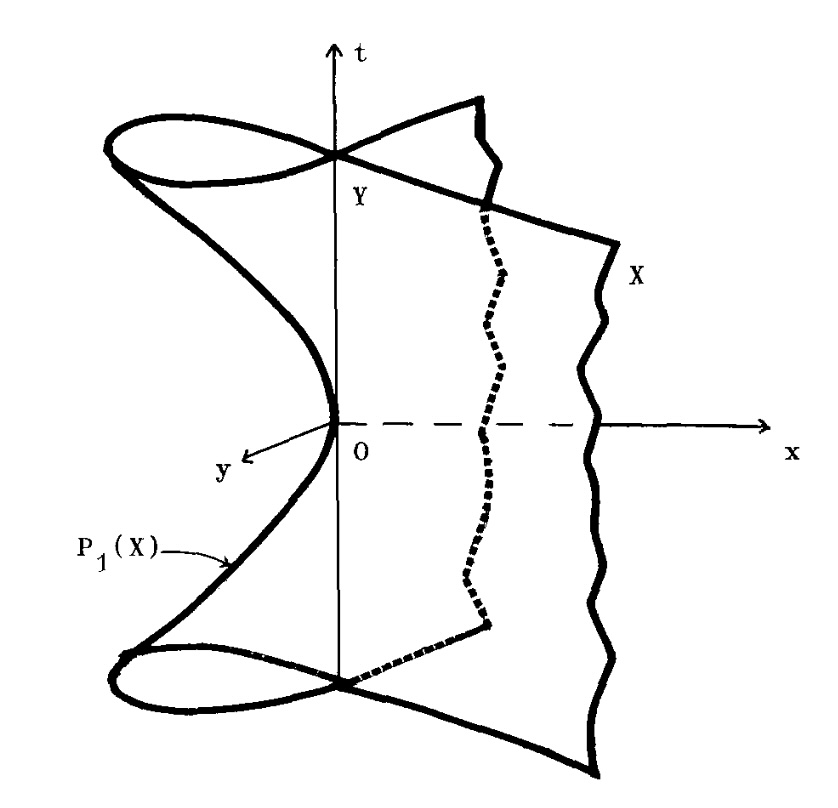} 
   \label{fig:1983-Node-degen}
  \caption{A family of nodes degenerating into a cusp}
  \end{figure}

\medskip
PPP: Do you also experiment by drawing? 

\medskip
BT: Yes, of course. I make drawings to fix some idea. I think it is necessary.

\medskip
PPP: I learnt from you a beautiful drawing\footnote{See Figure \ref{fig:1983-Node-degen}. It comes from the article  \cite{T 82}.} which appeared in your paper ``{\em Vari\'et\'es polaires II}'', with a node degenerating into a cusp, showing an example where Whitney's condition B is not satisfied.

\medskip
BT: I don't think I am the father of that example. I am sure that people like  Whitney, Thom, and Mather knew it. 

\medskip
PPP: Was it an example which was important for you in order to understand something? 

\medskip
BT: Yes. Because my main theorem in that paper\footnote{That is, \cite[Chapter V, Theorem 1.2]{T 82}.} is that the phenomenon of this example is general. I think I say it in the introduction. If you look at this picture and you understand it, then you understand the general theorem. 

\medskip
PPP: Do you think that there are other ways of experimenting in mathematics, than computing and drawing? 

\medskip
BT: That's an interesting question... I guess yes. There are more abstract ways. For example, one possibility is to detect a germ of equivalence of categories between two problems. The experiment is to try to find out whether it really is an equivalence of categories.  Where this fails, it may be very interesting. So, it is a much more abstract kind of experimentation, but I think that it occurs. 

\medskip
PPP: Do you have a precise example in mind? 

\medskip
BT: Yes. When I first became interested in the relationship between convexity and algebraic geometry, I was wondering how far it could go. 

\medskip
PPP: So, consciously you were trying to build a dictionary and you were looking for the places where it does not work. 

\medskip
BT: Yes, I thought it was interesting to understand precisely where it does not work. You have two worlds, the convex bodies on one side and algebraic geometry on the other side. You find phenomena which are, as phenomena, isomorphic. So, one attitude is to say, OK, they are both special cases of some general phenomenon which they incarnate in two different ways. The other attitude is to say that there is a dictionary between the two worlds. But in principle the dictionary cannot be perfect. So, where do the imperfections lie? I think that has an interest in itself. For example, to a convex lattice polytope you can associate a projective toric variety whose intersection theory contains much information about the polytope. But in the world of polytopes there is no concept of hyperplane section, which is essential in the toric proof of difficult results on convex polytopes. So McMullen enlarged\footnote{In the article \cite{Mc 89}.} the world of polytopes to a world where there is an analogue of the hyperplane section operation. In a more combinatorial way, June Huh and his collaborators have built cohomology theories attached to matroids\footnote{See \cite{H 18}.} in which the ``Hodge package'', which is based on hyperplane sections, holds and provides proofs of results on matroids which apply to many sorts of combinatorial questions. It is interesting to me that the first use of Hodge theory\footnote{In the guise of the negative definiteness of the intersection matrix of the exceptional divisor in a surface resolution.} to prove log concavity or convexity of a sequence was my proof\footnote{In \cite{T 77ter}.} of the inequality $\mu^{(i+1)}\mu^{(i-1)}\geq (\mu^{(i)})^2$ relating the Milnor numbers of successive sections by general vector subspaces of all dimension of an isolated singularity of a complex hypersurface $(X,0)\subset (\C^{n+1},0)$ defined by $f(z_0,\ldots ,z_n)=0$, where $\mu^{(n+1)}={\rm dim}_\C \C\{z_0,\ldots ,z_n\}/(\frac{\partial f}{\partial z_0},\ldots ,\frac{\partial f}{\partial z_n})$ is the usual Milnor number. This ubiquity of Hodge-theoretic/hyperplane sections ideas is remarkable in itself and not explained by some common superstructure. 

Of course, there is also the relation between the Hodge index theorem for surfaces and the Riemann hypothesis for curves in positive characteristic and its magnificent generalisations by Grothendieck and Deligne. I take great pleasure in the incarnation of the deep study of a simple phenomenon in {\em a priori} distant environments. For example when Yves Hellegouar'ch discovered\footnote{See \cite{H 02} and his Wikipedia page.} that you can associate an elliptic curve equation to a solution of Fermat's equation, a discovery which eventually led to the proof of Fermat's last theorem by Andrew Wiles.

\medskip
PPP: You thought seriously about the way mathematicians make discoveries...

\medskip
BT: I feel very concerned by the fact that there should be a lot of diversity in the way we discover and attack problems. I think that there is what I could call a Galois-Grothendieck way, which is to think that if you understand the structure which is behind it all,  if you can name the structure and study it, then all the problems become, not trivial but, you know, easy, or natural, whatever you want to call it. And then there is a much more experimental way, which I think is why different categories of problems are created. There is a category of problems which have to do with defining the structure, polishing the structure, applying the structure. There is another category which is to discover some phenomena, and then try to understand the phenomena like a physicist. For example, geometry of discriminants, that is a phenomenon. It follows certain rules because the deformation or the unfolding is miniversal. What are the rules? Once you have found the rules, do they suffice to characterise the phenomenon? 

Another example is the Grothendieck ring of algebraic varieties over a field, say. It is a very natural and useful object for anything concerning additive invariants, such as counting the number of points over a finite field, and of course there are many very interesting things in that direction. However, it is totally blind to the geometry of $X$ along a subvariety $Y$ when you write that $X$ is equivalent to the disjoint union of $(X\setminus Y)$ and $Y$. Contrast it with Thom's cobordism ring, which is another way of giving an algebraic structure to a set of varieties modulo some equivalence relation. It is based on the fact that a manifold with boundary has a geometry along its boundary and so there are conditions for a manifold to be the boundary of another. Two varieties are equivalent if they constitute the boundary of a third. It is also a structure, but much more to my taste because it is based on a phenomenon which can be explored geometrically, using Morse theory for example.

It is tempting to comment that this is purely subjective and that I am just saying that the right level to discover phenomena is the level at which I am more comfortable, but I believe there is more to it than that. I believe that the phenomena which we call geometric are in large part rooted in our atavic experience of the world and our perceptual system while those having to do with structure reflect more the way our brain works, its syntactic structure and the strong pulsions to organise and generalise which energise our thoughts and perceptions. Of course the two are not independent and that is what creates abstract Mathematics which has meaning for us. I am definitely not a platonist in the naive sense. I believe our mathematical concepts and objects are not shadows of some perfect inhuman objects but rather the expression of the interaction of our unconscious experience of the world and human pulsions, also largely unconscious, like the need to find causes, to organise, to compare, to predict, and many others. In short I believe that Mathematics is a human science and our impression that it is objective is due to the fact that all humans have very similar unconscious perceptions of the world and very similar unconscious mechanisms to process them.

\medskip
PPP: Do you think you combined both ways of making discoveries, or are you mainly a researcher of the second type? 

\medskip
BT: I think to be a very good mathematician you have to combine both. I feel I am more on the physicist's side. But I am an algebraic geometer in the sense that it makes me quite happy if I can find some general algebraic phenomenon which governs geometric phenomena, like the idealistic Bertini theorem which can be used to prove that Whitney conditions are satisfied outside of a closed analytic set of the small stratum and also that relative polar varieties are transversal to the kernel of the projection which defines them\footnote{See \cite[Chapter IV, Section 5]{T 82}.}.

\medskip
PPP: More, in the sense that you prefer doing that kind of mathematics? 

\medskip
BT: In the end it's not a preference. It's what comes to your mind when you think about something. 

\medskip
PPP: What you say is that, if you look back, you see that you did more the second type of mathematics. 

\medskip
BT: Yes. I have a lot of admiration and respect for mathematical activity of the first kind, because it creates very beautiful and important structures. But it turns out that I don't really work in that way. I always worked, I think, by fascination. Just like I was fascinated by Zariski's lectures, even without understanding anything. I was also fascinated by the course of Chevalley I told you about, again without understanding anything.  Like Zariski, Chevalley was obviously inhabited by what he was explaining. I think there is some kind of capacity of the mind to become interested in specific phenomena by a kind of act of faith. But to be creative you have to somehow identify yourself with some object or some phenomenon. Then, at least you ask questions. For instance, according to the historical book of Roquette\footnote{It is the book \cite{R 18}.}, the process leading to the formulation of Weil's conjectures started with work of Emil Artin, which was what we could call experimental number theory. In some special cases, of course, it was already quite elaborate, based on works of other algebraists. In my opinion, historians should also work to find the first sparks of various mathematical developments, it is very enlightening to try to see where things started and how they started. 

\medskip
PPP: It seems to me that the majority of mathematicians value much more the finished work, the theorem with the proof, and not the questions which were fertile. How do you interpret this aspect of our community? To ask a good question is in fact very difficult. 

\medskip
BT: Yes, we should not underestimate the value of questions. But on the other hand, we cannot live just with questions. It's a balance between two important things, and some mathematicians are more attracted to finding the most general statement and the most beautiful and simple proof or to show that some set of results is a manifestation of some more abstract result, of some underlying structure. This is what I called the Galois-Grothendieck way. As I said, I admire it, but I do not follow it. Every mathematician is somewhere between the two. I happen to be closer to the experimental side, where simple questions take birth, such as asking how much the geometry/topology of a hypersurface singularity determines the geometry /topology of its general hyperplane section.

\medskip
\section{Collaborations and seminars} \label{sect:collabsem}
\medskip

\medskip
PPP: Bernard, you have several collaborators. Can you describe how started some of those collaborations?

\medskip
BT: They all started with discussions, initially at the Centre de Math\'ematiques. It was very small at the time.  Basically, there was one big room, one smaller one and two small offices on the side. And so we were all more or less in the big room and we spent a lot of time just discussing what we were doing and explaining to the others. It was rather exceptional, because as I said, at that time there was no obsession with the thesis, we were just all doing mathematics together. Sometimes it was a bit more organised. I remember that, when Laudenbach gave a Cours Peccot, at the Coll\`ege de France, he wanted to rehearse, and I was his listener. So, in many sessions, he explained to me his ideas. It was not my subject of course, but it was very interesting, he was a very good explainer. That was the spirit.

With Monique collaboration was obvious, because we were both at Harvard at the same time, working on the same project. I admired very much Monique's mathematical creativity. With L\^e I can't remember exactly, we were talking about Milnor numbers... He had this manuscript of Thom which I discovered only later, relating polar curves and monodromy\footnote{\label{fn:Thommanuscr}This is the article \cite{T 72}.}. At that time I had also found a polar curve in my way\footnote{\label{fn:polarpopup}It appears for the first time in \cite[Proof of Proposition II.1.2, page 318]{T 73}, without a name. It is the curve denoted $\Gamma$ and defined by the vanishing of all the partial derivatives but one of the defining function of a hypersurface singularity.}. I was not interested by the monodromy, which was wrong, because ultimately it was topologically, using  monodromy, that my conjecture\footnote{That is, the conjecture ``$\sum \mu_i = \mu$'' discussed in Section \ref{sect:Harvard1}.} was proved. I still think that there is an algebraic proof. Anyway, we started discussing things, and little by little, sometimes something coalesced into a paper. 

\medskip
PPP: This is because you are very open in your discussions, you don't say ``I don't have to speak about this, because I want to publish it myself''. You always feel completely free to speak about all your ideas. 

\medskip
BT: Yes, absolutely. 

\medskip
PPP: This was also my impression as a PhD student, that you never told me not to speak with people. 

\medskip
BT: No, I don't like that. I know, one can on occasion regret it, but life is much more pleasant and simple if you are not permanently suspicious.

\medskip
PPP: You started very soon to organise a seminar, and you did it till now. I wonder in which measure those collaborations began related to the seminar. 

\medskip
BT: Ah, that may be the case. It started like this. When Monique and I came back from Harvard, we thought that we had to share what we had learnt. So, we started a seminar  at \'Ecole Polytechnique, because at that time there were many people interested in resolution of singularities. Risler and L\^e, but I think mostly Risler, wrote the notes from one week to the next. We just explained what we understood: maximal contact, Hilbert-Samuel stratification, this kind of things. The result was a booklet entitled ``{\em Quelques calculs utiles pour la r\'esolution des singularit\'es}''. It is on my web page. Did we start collaborations at that time? I don't remember. Perhaps that was the first time that I really met Giraud, because he came to the seminar and he asked questions. 

\medskip
PPP: So, you met Giraud only after coming back from Harvard? 

\medskip
BT: I cannot swear, but that's my impression. 

\medskip
PPP: Is it through him that you met Grothendieck? 

\medskip
BT: No, I met Grothendieck through L\^e, because L\^e was going to IHES and he took me along to this seminar, and then after two or three lectures I decided that it was not for me. 

\medskip
PPP: Was this before or after Harvard? 

\medskip
BT: I guess it was before. 

\medskip
PPP: Then you went to Thom's seminar, which you liked much more. 

\medskip
BT: Yes, but that was after the conference in Carg\`ese. Thom was there, Zariski was there, Hironaka was there of course, as well as Pham, \L ojasiewicz and Brieskorn, and many younger mathematicians. This meeting was really fantastic. Again, I have great debt to Pham, because shortly after Carg\`ese, two or three months later or so, in the fall, Thom organised a kind of mini-seminar. Pham was supposed to talk at that seminar. From what I have learnt, he told Thom ``let Bernard talk in my place''. So, I started to talk about what I understood concerning Milnor numbers. That is really where I met Thom. It was entirely thanks to Pham. I am glad that my work on Whitney stratifications has recently been used\footnote{See \cite{H-T 25} and its reference [11].} to develop his vision of Landau singularities in theoretical Physics.

\medskip
PPP: Afterwards, did you interact a lot mathematically with Thom? 

\medskip
BT: At the time, during the Thesis defence, you had to explain a mathematical paper in a field different from that of your thesis, chosen by an external adviser. I asked Thom if he could ``choose'' for me a paper of Coxeter on Phyllotaxis, which I found interesting\footnote{It is the paper \cite{C 72}.}. He accepted,  and so he was in the jury of my thesis in January 1973, with Hironaka and Verdier. I went to his seminar whenever I could. Sometimes I just went to visit him when I thought I had something to say. He was always very welcoming and very kind.

\medskip
PPP: Was he a good listener? 

\medskip
BT: Yes, indeed.

\medskip
PPP: So, he made important remarks? 

\medskip
BT: Sometimes yes. For example, he made me understand the difficulties of the real case. Because at that time I was all complex. On some occasion, because of a remark of Thom, I saw that there was a whole bunch of things I did not really understand.\par\noindent
In particular, he asked me several times to think about the catastrophist version of the Gibbs phase rule which he had enunciated and this, together with a question of Michel Herman, led me to think about the geometry of the discriminant of real miniversal deformations in a way which is explained in \cite{T 23}.\par\noindent
Once when we were talking about Whitney stratifications he explained to me his idea of using iterated jacobian extensions of the ideal of an affine complex space to derive the ideals of Whitney strata as is presented in \cite{T 67}. Afterwards he laughed and said: this is my Jugendtraum\footnote {Dream of youth, a reference to that of Kronecker.}. I believe he said that because it would give a precise relation between the position, in a sense to be defined using stratification ideas, of the map defining the space in the jet space and the geometry of the space itself.

\medskip
PPP: You spoke about Carg\`ese. If I am not wrong, the first big conference of Singularities was in Liverpool a little before and Carg\`ese was the second one. The volume of proceedings is wonderful, still now it is an important reference\footnote{It is the volume \cite{C 73}.}. Can you describe how crystallized the interest of so many people in Singularity Theory? 

\medskip
BT: The history is of course very complicated. It was Milnor's work on the exotic spheres and Hirzebruch's interpretation, then Brieskorn's work. That arouse a lot of interest among topologists, of course. Milnor's book\footnote{That is, the book \cite{M 68}.} played a very important role in attracting people to singularities, because it mixed algebraic singularity theory -- even if it is not so algebraic, the framework was algebraic -- with differential topology. That book was perhaps the tilting point leading to the search for algebraic invariants of singularities with geometrical meaning, other than the multiplicity, in high dimensions. Also, at the time, a number of people, especially in England, in Liverpool around C.T.C. Wall and in Warwick around Christopher Zeeman, were interested in the applications of the theory of singularities of mappings to differential topology and of course catastrophe theory. And you must add to that another aspect of the story, resolution of singularities, which was being done by Hironaka. Many people more or less knew that it was already done in the algebraic case, that it was interesting to understand the complex analytic case, and so on. In Liverpool it was much more the differential topology side of things which was emphasized. In Carg\`ese it was more of a mixture. A mixture of Milnor numbers, topology, monodromy, discriminants, stratifications, resolution, saturation, all these were in Carg\`ese. So, I think the fact that it became apparent that singularity theory is a crucible where many different ideas interact is the main source of the interest of people.

\medskip
PPP: You said that your interest in resolution came from your interaction with Hironaka. Concerning your interest for the differential topology side, I imagine that you read Milnor's book  immediately when it appeared. 

\medskip
BT: Yes, in fact, when Hironaka organised this mini-seminar in Bures, right after we met him, the first topic was Milnor's book, which he had in manuscript form. 

\medskip
PPP: Something which I realised rather recently, is that one says that Milnor numbers pop up from that book of Milnor. But without a name, they are already present when Thom imagines miniversal unfoldings. 

\medskip
BT: Yes, of course. 

\medskip
PPP: Was it you who realised that the same numbers occurred in both situations? 

\medskip
BT: It is Pham who realised it but I was certainly one of those who found it meaningful. I had almost comical discussions with Thom for years, because he kept denoting $\mu$ what is in fact $\mu -1$. At the time of Carg\`ese,  people were trying to prove the $\mu$-constant conjecture, which had been stated by Hironaka for curves. In general, it states that for a family of isolated singularities of complex hypersurfaces the constancy of the Milnor number implies the constancy of the embedded topological type. It was later famously proved using the $h$-cobordism theorem by L\^e-Ramanujam for hypersurfaces of dimension $\neq 2$ and by L\^e for plane curves\footnote{See the paper \cite{LR 76}, the case of plane curve singularities being treated in its Section 3.}. I realised that we did not even know if $\mu$ was a topological invariant. So, I proved that, which is easy\footnote{See \cite[Theorem 1.4]{T 73}}. L\^e proved around that time that even the monodromy is an invariant of the topological type. In times when ideas boil up, many simple things come out. And then at the end, it makes the thing whole. Yes, at the time this correspondence between dimension of the basis of the miniversal unfolding and the Milnor number fascinated me. I really thought there was something deep behind it. 

\medskip
PPP: There is a famous volume which was a product of your seminar, about Du Val surface singularities\footnote{It is the volume \cite{DPT 80}.}.  

\medskip
BT: That was at \'Ecole Polytechnique, I organised it with Demazure and Pinkham.  

\medskip
PPP: It was a  one year seminar. 

\medskip
BT: Exactly. I proposed this theme for the seminar, because I had just discovered Du Val's papers\footnote{That is, the papers \cite {DV 34}.} and again I was fascinated. I was always a little bit surprised by the correspondence between Coxeter-Dynkin diagrams and singularities, which Arnold was very fond of. When I saw Du Val's paper, I said in the end, yes, it is obvious. There's no mystery there, except that Du Val had a very specific viewpoint. He studied certain linear systems of curves in the projective plane, which corresponded to hyperplane sections of Del Pezzo surfaces, obtained by blowing up points in the projective plane. He presented their Picard group as the group of isometries of a certain polytope with respect to a metric given by the intersection form of the surface. So he made a completely algebro-geometric set of statements into statements about the symmetry of certain abstract polytopes. Then he told Coxeter ``I need to understand their symmetry groups'' and Coxeter came up with the Coxeter-Dynkin diagrams. All this is in a sequence of papers, two papers of Du Val, a paper of Coxeter, then another paper of Du Val, something like that. When I found that, I said ``Why don't we make a seminar to explain how these beautiful ideas were born?''. I am very fond, as I said earlier, of understanding the way ideas are born. In the case of Du Val, it was clearly Du Val and Coxeter who were at the root of this theory. And it's interesting that it did not really come from singularity theory. It came from the study of arithmetical properties of linear systems\footnote{Du Val was determining the intersection theory of surfaces obtained by blowing up points in the plane, interpreting minus the negative part of the intersection form of exceptional divisors as a metric in an auxiliary space. Those rational (Del Pezzo) surfaces were also resolutions of projective hypersurfaces with only rational singularities ``which do not affect the condition of adjunction''.}.

\medskip
PPP: Now that you directed your seminar for decades, do you see a change in style along time? 

\medskip
BT: Oh, yes. In the first years, let us say between mid-seventies and maybe nineties, many people came to the seminar and enjoyed it. Nowadays, the singularity seminar is not at all well-attended. I mean, the quality is there, but not any more the attendance. People say that there are too many seminars. I don't think that this is really the reason. It comes rather from the fact that we always had this seminar in the spirit of understanding ideas without any specific purpose. I remember that Schwartz told me once that when he was a young professor in Nancy, ``I was dean of the Faculty of Mathematics. That took me Monday morning from nine to eleven,  with the help of a secretary. Then, I had to teach a course on my research. And that was it, my professor duty for the week was done.'' So, he could think about many other things. The problem now is that this is like a dream. So, you cannot expect a PhD student or a young ma\^{\i}tre de conf\'erences to attend a seminar in which there is in general no hope that it is going to be immediately useful. But in those days people went to seminars out of curiosity. I think it's the style of the seminar which came out of fashion. 

\medskip
PPP: Last year, Norbert A'Campo told me that in the seventies, during the seminars people talked about what they did not understand, which was great for PhD students or for other persons looking for themes for research. He said that by contrast, nowadays people explain in complicated ways what they understand. Do you agree with this description of change of spirit? 

\medskip
BT:  Yes, it is another aspect of the change. Nowadays the evaluation system encourages speakers at seminars to try to make an impression.

\medskip
\section{PhD supervisions} \label{sect:phdsuperv}
\medskip

 \begin{figure}[h!] 
  \centering 
  \includegraphics[scale=0.25]{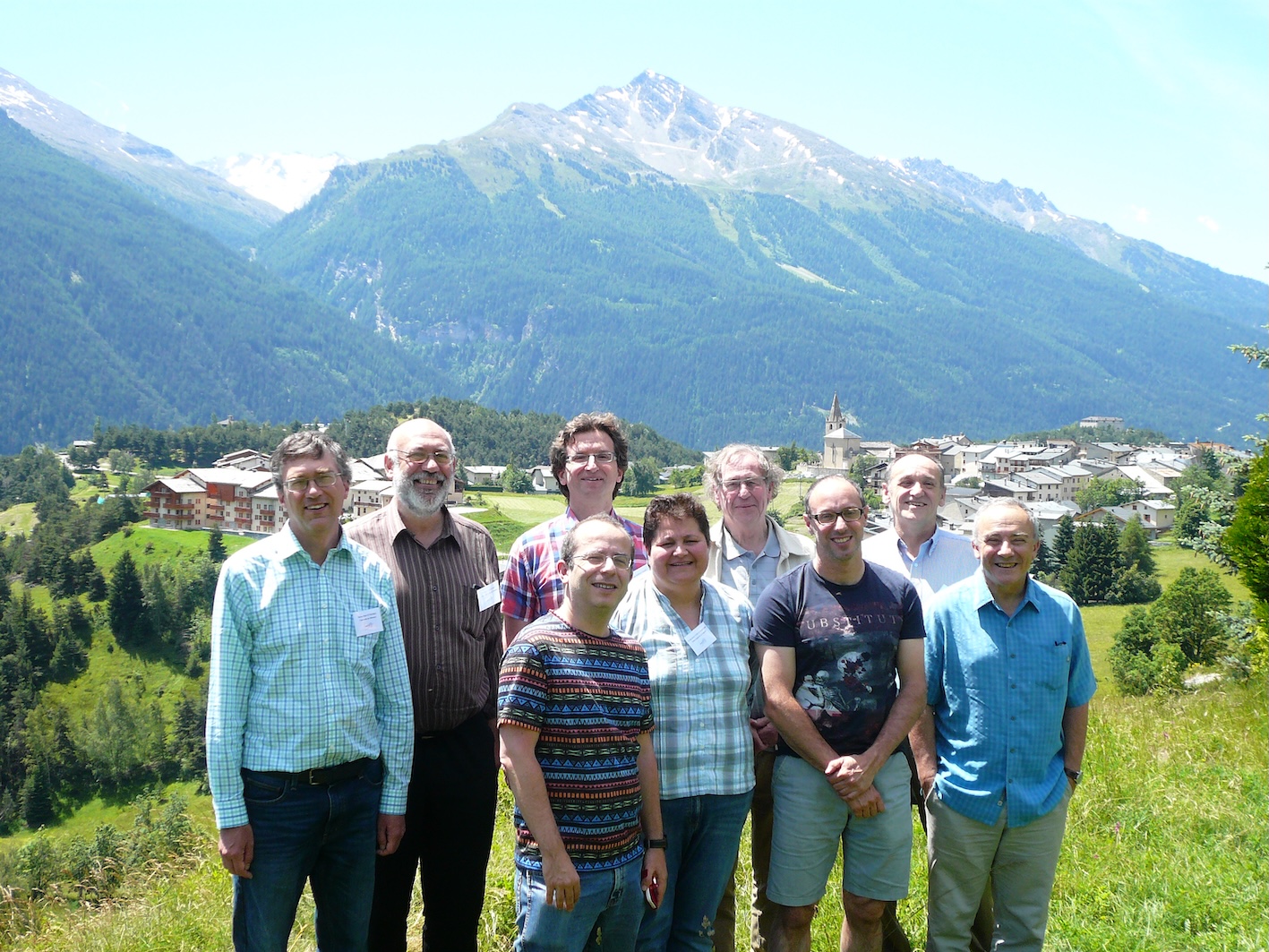} 
  \caption{Bernard with several former PhD students, friends and collaborators during the conference {\em Singular Landscapes} held in honor of his 70th birthday at Aussois in June 2015:  Dale Cutkosky, Mark Spivakovsky, Patrick Popescu-Pampu, Pedro Daniel Gonz\'alez P\'erez, Evelia Rosa Garc\'{\i}a Barroso, Bernard Teissier, Charles Favre, Ignacio Luengo, Fran\c{c}ois Loeser}
  \label{fig:Aussois2015}
  \end{figure}

\medskip
PPP: As you mention young persons, it is perhaps a good moment to speak about your interaction with your PhD students. You had many of them\footnote{Here they are, as well as the years of their defenses: Michel Merle (1974), Herwig Hauser (1980, codirection with Adrien Douady), Fran\c{c}ois Loeser (1983), Michel Vaqui\'e (1983), Evelia Rosa Garc\'{\i}a Barroso (1996), Kodjo Mawoussi (1996), Patrick Popescu-Pampu (2001), Pedro Daniel Gonz\'alez P\'erez (2000), Mohammad Moghaddam (2007), Arturo Giles Flores (2011), Ana Bel\'en de Felipe Paramio (2015, codirection with Evelia Rosa Garc\'{\i}a Barroso).}. How would you characterise your supervision style? 

\medskip
BT: I always thought that my duty to a PhD student was to introduce him to an atmosphere of research, in a certain topic.  It is the atmosphere that's important. By ``atmosphere'' I mean a set of questions and a set of things that people know, and if you think about them and you are bright, questions will pop up. Now, in some cases, like in your case or in that of Evelia, I had a specific question in mind, and I tried to share it. Every time, more than as a PhD problem, I see it as an interesting question that we can discuss. I don't have a list of PhD problems that is somewhere in my head. At any given time, I have a set of questions, some of them maybe not at all PhD questions. 

\medskip
PPP: You see them more as starting points for discussions, as entrances to a universe? 

\medskip
BT: Exactly, as entrances to a universe. Then, new questions pop up, and what I do is to say ``No, this is not a good problem for a thesis'' or ``Yes, that is a good idea''. 

\medskip
PPP: It seems to me that you introduced each of your PhD students to a universe and you let them develop there without exploring it at the same time. Was this something you consciously tried, not to ask questions which were your obsessions at that moment? 

\medskip
BT: I think you made a good description. I have obsessions at a given moment, but I don't want my students to be part of them. I don't think it's healthy. 

\medskip
PPP: I imagine that for some of your PhD students, there were moments when they were stuck, or perhaps depressed. Can you remember how you reacted?  What did you do to help them? 

\medskip
BT: I think in some occasions I said ``it's normal, don't be afraid, it's part of the process''. For some people it's more painful than for others. I never experienced that myself, because as I told you, I never really thought I had to defend my thesis before Monique and I had enough material already . But it was another era. Now, with all the pressure, all I can say is that it is normal to find it difficult.  I say ``Just work, you will make progress and you will see that a few months before your defence, suddenly you will realise that you understand much more than you thought was possible''. 

\medskip
PPP: Would you have some advice for people who start supervising PhD students? 

\medskip
BT: Advice number 1: {\em Choose your student, making it sure it's somebody you can make a friend of: don't take any students with whom you don't think you can have a confident relationship.} Advice number 2: {\em Introduce him to your way of thinking, to the set of questions that you have asked yourself, the set of things and in particular papers that you find interesting}. But, as I said, not to your current obsessions. Otherwise, it's not going to succeed. 

I always talked to my students a lot, probably more than the average. Because I enjoyed it, not out of kindness of my heart. I don't see them as a burden, that I have to see once a week and try to push to make progress on this or that. That is, my advice is to take in PhD only a person with whom you enjoy talking. And then talk! I think that then eventually, almost inevitably, it will produce a thesis. However, with the constraint of finishing a thesis in three or at most four years it is much more difficult to work that way. It is very tempting to give to your doctoral student a problem which you know how to solve, not to take risks. But risk-taking is a large part of truly original research.

\medskip
PPP: What about writing papers? What would be your advice about that? 

\medskip
BT: My only advice is to write them as late as possible. Don't rush. But of course there is a time constraint. So, it's often very difficult to write the first paper, the one you need to be published rapidly. 

\medskip
PPP: But ``don't rush'' means you are a little like Gauss, you prefer to publish ``pauca sed matura''\footnote{``Few, but ripe''.}?

\medskip
BT: Yes, except that in the end you have to be realistic. I would say this: you need to take the advice not only of your thesis advisor, but also of two or three other people. Maybe other PhD students, or other professors. You know this well, we have now in France a system called {\em Comit\'e de Suivi Individuel de Th\`ese}. It obliges each PhD student to have a secondary consultative advisor, and I think that is very good. The system is a little too formal, you have to organise a meeting once a year, and essentially the question which goes around is ``do you have a conflict with your thesis advisor?''

\medskip
PPP: This is an important question, indeed. 

\medskip
BT:  It is an important question, but it is not the only question worth asking. I think the role of the secondary advisor is also mathematical, which does not seem to come to the minds of the people who write the rules. I think it is important to have not only your thesis advisor, but also somebody else to whom you can talk sometimes more freely. You could say ``I have done this, what do you think?'' Or, ``I have written this, what do you think?'' This helps the beginner gaining confidence. That somebody who is not directly involved looks at your first writing or project or whatever else and says yes, this looks interesting to me. I mean, I think that any mathematician can detect even outside of his own speciality whether some paper has the flavour of something interesting or not. So, I am very much in favour of not remaining each one in his own corridor.

\medskip
\section{Important conferences}   \label{sect:impmeetings}
\medskip

 \begin{figure}[h!] 
  \centering 
  \includegraphics[scale=0.6]{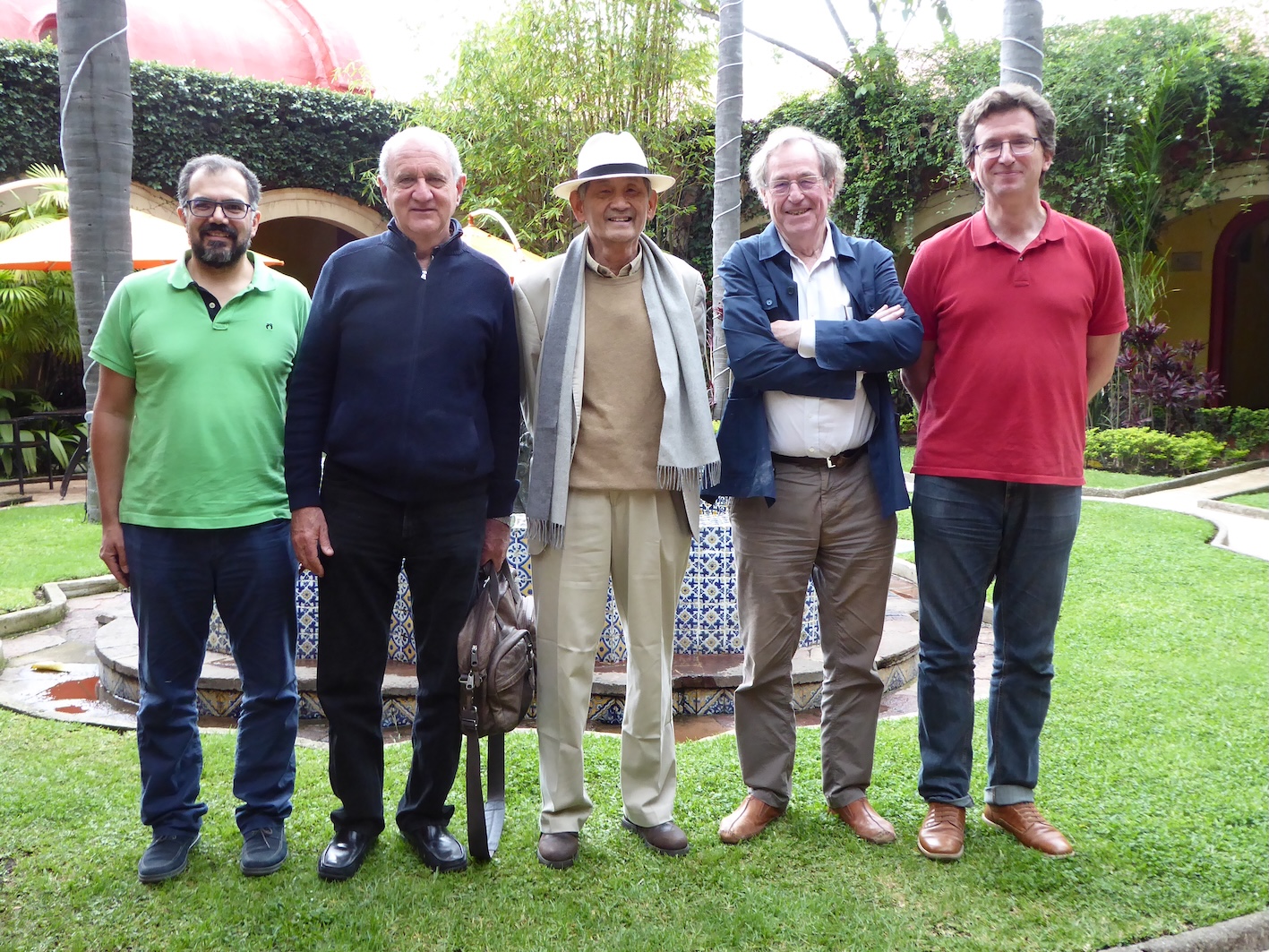} 
  \caption{Jos\'e Luis Cisneros-Molina, Jos\'e Seade, L\^e D\~{u}ng Tr\'ang, Bernard Teissier and Patrick Popescu-Pampu in June 2018}
  \label{fig:Cuernavaca2018}
  \end{figure}

\begin{figure}[h!] 
  \centering 
  \includegraphics[scale=0.8]{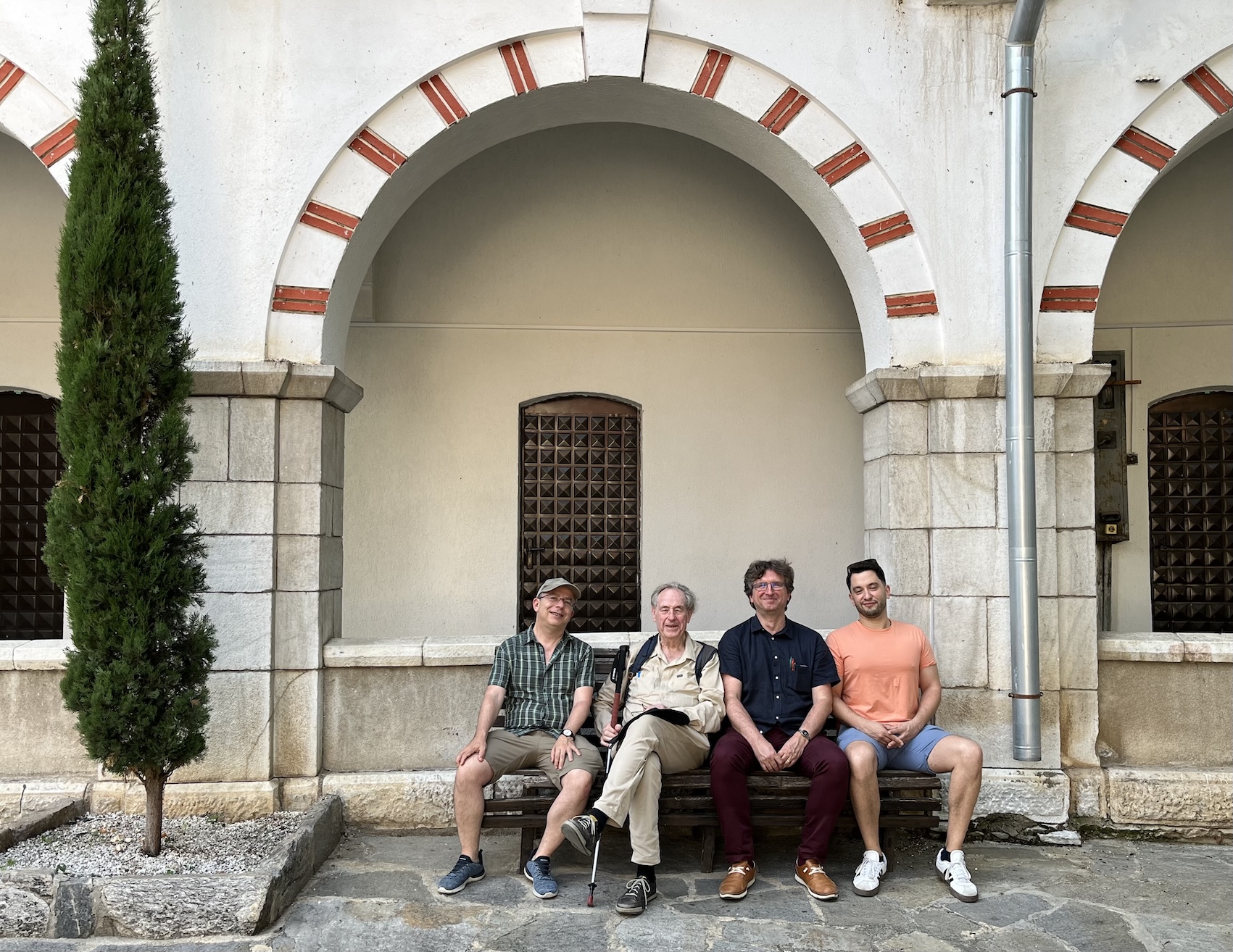} 
  \caption{Pedro Daniel Gonz\'alez P\'erez, Bernard Teissier, Patrick Popescu-Pampu and Antoni Rangachev  in June 2025}
  \label{fig:Bulgaria2025}
  \end{figure} 

\medskip
PPP: Could you tell me which conferences were particularly important for you? 

\medskip
BT: I went to many conferences. Not all of them were very significant for me. But of course, professionally they were all useful, because they allowed me to discuss with many people\footnote{In Figures \ref{fig:Cuernavaca2018} and \ref{fig:Bulgaria2025} are shown photographs taken during the {\em International School on Singularities and Lipschitz Geometry}, which took place in Cuernavaca, M\'exico, from June 11th to 22nd, 2018 and during the conference {\em Singularities in Plovdiv}, which took place in Plovdiv, Bulgaria, from June  8th to 13th, 2025.}. 

The first international meeting where I was invited is, I think, a CIME meeting in 1969 on Lago di Como. Zariski was there, as well as Benjamino Segre, Seidenberg and Mumford, but most of it went over my head. I understood a little more the next year at the 1970 Nordic Summer School in Oslo. Hironaka was there but had to leave early because of a death in his family. The first meeting which I attended where I fully participated was that in Carg\`ese in 1972. This was obviously enormously important for me, because as I told you many of the mathematicians interested in singularities were there, like Zariski, Hironaka, Thom, \L ojasiewicz, Brieskorn, and Pham, of course. The atmosphere was very free, there were lots of discussions and I learnt a lot from the lectures. So this was perhaps for me  the most important mathematical meeting which I attended. 

I think also of the meeting in Oslo, the Nordic Summer School in 1976. Once again there was a great variety of mathematicians and I was also in a special set of mind, because I was giving a course and I had to prepare it, but at the same time I had to listen to other people. So this is one thing which stays in my mind, perhaps because I was younger. 

I can think also of the second AMS Arcata meeting on singularities, in 1981. The atmosphere was relaxed, in a nice environment, with a concentration of mathematicians interested in singularities.

Then there was a meeting in Japan, in 1993, which was organised by Hironaka. A lot of mathematicians interested in singularities met there, for a long period of time. Not many people stayed so long, but many people stayed several weeks. There were not many lectures, but there was a lot of intellectual exchange between the participants. Because it was very free, you could walk, or you could talk, it was quite different from the usual schedule, when you plan lectures. We also planned lectures, but maybe one day in advance, and there were many discussions.  I made a lot of subconscious progress during that period, because it was very stimulating to be in a quiet environment, surrounded by many mathematicians for a long period of time. And also to get to know better people like Shigefumi Mori, Mutsuo Oka, Jos\'e Manuel Aroca and Miles Reid. 

There was also the meeting in Oberg\"urgl organised in 1997 by Herwig Hauser to honour the memory of Zariski. It was very lively and varied. Abhyankar was there, also Hironaka, Monique and L\^e as well as Dan Abramovich, Vincent Cossart, David Cox, Franz-Viktor Kuhlmann, Joseph Lipman, Franz Oort, Orlando Villamayor and many others.

I would also like to point out a whole series of Oberwolfach workshops. I appreciate them particularly, because I met there people which are different from those which I meet usually. Those meetings are very well organised in such a way that people exchange conjectures and ideas very freely. They are always very productive for me. 

Then there is the meeting organised in Delphi by Barry Mazur and Apostolos Doxiadis in 2007. The topic was {\em mathematics and narrative}. It produced a book called ``{\em Circles disturbed}''\footnote{This is the book \cite{CD 12}.}, which is a reference to Archimedes' sentence ``Don't disturb my circles'', which he said to the soldier who killed him. The idea was to explore the analogy between writing a proof and writing a story. This may sound a little bit crazy, but there are a lot of interesting ideas that came out of this meeting.  Interesting for the philosopher of mathematics. What happens when a mathematician writes? What does rigour mean in a novel? That is an interesting question, because rigour in that case means, I think, that the characters act according to their character/personality. Except maybe in detective stories, you cannot have a person with a certain set of characteristics, who suddenly becomes a completely different person. 

\medskip
PPP: It is very difficult to translate what you just said, for instance to French, where we have two words for the English ``character'', namely ``personnage'' and ``caract\`ere''. 

\medskip
BT: Absolutely. I found interesting to think in this way. That also struck me probably because it was a very original workshop and it was an occasion to meet people from completely different horizons of intellectual life, like people who study narration as a topic. I remember that one of them told us that when there is a film which is a blockbuster, the experts in the cinema industry follow the time-line of the action, and then produce another film with an isometric time-line. This is something about which I had absolutely no idea. Altogether, it was a very interesting meeting.  Apart from the fact that I learnt why Delphi prophets stopped prophetising, which has always been a question of interest for me. The alleged reason is that there was an earthquake. And the fumes on which the Delphi oracle sat while making predictions ceased to exist. I am quite interested in the role of oracles or far-sighted people in general, like Keynes and Szilard. 

\medskip
PPP: By the way, whom do you admire most in the history of mankind, from any domain? 

\medskip
BT: That is a difficult question. If I go back in time, I think about Archimedes, and of course the Greek philosophers, basically all of them. Then Maimonides, Fontenelle, Legendre, Gauss, Galois and so on. As I said, I also have a great admiration for people who have the capacity to guess the future to some extent, who have a very refined perception of human nature and social phenomena. They are able to say that if we do this, then this will happen, which is not obvious, but which is obvious to them. I think that it is very important to have such people. At any given time, there are very few of them. Usually we don't listen to them, but {\em a posteriori} we realise that they predicted with great accuracy what happened afterwards.

\medskip
\section{Member of Bourbaki}  \label{sect:bourbakimemb}

\begin{figure}[h!] 
  \centering 
  \includegraphics[scale=0.7]{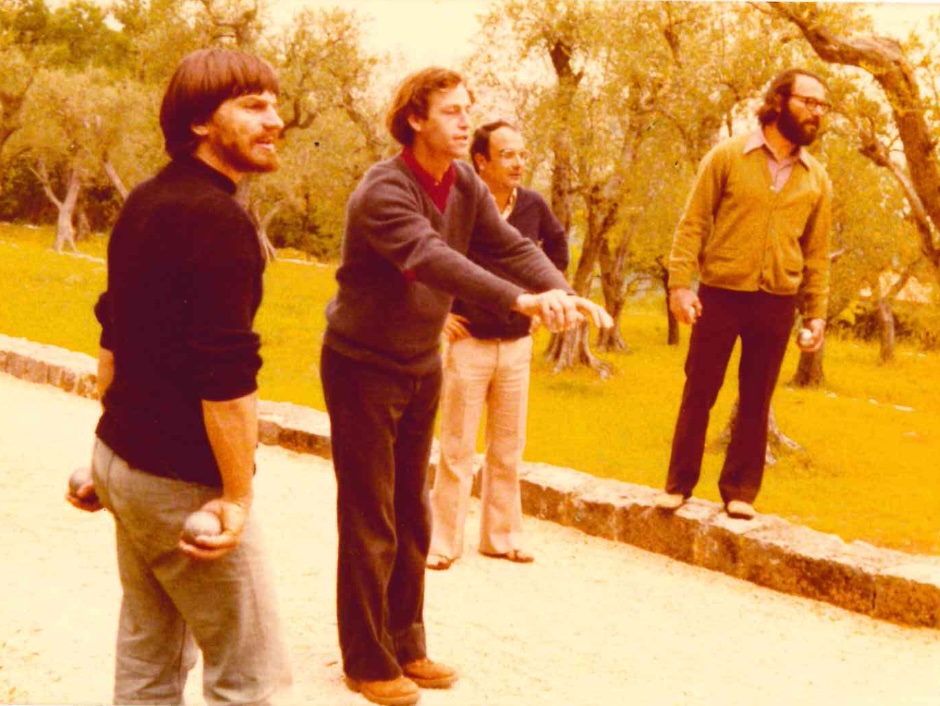} 
  \caption{Bernard playing p\'etanque during a Bourbaki meeting in Provence near Grasse, around 1975, with Jean-Louis Verdier, Andr\'e Gramain and Hyman Bass.}
   \label{fig:Bourbaki-petanque}
  \end{figure}

\medskip
PPP: Starting with the mid-seventies, you entered the Bourbaki group. I am interested by anything you can tell about this experience. For instance, how you were called in the group, or how was organised the work. 

\medskip
BT: How I was called I really don't know, but I suspect that it was because I gave some talks at the Douady-Verdier seminar at \'Ecole Normale in the early seventies. At least, Douady seemed to like them. Then I was invited to take part in one of the meetings as a ``cobaye'' as they said, that is as a guinea pig. So I went, it was in Saint-R\'emi-de-Provence. 

\medskip
PPP: Which year? 

\medskip
BT: Probably 1974. The principle of this guinea pig is that if you show some interest in what Bourbaki is reading, then they invite you again. If you still show interest, then you are invited to become a member. This is what happened with me. I found the experience very interesting, because it is captivating to see a group of people trying to clarify ideas as much as possible.  I suspect that Bourbaki hoped to write a volume on complex analytic geometry, perhaps even resolution of singularities and that this is the reason they invited me. In that sense they were disappointed. But at that time, Bourbaki was disappointed in many projects actually, because much of the time was spent with the {\em Grande \'Edition}, which was supposed to be the definitive edition. Nevertheless, I continued to take part. For me all those years were very instructive. First of all, you see a lot of mathematics and secondly, you are in this process of trying to find the good way to express things, the good objects to be defined, and so on.

\medskip
PPP: How started the work on a precise topic? By reading some articles or books? Or trying to do everything from scratch? 

\medskip
BT: No. Let me explain how Bourbaki worked. The group decides that some topic is a good target for a volume or a series of volumes. That is the subject of protracted discussions. Is it a good idea to make a volume on that topic? If the decision is yes, then Bourbaki assigns the subject to one member, who is asked to make a first draft of a first chapter in the subject. 

\medskip
PPP: That member is completely free to look at literature? 

\medskip
BT: Yes, of course. 

\medskip
PPP: And there was the rule never to include a bibliography in Bourbaki...

\medskip
BT: That's right. I can't say I know why this rule was set in place, because this was long before my time. But I understand it goes with the unanimity rule. In the end you see that there is an aim to it. Especially that many exercises are taken from published papers and if you want for each exercise to give a reference, then it has no end. 

\medskip
PPP: What happened after the first person wrote a draft? 

\medskip
BT: It was read at one of the meetings, it was criticized, and then another person or sometimes the same one is asked to write a second draft. And so on. It can take a number of years.

\begin{figure}[h!] 
  \centering 
  \includegraphics[scale=0.7]{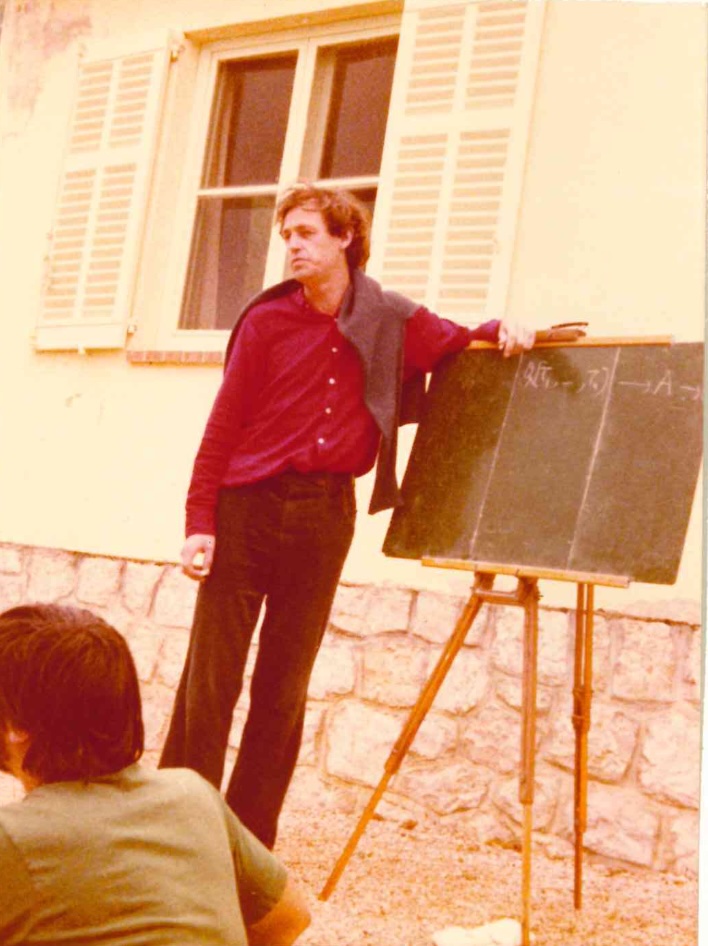} 
  \caption{Bernard during the same Bourbaki meeting, session on commutative algebra}
   \label{fig:Bernard-Bourbaki}
  \end{figure} 
\medskip

\medskip
PPP: When you say that it was read, was it done simply orally or was it a blackboard talk? 

\medskip
BT: No, purely orally. 

\medskip
PPP: I saw some pictures of Bourbaki meetings with blackboards, so I wondered how they were used. 

\medskip
BT: At the blackboard were explained some technical points, or maybe some idea to improve the text. Normally, all Bourbaki members are sitting around a table with the text in hand and someone reads.

\medskip
PPP: It was read also with Flaubert's idea of detecting what did not sound well, or the repetitions? 

\medskip
BT: Yes, there was also this aspect.

\medskip
PPP: How would you describe the philosophy of the  Bourbaki enterprise? 

\medskip
BT: I view Bourbaki as a public service. People with competence and desire write a treatise which is supposed to be a toolkit for working mathematicians. I felt that very strongly. We tried to make the best possible tools for the working mathematician with a problem in a domain covered by Bourbaki. We looked for good definitions, good proofs. 

\medskip
PPP: Would you say that the best tool is the most general one, or is this caricatural? 

\medskip
BT: I think this is caricatural. 

\medskip
PPP: How did you discuss about best possible tools? Did you look at applications and tried to see if it works in many applications? 

\medskip
BT: No, not really. In the course of the discussion, someone objected to what was written and said ``it would be better to put it this way''. It was not in terms of applications, but more in terms of writing. Of course, the applications are taken into account but what guides the writing is to seek the best path to the results, in terms of eliminating useless assumptions, introducing the useful notions in the clearest way, etc. It is not at all trying to build a metatheory, but much more eliminating the part of the environment in which the theory was born and which, {\em a posteriori,} turns out to be unnecessary. And also, if necessary, add some notions which clarify the structure of the theory.  It is more like removing diamonds from their gangue and polishing them.

\medskip
PPP: And what did you think about the users? Did you imagine them looking for a precise theorem, or as people who would learn a whole subject through the book? 

\medskip
BT: I think both. 

\medskip
PPP: Were you very careful about the global structure of the treatise? 

\medskip
BT: Yes, there was a lot of discussion about which thing should go first. But also for each specialized topic, why things should be put this way. And it was not always a quest of generality. The main objective was clarity. Sometimes, if you try to make things hyper-general, then you lose the meaning too much. 

\medskip
PPP: A legend tells that Bourbaki started because of problems with Stokes' theorem and the way it was written and proved in a textbook of Goursat. Did you also experiment with students, in order to see if it can be used in teaching? 

\medskip
BT: No. 

\medskip
PPP: So, this preoccupation was lost in the meantime, the fact that the treatise was written for teaching at university. 

\medskip
BT: I think that it was lost, and besides I am not sure that it was the problem of teaching which interested the first members of Bourbaki, it was really the problem of rigour. If you want to have a good version of Stokes' theorem, you must know what a boundary is, you must know quite a lot of things about which people were sometimes quite careless. I do not own the treatise of Goursat, but I own that of Valiron, which came out 10 years later or so, and it is true that I found it lacking in rigour many times. On the other hand, it is full of meaning. I like the book of Valiron, because despite its relative lack of rigour, it tries to be meaningful. 

\medskip
PPP: Did your writing practice evolve under the influence of Bourbaki? 

\medskip
BT: I think so, yes. 

\medskip
PPP: Could you describe the changes? 

\medskip
BT: It is difficult to say... I never had any indication on how to write. Because I never really had a thesis advisor in the modern sense. I don't think Hironaka read very critically what Monique and I wrote, except that he wanted to know what was in it. So he would not nitpick sentences or things like that, like I did with some of my students. So I am a totally self-made author, with many defects I think. I cannot really point them out but I think my writing skills improved with Bourbaki. I learnt to be more demanding on rigour.

\medskip
\section{Editorial work}  \label{sect:edwork}
\medskip

\medskip
PPP: Let us pass now to another of your interests, the editorial work. You were editor of various journals and collections of books. Could you list and compare those activities? 

\medskip
BT: It began with Ast\'erisque. I cannot remember who asked me to become editor of Ast\'erisque. When I started being editor-in-chief, I did not know much, but fortunately there was a very good secretary. Between the two of us, we did everything, meaning receiving the manuscripts, choosing the referee, handling the correspondence with the referees, even go to the printer. I discussed with the editorial board whenever I needed advice for the choice of referees. It was a lot of work, beyond what is reasonable. But it was interesting, because I saw a lot of good mathematics pass before my eyes.  I had the impression that it was useful. Ast\'erisque started a little like a French counterpart of the Springer Lecture Notes in Mathematics, which at that time was already dominant in the publication of manuscripts, conference proceedings and so on. I became quite interested in the work, because it has both an intellectual and a practical side. One has to make many practical choices, even the cover and the font.  Also, it was a period of growth for Ast\'erisque, which was still young, and that was encouraging. 

Then, Michel Herman left the editorship of the Annales de l'\'Ecole Normale Sup\'erieure. He came and asked me to be his successor. I had a lot of respect and friendship for him, so I accepted, although I thought it was a lot of work. In the end it was not so much work, because there was a full-time secretary (for Ast\'erisque I had a half-time secretary). Also the publisher Gauthier-Villars took on the practical part of the work. So, it was much less demanding, it was pure intellectual work. But it was also quite interesting. 

\medskip
PPP: After that, you passed to Springer Lecture Notes in Mathematics. 

\medskip
BT: Yes. In the meantime, I had become known to Springer because with Genevi\`eve Sureau, the co-director of the network of French Mathematical libraries (RNBM), I had organised a consortium between Springer journals and French libraries. I was well-known in Heidelberg for that reason. So, Catriona Byrne asked me one day if I would accept to be editor-in-chief of the collection in the company of Floris Takens and Albrecht Dold. Unfortunately, they both died soon, around 2010. I still like the job very much, because it has all the advantages and none of the drawbacks. We have an extremely efficient secretary, Ute McCrory. She does not let anything pass, for instance when I don't answer she asks me again and again. I have the impression that it is really useful. With Jean-Michel Morel, my co-editor-in-chief, we have exactly the same philosophy, which makes us agree on everything after discussions. He covers large parts of applied mathematics. For the last three or four years we have had an editorial board to whom we can delegate the treatment of some manuscripts. There is still some work, of course, but I enjoy very much the diversity of the mathematical landscape. 

I have the impression in all these editorial jobs, that the community of mathematics is like an orchestra without a leader. One of my friends, G\"unter Ewald, who is also a theologian, once wrote a paper about the brain, and that was his image, of ``an orchestra without a conductor''.  I don't know if neurophysiologists agree to that description. But I think it fits the community of mathematicians very well. The miracle is that we have no partition, no leader, and we make music. I am very happy to be part of that. Mostly, that's what satisfies me in the job of mathematician, it is to be part of this orchestra. Then, doing editorial work is more or less like letting the music be heard. I enjoy that part too. Sometimes there are manuscripts in very niche subjects. But if we get persuaded that they are good mathematics, then we publish them. I don't think many journals do that, because they seek a high impact factor, and we care less about this kind of thing. 

\medskip
PPP: Do you have some advice for starting editors, for instance about the way to choose referees? 

\medskip
BT: At the Lecture Notes we have the privilege that we can ask several referees at the same time. Many times I have to treat manuscripts in subjects on which I am very far from being an expert. I spend a lot of time on MathSciNet and Zentralblatt to see who publishes with whom, who reported on whom, who published what. In this way I get an impression of the frame of mind of the person I think about as a possible referee.  After a while, I make a choice. I look at six names maybe, then I make my decision. 

\medskip
PPP: You don't ask to specialists of the field of the manuscript to help you choosing the referees? 

\medskip
BT: It depends. If I can, yes of course. In the last few years we have recruited associate editors who can organise the refereeing of a manuscript and then make a recommendation. This alleviates the task very significantly. But a part of the manuscripts we receive do not fall within the field of any of the associate editors. Then it takes a little while to make sure that the person chosen as referee is not too close to the author and not too far from the subject. 

\medskip
PPP: Do you expect your referees to read line by line and to detect errors? 

\medskip
BT: No. Of course, if they do it, we are very happy. Mostly we expect them to look at the theorems, to see if they are interesting and if they are probably true. The real expert doesn't need to read the proof line by line.  He looks at the statement, at the way it is formulated, he looks maybe at several lemmas around it, and with that he can make an opinion. The main question is to know if the text is interesting and truly a monograph on the subject. This is what we expect the referees to concentrate on but of course, being mathematicians, they form an opinion about the correctness and sometimes find mistakes. 

\medskip
PPP: Do you have other advice for a starting editor? 

\medskip
BT: You must respect the author and the referee. Sometimes there is a disagreement between them. Then the editor can play a critical role. If you respect both of them, then in the end an agreement will be reached. Most of all, you have to get your own opinion. Once you have at least two referees' opinions, you have to get your own opinion by looking at the manuscript and at the reports. If I detect some problem, it is my job to make the final decision. There is no rule for accepting. Some journals ask for the unanimity of referees. I think that is not a good rule, because it may happen that you choose two or three referees and one of them does not understand what is going on.

\medskip
\section{Administration of mathematical libraries}  \label{sect:admres}
\medskip

\medskip
PPP: Let us pass now to your involvement in the administration of mathematical libraries. 

\medskip
BT: I started that when I arrived in 1987 at the \'Ecole Normale Sup\'erieure. Then Michel Brou\'e asked me to take charge of the library of the mathematics department. I had never done anything like that before, but I thought it could be interesting.  It was indeed interesting, because this was the time when mathematical libraries got digitalized. And the librarian, Madame Martin, was about to retire, but when I discussed the subject with her, she was very enthusiastic. So we did what was needed to computerize the library. That got me in touch with the other librarians. At that time the R\'eseau National des Biblioth\`eques de Math\'ematiques (RNBM), which I have already mentioned and which existed already for several years, had the very good idea that all mathematical libraries should use the same system, Texto, because that way they could exchange files. We did therefore use Texto. So I got involved not only with the library at the \'Ecole Normale, but also with the RNBM. I found a very nice community, made up of very competent people. They were very devoted to their task. There was a special spirit, due to the fact that the librarians had direct contact with the researchers. There was inside the RNBM a kind of parity between librarians and mathematicians. Which was not at all the case in other fields. The RNBM had a director, who was a librarian, and a co-director who was a mathematician. 

At the time when I had the idea of making a consortium with Springer in order to limit the prices of subscriptions, I was asked to become Scientific Director of the RNBM while the director was Genevi\`eve Sureau, the head Mathematics librarian at Orsay. That was a very interesting period. I explained my idea to Catriona Byrne and as a result Genevi\`eve and I were both invited to Heidelberg. We discussed this with many people in the hierarchy of Springer. We started negociations, which took a long time. At some point, one of the vice-directors of Springer said ``OK, we will give you free access for a year and a half, then we negotiate''. This was very clever. All mathematical libraries in France got free access to a list of Springer mathematics journals which we at the RNBM chose and of course the libraries became addicted. Nevertheless, at the end we had a very good negotiation and a good consortium agreement, which continues to this day.  I believe that everybody is happy with it, except that there have been periods when Springer was more aggressive, for instance trying to impose the whole mathematics catalog.  That's mostly coming from the fact that Springer became more commercial. The RNBM later negociated similar agreements with the European Mathematical Society Publishing House and the Soci\'et\'e Math\'ematique de France.

\medskip
PPP: If you think at what you saw since you started doing research, what institutional evolution of mathematics do you judge positive? 

\medskip
BT: Something is definitely positive, the creation of mathematical laboratories. In fact, the Centre de Math\'ematiques de l'\'Ecole Polytechnique, created in 1966, was the first official mathematical laboratory in France.  A few years later it was imitated at the \'Ecole Normale, by Douady and Verdier. Before that, graduate students had no office, no desk, they had to work at home. Even professors had to share offices. Little by little this infrastructure was put in place everywhere. Probably it started to develop really after 1980, because the number af mathematicians and of mathematics students grew. 
One should say that a lot of development of mathematics started in late 1960s and early 1970s because of the bomb. The government had put an incredible amount of energy in building an atomic bomb. Part of this energy was to create positions for young researchers like me. Politicians had this view, common with the United States, that you give money in the hope that some of it will be useful after a while for what you want. It worked very well in the United States and it worked relatively well in France. A whole generation of mathematicians was recruited at that time and did mathematics. For the majority of them, their research had nothing to do with the bomb. The government did probably not realise immediately that if you create a pool of young researchers, then you have to give them the means to work: offices, chairs, tables, computers... That has definitely improved. 

\medskip
PPP: What evolved in a negative way? 

\medskip
BT: First and foremost the salaries of the young researchers. That is a catastrophy, not only in mathematics of course. Also the number of positions of researchers. 

\medskip
PPP: Do you think that politicians understood the need for pure mathematics? 

\medskip
BT: What I think happened is this. During the war, politicians understood that basic science wins wars, especially because of the atomic bomb. The atomic bomb is an incredible succession of events  starting from very very pure research, and using a lot of  pure Mathematics and Physics in the development of the Manhattan project.  There are many other examples, with radar, sonar, the treatment of signals, image analysis, cryptography, and so on. Then came the cold war. So, politicians  injected a lot of money in pure science, almost indiscriminately. Just with the idea that it was a strategic investment in the context of the cold war. Basic research prospered on that perception of the politicians until the fall of the Berlin wall. After that, it did not seem so necessary to invest in basic science, because politicians had this notion, which is wrong of course, that the cold war was won, therefore they didn't need to invest so much in basic research. From that time on, the funding decreased, with a respite in France in the early 1990's because Mitterand had good ministers and advisers and now it has begun to decrease dramatically, in most western countries. That is my view. I may be wrong, but I think it is basically correct.  

Now, what politicians don't realise is that we have another enemy, the climate change, for which we are going to need a lot of basic science, because there are many things about climate and its change which we do not understand. At the same time, politicians refuse to fund basic research. They don't say they refuse, of course, they just let it die. For me, this is the main evolution concerning scientific research. Of course, mathematics is concerned. We have many politicians in France who do not have the slightest idea of what mathematical research is, or even a basic understanding of Mathematics and its applications. In France, when De Gaulle decided to make the bomb an enormous amount of resources were directed to basic science because De Gaulle's advisors knew that it was indispensable. I think my job at the CNRS is a consequence of that. You know, De Gaulle had as one of his advisors Pierre Lelong, who was a very good mathematician, and I think he listened to his advice. I wonder how many of the recent presidents or prime ministers listened to the advice of excellent basic or applied scientists. 

\medskip
PPP: What improvements would you suggest? 

\medskip
BT: I don't think one can improve without changing the perception of the political class. Politicians are currently all excited about Artificial Intelligence, but how many realise that the idea of training neural networks was first explored and developed in the 1970's by people who were mostly mathematicians and theoretical computer scientists, as well as some biologists doing fundamental research? Nowadays the time these people spent pursuing these radically new ideas in addition to their ongoing research would be spent writing grant applications. 

\medskip
PPP: You think that the main problem is one of communication? 

\medskip
BT: I don't like to call it communication, because this term has connotations which I do not like. 

\medskip
PPP: It is a profession...

\medskip
BT: Unfortunately yes. And it is not a profession which is able to help our cause. I proposed when I was director of the Institute in Chevaleret, to take in our lab some students of Sciences Po, so that they can see how a lab works. I had the idea that later, when they are politicians, they will somehow remember that people who do basic research are not there because there is light and heating, as Sarkozy said in one of his famous speeches. I think there is this deep lack of respect of the political class towards pure scientists, because they do not understand what they do and are unable to imagine the applications because they have not studied history properly. I doubt very much that among the courses in the curriculum of the \'Ecole Nationale d'Administration, of which many of our politicians are alumni, there is one which extols the historic role of basic science, whether in sociology,  industry, or for war. This reminds me of the book {\em The rape of the masses: the psychology of totalitarian political propaganda}, published by Labour Book Service at London in 1940. It appeared in French in 1939 with the title {\em Le viol des foules par la propagande politique}. Its author, Serge Tchakotine, was a student of Pavlov and a friend of Einstein. It would be useful to make a 21th Century version of this book with the same goals.

\medskip
PPP: And what happened, did those students come? 

\medskip
BT: No, because I couldn't do it. I did not stay director for long, because the ministry did some things I didn't like. I think there is one small improvement, that the Universit\'e Paris Sorbonne has a mixed Master with Sciences Po. Basically, we have to change the viewpoint of a sufficiently large part of our political class and even of our industrial class. But in industry the viewpoints are more varied. There are people in industry who understand the need for basic research. A lot of basic research in France is done in public laboratories. Such laboratories depend on public funding, which is missing because politicians do not see clearly why they should give funding apart from the research which has obvious short term applications. Now that we are again in a state of war, with climate change, we should fund basic research much more because it is a strategic investment for this war. A country which has understood this is China and I expect that in some years our politicians will say that nobody could have seen its dominance in fighting that war coming.

\medskip
\section{Evolution of research interests} \label{sect:evresint}
\medskip

\medskip
PPP:  Let us come back to your research. How would you describe in broad lines the evolution of your research interests, since your PhD till now? Are there main periods?  

\medskip
BT: There was first this period when Pham and I wrote our paper\footnote{See Section \ref{sect:firstart}.}. This was a short period of time. 

Then I started working with Monique on Hironaka's book.  We worked on it roughly from 70 to 74, when Hironaka was no longer so much in contact with us. My own interest went somehow astray. He asked us to give what we had written to the young Spanish mathematicians Jos\'e Manuel Aroca and Jos\'e Luis Vicente, to finish it. It was always important for Hironaka to have help for the redaction. I think he suffered very much when he wrote his original paper on resolution, and he wanted to share the burden with others. Of course, for us it was a fantastic opportunity. 

Then, I became more and more involved in polar varieties and studying singularities from a different viewpoint than resolution. That lasted a good while, starting in Carg\`ese in 1973, until maybe the La R\'abida paper ``{\em Vari\'et\'es polaires II}''\footnote{It is the paper \cite{T 82}.} in 1982. Already at that time I was also interested in this connection between combinatorics, convex bodies and toric varieties. That was also a period when I worked hard to characterise the equality case in what is now called {\em Khovanskii-Teissier inequalities}, as I had done for the Minkowski inequalities for multiplicities. I failed, but this effort produced the paper on {\em Bonnesen-type inequalities}\footnote{It is the paper \cite{T 82bis}.}. The result I was looking for was obtained by Boucksom, Favre, and Jonsson\footnote{In the paper \cite{BFJ 09}.} and later with minimal hypotheses by Cutkosky\footnote{See the papers \cite{C 15} and \cite{C 24}.}. Then I gave up and concentrated on writing my La R\'abida paper.  After that, my children were born and I became less active in research for a while. 

\medskip
PPP: There was again a paper with L\^e about polar varieties, at the end of the 80s\footnote{It is the article \cite{LT 88}.}. 

\medskip
BT: That's right. 

\medskip
PPP: So, you kept working on polar varieties in the eighties. 

\medskip
BT: This was not a breakthrough, but it was a good thing. It was an attempt  to insert the theory of local polar varieties in the mainstream. 

\medskip
PPP: You had also some papers on stratifications, at the Polish International Congress of Maths and the paper in Publications de l'IHES, at the end of the 80s. So you had also in the 80s a strong interest in stratifications? 

\medskip
BT: I already had a strong interest in them when I worked on Whitney stratifications in the 70s. You could say that my interest on stratifications widened. The paper in the proceedings of the ICM 83  was written just before my children were born. In the late 1980's I became interested in a problem of triangulation, which came from questions of Tr\`eves, which I had already explored before from a more special viewpoint\footnote{See \cite{T 83} and \cite{T 89}.}. I realised that there were many problems which people were struggling with in real analytic geometry, which could be solved by proving a theorem on triangulations. I used this later in a paper which was written with Brasselet over a long period of time about relative integration of subanalytic forms, again motivated by a question of Tr\`eves\footnote{It is the paper \cite{BT 11} }.

\medskip
PPP: The start of my PhD with you, in 1998, seems to me to be more or less the beginning of the period which comes till now, of your interest in a mixture of toric geometry and valuation theory. Is it correct to describe it like this? 

\medskip
BT: Yes. It started perhaps when we had the meeting in Hokkaido, in 1993 I have told you about.  Hironaka was born in 1931. So, around 1990, maybe a little before, I started telling him ``we should make some international conference for your 60th birthday''. It was absolutely out of question for him, he almost became angry. He said ``I want something like Finland''. So, he had a good memory of our Finland meeting\footnote{See Section \ref{sect:firstmeetsing}.}. Of course, we did not know how to organise a similar meeting. I told you, Finland would be impossible nowadays, and it was already impossible in 1990. But Hironaka found a solution, because he had a foundation. In a way, it was exactly like Finland, except that it lasted much longer, for two month, there were many more people and it was not devoted to a single subject. But we were both times in holiday resorts.

\medskip
PPP: What was the theme of the conference? 

\medskip
BT: None. Speaking about singularities. 

\medskip
PPP: Did this give birth to some volume? 

\medskip
BT: No, just like Finland. 

\medskip
PPP: Again, it is a good memory for you. 

\medskip
BT: Yes. It was fantastic. There were many interesting lectures, many interesting people, my children loved it, Maryvonne too of course. For some reason, I was free in mind. And there I started thinking seriously about this correspondence, which maybe had been on my mind for a while. 

\medskip
PPP: Which correspondence? 

\medskip
BT: The correspondence between convex sets and toric geometry. You know, I had proved these Minkowski-type inequalities for multiplicities and I knew that for monomial ideals multiplicity was essentially the covolume of a coconvex set. I had also written something about Bonnesen-type inequalities in algebraic geometry. In Hokkaido I started thinking that maybe there is something more to it. I thought that toric geometry had more to give to singularity theory but I needed a spark. Of course, the spark was the paper\footnote{It is the paper \cite{GT 00}.} with Rebecca Goldin, which directed me towards using toric geometry for resolution of singularities . 

\medskip
PPP: You published it in 2000. You say that you wrote it at the beginning of the 90s? 

\medskip
BT: No, no, the interval is shorter. The history of this paper is interesting. I was supposed to give a lecture at the singularities seminar. I forgot the exact date, around 1995-96.  The day before I had no idea what I should talk about. Then I came up with this idea of the simultaneous toric resolution of the monomial curve with the branch  having the same semigroup. I computed some examples, and then the next day I lectured on them. 

\medskip
PPP: But the roots are much older, because you wrote about that monomial curve in the Appendix to Zariski's book at the beginning of the 70s. It means that twenty years later you understood something new about the relation of a branch with its monomial curve. 

\medskip
BT: Absolutely. I should have thought of it before. But when I wrote the Appendix to Zariski's book, the emphasis was on moduli.  At the time we discussed a lot, in part because he was ill and I had to give some of the lectures in his place, and of course when he explained them to me I understood more. Finally he asked me to write this appendix and it is the paper I had the most pleasure in writing. Everything flowed naturally from the idea that a branch is a deformation of the monomial curve with the same semigroup. Then, twenty five years later or so, as you say, for some reason I understood something new, probably because I was desperate to find something to say at the seminar.

\medskip
PPP: The paper with Goldin was dedicated to Hironaka. Why? 

\medskip
BT: Because it proposes a method of resolution which is completely different from his. 

\medskip
PPP: It was a kind of humour? 

\medskip
BT: Yes, but also an expression of gratitude for all that he has taught me, and of appreciation of his own humour.

\medskip
PPP: Then you built up a whole program to attack the same thing in any dimension. 

\medskip
BT: Yes, exactly. 

\medskip
PPP: It seems to me that this topic is what most obsessed you, compared with the other ones. Is it true? 

\medskip
BT: I would not say that. It is a more complex program, that is for sure. 

\medskip
PPP: Have you thought from the start that it would last so long?

\medskip
BT: I think that I knew from the beginning that it would take much time, that it had to be complicated. I was convinced of several things. First of all, that one should go back to local uniformization. Secondly, that the key was in the semigroup. In characteristic zero, you can do without it, but if you want a proof which works in all characteristics, even for branches, it has to use the semigroup. The idea is that the semigroup rules all and therefore toric geometry also, even if one has to generalise it to infinite embedding dimension. I tried then to understand how complicated the semigroup could be. I discussed this question with Cutkosky and we wrote two papers on it\footnote{Those are the articles \cite{CT 08} and \cite{CT 10}.}. We realised that the semigroup can be very complicated. This convinced me that my program had to be complicated, that it could not be finished in a few months. I started trying to prove things which I could prove. Finally I found this idea that I can tame the complexity of the semigroup by using some very specific finitely generated approximations. The semigroup can be as complicated as it wants, even a finitely generated semigroup can be extremely complicated, but at the end it does not matter. And that is what I understood very slowly. At first I was obsessed with how complicated a finitely generated semigroup can be, but then I realised ``no, if it is finitely generated, then I can do my program easily by invoking regular subdivisions of rational fans, even if the combinatorics are terrible, so let's go back to approximating a semigroup by finitely generated semigroups''. And I had already described, in my first paper on valuations, a model for that with infinitely singular valuations in the plane.

\medskip
\section{Details about the Carg\`ese paper} \label{sect:detres}
\medskip

\medskip
PPP: Let us speak now in more detail about each period of your research career in mathematics. We already spoke about your paper with Pham\footnote{In Section \ref{sect:firstart}.}. What followed it? 

\medskip
BT: In Harvard I started working with Monique Lejeune-Jalabert on Hironaka's book about resolution of singularities. But very quickly we concentrated on the stratification by the constancy of the Hilbert-Samuel functions of the local algebras and on the notion of integral dependence of ideals, which had been taught to us by Hironaka. And we made some progress towards understanding what Hironaka needed for his proof of resolution. That had two aspects. One is known as {\em the continuity of maximal contact}, which means that if we have a singular subvariety of some smooth ambient space and a non-singular subspace which has maximal contact with it at one point, then necessarily it contains the Hilbert-Samuel stratum of that point and it continues to have maximal contact at nearby points of this stratum. This was one important step in Hironaka's proof. One of the things we did was to develop what we called a {\em theory of installations}, which is a bit technical. It is a generalisation of what one has to do in order to define the Newton polygon. We used that theory to prove the continuity of the maximal contact. We discovered that this theory, which appeared in our joint thesis\footnote{Parts a) are common to the theses of Monique Lejeune-Jalabert and Bernard Teissier. Bernard's thesis is \cite{T 73these}.}, could also be used to study flatness. When we told this to Hironaka, he immediately saw that he could use it to prove a flattening theorem which was at the time on his mind. So, we wrote with him a paper about local flattening. Later,  Hironaka found a much better proof, based on his division theorem\footnote{See \cite[Section 4]{H 73}.}. 

At the same time as I was working on Hironaka's program for resolution of singularities of complex analytic spaces, I was fascinated by some results about versal deformations, especially by the geometry of the discriminants. I thought a lot about it. Thinking about Zariski's conjecture on constancy of multiplicity led to my Carg\`ese paper\footnote{It is the article \cite{T 73}.}. The main ingredients of that paper were  plane sections,  preoccupation with the geometry of discriminants and integral dependence. Some people say that I introduced there integral dependence into equisingularity theory. I certainly developed it, but the roots are in ideas of Hironaka. The basic reason integral dependence comes in is that proving equisingularity amounts to proving that some partial derivatives of a function are integrally dependent of other partial derivatives multiplied by coordinates\footnote{This is called {\em condition $(c)$} in \cite[Section 2.5]{T 75}}. More generally, integral dependence of ideals is a useful algebraic concept for comparing the modulus of an analytic function to that of other analytic functions. It allows you to deal algebraically with geometric limits, for example limits of tangent hyperplanes if you compare partial derivatives or \L ojasiewicz inequalities if you compare functions.

The other result of the Carg\`ese paper which I think is important, because of its legacy, is the so-called {\em principle of specialisation of integral dependence}\footnote{It appears in Section 3 of the paper, and the name is introduced in \cite[Section 5]{T 80}.}. It says that under the equimultiplicity condition, integral dependence in a generic fiber induces integral dependence on the whole space, and thus also on the special fibers. I called this a ``principle'', because I saw it as a kind of converse to Poncelet's {\em principle of conservation of number}. Poncelet's principle says that under some conditions, if you deform a configuration, the numerical characters remain constant. The number of intersection points for example. The general idea, which I find very clearly stated by various philosophers, is that if you have discrete invariants of a configuration which deforms in a continuous way, then the discrete invariants must remain constant. The opposition between continuous and discrete is, I think, an interesting philosophical question. I must say that Poncelet was very well aware that his principle was not always true.  Some people seem to think that he was mistaken. But no, if you read Poncelet, it is quite clear that he understood that there were exceptions. 

Now, what is the principle of specialisation of integral dependence? It says that some geometric condition, the fact that in some one-parameter family of divisors you do not have a geometric accident, that is, vertical components, whose images in the parameter space are points, is actually equivalent to a numerical condition. Poncelet said that unless you have some special degeneration in a family of configurations, then the associated numbers remain constant. I say that if the number, that is, the multiplicity, remains constant, then you have no geometric accident, and this implies the specialisation of integral dependence. More generally, the idea that ``no geometric accident'' in an equidimensional family of spaces means that all irreducible components of the exceptional divisor of the blowing-up of an ideal in some auxiliary space map onto the parameter space, and that this can be insured by an equimultiplicity condition, lies behind many results of the type ``constancy of some multiplicities implies equisingularity'', typically ``equimultiplicity of polar varieties implies Whitney conditions''. \par\noindent

\par\medskip
PPP: It is a converse to Poncelet's principle. 

\medskip
BT: Yes, that is why I also called it a principle.

\medskip
PPP: Can you describe a little better the evolution of the ideas which led to the Carg\`ese paper?

\medskip
BT: The key, which made me think I made real progress, came in the train, coming back from Roscoff in the spring of 1972. It is what I call the {\em idealistic Bertini theorem}\footnote{It first appeared under this name in the notes \cite{T 81} of lectures given in Nancy in 1980.}, which is a result I like very much.  Let me describe it in the simplest form. If you have a singular point on a hypersurface and you take a generic hyperplane through that singular point, then you have two different Jacobian ideals on the hyperplane. You can look at the restriction to the hyperplane of the ideal of partial derivatives of the defining function of the hypersurface, and you can also look at the ideal of partial derivatives of the restriction to the hyperplane of the defining function, which is contained in the previous one. Now, a natural extension of Sard's theorem, or Bertini's theorem, says that for a generic hyperplane, these two ideals have the same radical.\footnote{See \cite[Remark I.1.2]{T 73}.} What my theorem says is that  they have the same integral closure, which is a big difference. Technically, this is the key to the whole paper. I remember clearly that once I saw this statement, I knew I had some key, and that all I had to do was to turn the key and open the door. 

\medskip
PPP: What were the initial questions which led to the paper? 

\medskip
BT: Initially, I was thinking about Zariski's question about equimultiplicity\footnote{It is \cite[Question A]{Z 71ter}.}: {\em is it true that if you have two singularities with the same topology, then they must have the same multiplicity?} The question is still open. I looked instead at what I believed was an easier question:  {\em in a family of hypersurfaces, if the topology is constant, is then the multiplicity also constant?}

\medskip
PPP: This is what Javier Fern\'andez de Bobadilla and Tomasz Pe\l ka proved recently\footnote{In the article \cite{FP 24}.}. 

\medskip
BT: Yes, exactly. 

\medskip
PPP: I did not know it was your conjecture. 

\medskip
BT: Yes, but quite obvious starting from Zariski's conjecture. I knew that if the topology is constant, then the Milnor number is constant.  I thought that maybe it is enough to prove that if the Milnor number is constant, then the Milnor number of the generic hyperplane sections is also constant. Because then, when you go down, through sections of higher and higher codimension, you end up with the multiplicity minus one. So, my basic idea was to prove that stronger statement. I had therefore to understand the Milnor number of a generic hyperplane section. And that's where the idealistic Bertini theorem I just mentioned comes in. The first step is to understand the difference between the Jacobian ideal of a function restricted to a hyperplane section and the Jacobian ideal of the restriction of the function to the hyperplane section.  This contains everything. If you look at it, it brings in the polar curve\footnote{\label{fn:absversrel} See footnote \ref{fn:polarpopup}. The term {\em local polar variety} first appears in the paper \cite{LT 81} with L\^e. The distinction between {\em relative} and {\em absolute} polar varieties was introduced in \cite[Chapter IV, Definition 1.4]{T 82} in order to avoid confusion between the polar loci {\em \`a la} Todd and the generalisations of the polar curve studied there.}. This is how I met the polar curve for the first time. Because, when you restrict to the hyperplane, immediately you see that the ideal you get is the restriction to the hyperplane of the ideal of a generic polar curve. So everything stems from this result, basically. Then, all I had to do was to fill in the details. This brings in the mixed multiplicities, because I had to understand how the Milnor number of a generic hyperplane section could be defined algebraically. 

\medskip
PPP: In terms of a defining function of the hyperplane and the initial function.

\medskip
BT: Yes. That is why I introduced the mixed multiplicities\footnote{See \cite[Section I.2]{T 73}.}. 

\medskip
PPP: When you say that you wanted to define the Milnor number algebraically, it seems that you are there in a kind of dictionary, sometimes you think topologically, sometimes algebraically. Can you be more precise, at which moment did you understand that you have to think algebraically?

\medskip
BT: I understood that very soon. If you think of the idealistic Bertini theorem, it is an algebraic theorem which has no topological counterpart. When I saw that, in order to understand the Milnor number of a hyperplane section, I had to understand the difference between the two Jacobian ideals, and this led to the fact that they have the same integral closure. There is a priori no topological counterpart to this. You can interpret it geometrically, that is true, but it is not a topological fact.  If you take the holomorphic function whose zero level is your hypersurface, you cut by the generic hyperplane, you look at the tangent hyperplanes to the level hypersurfaces at the intersection points when they tend to the singular point, then the theorem says that the limits of those tangent hyperplanes is always distinct from the hyperplane by which we have cut\footnote{See \cite[Paragraph I.2.8]{T 73}.}. So, the tangent hyperplanes to Milnor fibers at the points of the section do not flatten towards the cutting hyperplane as you approach the singular point. I call this kind of result {\em refined transversality theorems} because the space to which you have to be transversal moves with the points where you look. I told this once to Deligne, and his answer was ``of course, it is the contrary which would be surprising''. Which is true, but does not provide a proof.

\medskip
\section{Equisingularity problems} \label{sect:equisingpbms}
\medskip

\medskip
PPP: At the same period you wrote several papers about equisingularity\footnote{\label{fn:papersequis}These are the papers \cite{T 74}, \cite{T 75}, \cite{T 77bis}, \cite{T 80I}, \cite{T 80II}.}. Did all of them come from the same constellation of ideas? 

\medskip
BT: Yes, absolutely. 

\medskip
PPP: Can you describe a little the relation between these papers? 

\medskip
BT: The paper ``{\em Introduction to equisingularity problems}''\footnote{That is, the article \cite{T 75}.} is something I was asked to write. Basically, at the time I was really fascinated by the polar invariants, which came out of my work on the polar curve. As I said, I was not aware of Thom's manuscript.

\medskip
PPP: It was L\^e who told you this? 

\medskip
BT: Yes, but much later. So, I wanted to express this idea, which was new to me, that equisingularity was a matter of cutting. I mean, if you want to have strong equisingularity invariants, you had to include the invariants of sections of your singularity by linear spaces through the singular point. This is because I had proved that $\mu^*$-constant implies Whitney's conditions -- $\mu^*$ is just the sequence of Milnor numbers of generic plane sections\footnote{See \cite[Definition I.1.5]{T 73}.} -- so I wanted to popularize this idea. And, at the same time, I was fascinated by the interaction between the polar curve and the discriminant. Because my interest in the discriminant and my interest in the polar curve were in the beginning different things to me. I had introduced polar curves essentially for the study of Whitney conditions and of the idealistic Bertini theorem, and my interest in discriminants was of different nature. It was centered around my conjecture $\sum_i\mu_i=\mu$ and that of non-singularity of the $\mu$-constant stratum\footnote{This is \cite[Conjecture III.2.11, page 359]{T 73}, which was later disproved by Luengo in \cite{L 87}.}, which I had defined as an analytic subspace of the discriminant of the miniversal deformation. So, suddenly I realised that the two things are connected. For me, this was important, it was a connection I had not been aware of before. The whole paper is in fact somehow written with the objective of understanding this connection. 

\medskip
PPP: The whole paper ``{\em An introduction to equisingularity problems}''? 

\medskip
BT: Yes, and also ``{\em The hunting of invariants in the geometry of discriminants}''\footnote{This is paper \cite{T 77bis}.}. The ``{\em Introduction}'' is a kind of precursor of the ``{\em Hunting}'', in which in particular I could interpret the polar invariants in terms of the Newton polygon of a plane section of the discriminant\footnote{See \cite[Theorem of page 670]{T 77bis}.}. For plane branches this Newton polygon determines the equisingularity class of the branch, that is,  its topological type. This made me happy with this interpretation, because it seemed to me to provide a concrete example of what Thom called {\em reconstruction of the organizing center}, a somewhat elusive concept suggesting that one could learn a lot about a function from the discriminant of its miniversal unfolding.

\medskip
PPP: Therefore, both articles are developments of ideas you had while writing the Carg\`ese paper. 

\medskip
BT: Yes. 

\medskip
PPP: And you entered into equisingularity problems by working on Zariski's equimultiplicity question. 

\medskip
BT: Yes. This ended up in a failure, because my conjecture that $\mu$-constancy implies Whitney conditions is false, as proved by a counterexample of Brian\c{c}on and Speder\footnote{This counterexample was published in \cite{BS 75}. See the remarks of \cite[page 632]{T 75}.}. The positive result was the characterisation of Whitney conditions in a family of hypersurfaces with isolated singularities. Then, when I realised that you could characterise algebraically Whitney conditions, I had to try to do it in a more general context than for families of hypersurfaces.

\medskip
PPP:  Could you speak also about your other papers on equisingularity\footnote{That is, the remaining ones listed in footnote \ref{fn:papersequis}: \cite{T 74},  \cite{T 80I}, \cite{T 80II}.}?

\medskip
BT: In ``{\em D\'eformations \`a type  topologique constant I et II}''\footnote{That is, in  \cite{T 74}.} I tried to explain the relation between the constancy of the Milnor number in a family and the equimultiplicity stratum of the origin in the discriminant of the miniversal deformation of an isolated singularity of hypersurface, the $\mu$-constant stratum. It was in a seminar on moduli problems, so I took that viewpoint. Later, in ``{\em R\'esolution simultan\'ee I et II}''\footnote{That is, in \cite{T 80I} and \cite{T 80II}.}, I tried to clarify the concept of {\em simultaneous resolution} in a family of singularities because there are several possible definitions. I found that under certain hypotheses the strongest one implies Whitney conditions. In those papers, using simultaneous resolution for curves, I tried to give an algebraic proof of the fact that the constancy of $\mu$ in a family of reduced plane curves implies equisingularity. It did not quite work out. The paper  ``{\em R\'esolution simultan\'ee II}'' contains a proof of the principle of specialisation of integral dependence and also a proof of what is known as the {\em Pl\"ucker-Teissier formula} for the computation of the degree of the dual of a projective hypersurface with isolated singularities. In my mind at the time, for a family of projective plane curves, the constancy of the degrees and of the degrees of the duals was an equisingularity condition. In any case, the notions of {\em weak} and {\em strong simultaneous resolutions}, which I introduced there, have had some use, and were very nicely illustrated, in the case of surfaces, by Laufer\footnote{In \cite{L 83} and \cite{L 86}.}.

\medskip
\section{The collaboration with L\^e}  \label{sect:collable}
\medskip

\medskip
PPP You wrote several papers with L\^e. How did your collaboration  start? Perhaps from a common interest for polar curves?

\medskip
BT: As I said, L\^e had Thom's manuscript\footnote{See footnote \ref{fn:Thommanuscr}.}, which is an attempt to prove the finiteness of monodromy. Thom's idea was to prove the conjecture that the homological monodromy was of finite order by looking at the relative vanishing cycles with respect to a hyperplane. The idea was that if you can prove the finiteness of the monodromy for these relative vanishing cycles, then you can conclude by induction. 

So, how do you study the action of monodromy on the relative vanishing cycles? One of Thom's beautiful ideas is that those relative vanishing cycles are gradient cells on the Milnor fiber attached to the real part of the defining function of the hyperplane. So, it suffices to understand the monodromy of the critical points of that real part. But those points are exactly the intersection points of the Milnor fiber with the polar curve. This is what Thom had seen, so he introduced in that paper the polar curve and he started to study the monodromy using it. Where he goes wrong, is that he thought that all the polar invariants, which correspond to the vanishing heights of these gradient cells, are equal. That is false. He had some funny argument, thinking that one may exchange them in some cycle. But it doesn't work. 

So, L\^e was in possession of that paper. He developed Thom's idea of using the polar curves. He looked carefully at the monodromy of the various components of the polar curve, which Thom thought was independent of the component, and he proved the quasi-finiteness of the monodromy\footnote{The quasi-finiteness of the monodromy means that the eigenvalues of its action on homology are roots of unity.}, which is true. So, L\^e was interested in the polar curve for that reason. 

I was interested in the polar curve for different reasons, which I explained to you before. L\^e and I saw that my polar curve was Thom's polar curve. When we started working together, my theory of mixed multiplicities gave a formula for the multiplicity of the "absolute" polar variety of dimension $i$ of a hypersurface $(X,0)\subset (\C^{n+1},0)$ with isolated singularity, which is the closure of the set of non singular points of the hypersurface where the tangent space contains a given general linear space of codimension $i+1$. It is the sum $\mu^{(n+1-i)}+\mu^{(n-i)}$ of two consecutive sectional Milnor numbers, while the multiplicity of Thom's polar curve is $\mu^{(n)}$. When you make the alternating sum of the multiplicities of absolute polar varieties, everything vanishes except the end terms. Then we realised that this is the Euler obstruction, a notion we were discussing at the time. It was known, from work of Kashiwara, to be one plus or minus the Milnor number $\mu^{(n)}$ of a general hyperplane section. When we saw that, we thought immediately that there had to be a formula for the Euler obstruction in terms of the multiplicities of the polar varieties. So, we set out to explain this. We started actually from a paper by Gonz\'alez Sprinberg\footnote{It is the article \cite{GS 79}.}, who gave a formula for the Euler obstruction in terms of degrees of line bundles, which was a big progress I think.  The geometric interpretation of this formula showed us that we could find a more precise formula, using multiplicities of polar varieties.  Degrees of line bundles on some Nash modification is very abstract, not directly computable. Multiplicities of polar varieties are computable. So, we thought about this together, we introduced the notion of what I called {\em absolute local polar varieties}.  I wanted to make a distinction which was not always clear. If we go back to projective geometry, if we look at Todd's polar loci, they are the sets of points on a projective variety where the tangent space satisfies some linear conditions, Schubert-type conditions. They are subvarieties of the projective variety. Whereas in Carg\`ese my local polar varieties are in the ambient space. So, I wanted to make clear, when you say ``polar something'', which one you are talking about. I called Todd's polar loci ``absolute'' and the others, the ambient ones ``relative'', because they correspond to a map. They are the closures of the sets of points where the tangent spaces at non singular points of the fibers of the map satisfy some Schubert conditions. The relative polar curves are relative in that sense. Their points are those where the tangent space to the Milnor fiber satisfies some condition, namely, to be parallel to a given generic hyperplane. We introduced absolute polar curves because they correspond to the formula for the Euler obstruction.\par

\medskip
PPP: Which paper was it? 

\medskip
BT: This was the paper in Annals of Maths\footnote{It is the paper \cite{LT 81}.}. Let me explain the philosophical motivation of our formula for the Euler obstruction in terms of polar invariants. If you look at algebraic geometry in general, if you want numerical invariants, you look at characteristic classes, which are homology or cohomology classes. Our collective mindset at the Centre de Math\'ematiques de l'\'Ecole Polytechnique was that this was not really geometry. Geometry is looking at concrete subspaces, not at homology classes.  The idea of this paper is to give a formula for the Euler obstruction, which is an important invariant in the theory of Chern classes. We did not care very much about Chern classes, but we thought that the essence of the Chern class computation is the Euler obstruction. So, we wanted to understand geometrically the Euler obstruction. As I said, if you look at the alternating sum of multiplicities of absolute polar varieties of a singularity, you find the Euler obstruction. Then we think that we have a geometric interpretation. It's a general frame of mind: geometric interpretation versus cohomology. People tend to think that whenever you pronounce ``cohomology'' or ``homology'', this is topological. It is true in a sense, but if it is topological, it is not sufficiently geometric for us. Later, thanks to Brasselet, we discovered the beautiful construction of characteristic classes of singular varieties using vector fields created by M.H. Schwartz in the 1960s.\footnote{See the book \cite{S 00}.}

\medskip
PPP: You also worked with L\^e on limits of tangent hyperplanes. How did you get interested in that theme? 

\medskip
BT: It is for the following reason. When you want to define absolute polar varieties, the best way is to look at the {\em conormal space} of a singularity. In order to explain why, let me go back a little. Polar varieties are essentially defined as vanishing of minors of partial derivatives, because Schubert conditions are of that nature. But that is difficult to interpret algebraically. For example, you don't even have a direct proof that the polar varieties are generically reduced. I said before that when I was thinking about the Carg\`ese paper, I did not know how to prove that the polar curve is reduced. It was very frustrating.  The good way to prove that the polar varieties are reduced is to look at the conormal space,  the space of limits of tangent hyperplanes, which lives in the ambient variety times the dual projective space, and there you cut the conormal space by the inverse images of linear subspaces of the dual projective space. If you do that, you are actually writing down Schubert conditions. But, up there you can use transversality. This is one of the key ideas. You look at the conormal space, you know that it is reduced, you can stratify it using Whitney conditions, so that the fiber over the origin is a union of strata. You imagine that you cut it by a generic linear space. It is an open condition that it is transversal to the strata of the special fiber. Then this propagates to the nearby strata because of the Whitney conditions. This tells you that the hyperplane section is transversal to everything. So, it has to be at least generically reduced. Now your polar variety is just the image of that intersection in the starting variety, by a map which is almost everywhere an isomorphism. So you know that your polar variety is at least generically reduced. This works also for relative polar varieties, using the relative conormal space. Later\footnote{In \cite{T 81}.}, I made the same argument with the relative Semple-Nash modification. This is enough, because in order to compute multiplicities, you don't care if it has embedded components. Since in general hyperplane sections of a reduced space are not reduced, I should take the reduced part.

\medskip
PPP: Does it have embedded components? 

\medskip
BT: I think it is probable. In any case, the result we need is that they are generically reduced. It is a conjecture I wrote somewhere, that {\em the Whitney equisingularity type of the reduced part of polar varieties depends only on the analytic type of the singularity}. But that is not proved. 

\medskip
PPP: Is what you just explained about the conormal space contained in the Carg\`ese paper? 

\medskip
BT: No, it came later, it is at the end of Chapter II of the La R\'abida paper\footnote{That is, in the paper \cite{T 82}.}. Actually, L\^e and I started defining local polar varieties using the Semple-Nash modification and the Schubert varieties in the Grassmannian. It was Merle and Henry who suggested to us during discussions that it was much better to look at the conormal spaces, because their intersection theory is much simpler and also their behaviour under projection is better than for the Nash blowing-up as they explained in \cite{H-M 83}. L\^e and I worked out the consequences for projections of the conormal view point in \cite{LT 88}.

\medskip
PPP: At some moment you began to think about those phenomena in terms of symplectic or contact geometry and Lagrangian or Legendrian subspaces. Can you explain why? 

\medskip
BT: The space of limits of tangent hyperplanes is the fiber of the conormal. So, in order to have an understanding of the definition of polar varieties which I just explained, you need to understand the fibers of the conormal. Then I realised, I was at Harvard at that time, that the fiber of the conormal is a union of subvarieties which are projectively dual to some subvarieties of the projectivized  tangent cone.  The contact structure of the conormal space is what you need in order to see this duality. That is how it came in. L\^e had previously shown, for surfaces and with a different argument, that the limits of tangent hyperplanes to a surface consisted of the dual of the projectivized tangent cone, plus a number of linear spaces, which were dual to lines on the tangent cone. Then, what we did was to generalise this fact. At that time, people came to me and said ``Is it true that the limits of tangent hyperplanes are the tangent hyperplanes to the tangent cone?''. I said ``No, completely false in general, but L\^e and I know exactly what it is''. 

 Thinking in contact and symplectic terms gave us this duality between some irreducible components. This I had learnt from Pham and it immediately stuck in my mind. Pham wrote a book about D-modules\footnote{It is the book \cite{P 79}.}. It begins with a presentation of symplectic geometry, for instance the symplectic structure on the cotangent bundle. There it says that a Lagrangian variety which is irreducible is the conormal of its image in the base. When L\^e and I were thinking all the time about conormals, I said ``This is what we need! We have just to find many irreducible Lagrangian varieties''. And it worked. Moreover, we could show that Whitney conditions were equivalent to some auxiliary varieties being relatively Lagrangian over the small stratum.

\medskip
PPP: What about your paper with L\^e concerning vanishing cycles and Whitney conditions\footnote{It is the paper \cite{LT 88}.}?

\medskip
BT: There, we tried to put in order the whole theory of Whitney stratifications, topological triviality, what I called {\em total topological type}, which means the topological type of the singularity itself plus the topological types of general plane sections of various dimensions through the singular point. We tried to make a coherent presentation of the whole theory. For example, there was no proof of the fact that the constancy of the local total topological type  implies the Whitney conditions. We knew it had to be true, but there was no written proof. This resulted in a long and very complicated paper, I don't think many people have read it. But one of the good things of that paper is a complete proof of a formula which relates the multiplicities of local polar varieties in an arbitrary complex analytic Whitney-stratified space at a point of a small stratum with the local vanishing Euler characteristics of a big stratum along the small stratum. What we wanted to prove was that the constancy along the small stratum of the total local topological types of the adjacent strata -- which could be expressed by vanishing Euler characteristics -- was equivalent to the Whitney conditions for these adjacent strata along the small stratum. We were thinking that this was an important result in complex-analytic geometry, the proper converse to the Thom-Mather theorem of topological triviality along the strata of a Whitney stratification\footnote{See \cite[Section 11]{M 12} and \cite{C 74}.}. The paper was dedicated to proving that result in a way which is coherent with everything, including the existence of a canonical minimal Whitney stratification. It was very long because we tried to start from scratch. As soon as you want to compute multiplicities or do something like Morse theory in complex analytic geometry, you have to be very careful about which neighbourhoods you work with: you use polydisks whose boundaries intersect the strata in a nice way and so on. This brings a lot of technical burden.

\medskip
\section{The interest for polar varieties} \label{sect:intplar}
\medskip

\medskip
PPP: Let us speak more about polar curves. You said that they fell upon you while writing the Carg\`ese paper. Then, they followed you for decades. Can you explain the evolution of your thoughts about them? We could start from your paper ``{\em Vari\'et\'es polaires I}''\footnote{This is the article \cite{T 77}.}. 

\medskip
BT: That paper came out of my desire to understand the relation between Milnor numbers and polar invariants\footnote{\label{polinv}  Assume that $f(z_0, \dots, z_n)$ is a convergent power series without constant term and with isolated critical point at the origin. Denote by $X_0 \hookrightarrow \mathbb{C}^{n+1}$ the hypersurface singularity defined by it and by $\Gamma_q$ the irreducible components of a generic polar curve of $f$. The multiplicity of $\Gamma_q$ is denoted $m_q$ and its intersection number with $X_0$ is denoted $e_q + m_q$. The {\em polar invariants} of $X_0$ are the quotients $e_q/m_q$. For details, see \cite[Section 1.3]{T 77}.}. It turned out that it gave me much more. 

\medskip
PPP: What do you mean by ``understand the polar invariants''? Did you have a more precise question? 

\medskip
BT: ``{\em Vari\'et\'es polaires I}'' was also born from the fact that I wanted to understand the Milnor numbers of sections. I thought that the Thom-Sebastiani sum\footnote{The {\em Thom-Sebastiani sum} of two functions $f(x)$ and $g(y)$ is $f(x) + g(y)$. Here $x$ and $y$ denote disjoint multisets of variables. The name makes reference to the article \cite{ST 71}, in which the homological monodromy of the Thom-Sebastiani sum was expressed in terms of the homological monodromy of its components.} is a nice experimental ground, because you have distinct analytic types, with the same topology, just by replacing $f(x)+g(y)$ by $u(x)f(x)+v(y)g(y)$, where $u(x)$ and $v(y)$ are units. I thought ``{\em let me see if I can compute the Milnor number of a general hyperplane section of a Thom-Sebastiani sum in terms of the Milnor numbers of the components}''. It turned out that you could not. Polar invariants have numerators and denominators, the sum of the numerators is the Milnor number and the sum of denominators is the Milnor number of a general hyperplane section. When you try to compute the Milnor number of a generic hyperplane section of a Thom-Sebastiani sum, the numerators and denominators of the polar invariants come up. They come up in a very elementary way, so to speak. I was very happy, because I had a formula for the Milnor number of a generic hyperplane section of a Thom-Sebastiani sum, but not in terms of the Milnor numbers, instead in terms of the polar invariants. That was for me very illuminating about the role of the polar invariants. If you think of them as vanishing rates of the size of (relative) vanishing cycles this algebraic computation tells you something quite geometric of the same general nature as the Thom-Sebastiani theorem. This vision led me to the conjecture about the minimal exponents of isolated singularities of hypersurfaces which was recently proved by Dirks and Musta\c{t}\u{a}\footnote{In the article \cite{DM 23}.}.

\medskip
PPP: Therefore, the initial question was not about polar curves, it was about Milnor numbers of Thom-Sebastiani sums, and you had the surprise to see that polar invariants come into the story. 

\medskip
BT: Yes, exactly. So, I was very happy with that, and also I realised that the way they intervene could be summed up in a product operation on Newton polygons.  That is why I introduced a notation for Newton polygons\footnote{\label{fn:notnewtonpol} Each Newton polygon can be written as a Minkowski sum $\sum_i P_i$ of {\em elementary Newton polygons}, which have at most one compact edge. The notation concerns those elementary Newton polygons: $P_i =\Teiss{ \ell(P_i)}{ h(P_i) } $ means that the compact edge of $P_i$ joins the points $(\ell(P_i),0)$ and $(0, h(P_i))$.}. They have obviously a sum, their {\em Minkowski sum}. But it turns out that they also have a product\footnote{The product is distributive with respect to Minkowski sum and on elementary Newton polygons it is defined by the rule $P \star Q = \Teiss{ \ell(P) \ell(Q)}{ \min \left(\ell(P) h(Q), \ell(Q) h(P)\right)}$. For details, see \cite[Definition 2.2]{T 80}.}. Apart from the fact that there is no multiplicative unit, they form a semiring. This product is exactly what is needed in order to understand the Milnor number of the generic hyperplane section of the Thom-Sebastiani sum. So, I introduced this notation and I introduced the product of Newton polygons, and then the theorem was very nice, it said that the Jacobian Newton polygon of the Thom-Sebastiani sum was the product of the Jacobian Newton polygons of the components\footnote{\label{JNpolyg}The {\em Jacobian Newton polygon} of $f \in \C\{z_1, \dots, z_n\}$ is the Newton polygon of the discriminant germ of a morphism of the form $(f,\ell): (\C^n, 0) \to (\C^2, 0)$, where $\ell$ is a generic linear form.}. To me, it has the flavor of something which could be related to Hodge theory, to the mixed Hodge structures on the Milnor fibers. There is an article explaining this link, by Steenbrink and Zucker\footnote{It is the article \cite{SZ 87}.}. Once again, my philosophy is that polar invariants are geometry, Hodge structures are cohomology. If you can link them, very nice. If you cannot, I stay on the side of geometry, vanishing rates of gradient cells in the Milnor fiber and curvature. This line of thinking led to a proof of the fact that in a $\mu^{*}$-constant family of hypersurfaces the Jacobian Newton polygon remains constant although the number of irreducible components of the polar curves, and therefore of the discriminant curve may vary, as Pham had shown\footnote{See Example 3 of \cite{P 71}.}.\par
As it turned out, with this viewpoint, the simplest Thom-Sebastiani operation, $f(z_0,\ldots ,z_n)+w^N$ gave me a key to understand the topological sufficiency of jets, the \L ojasiewicz exponents, and to show that any isolated hypersurface singularity is a general hyperplane section of another isolated hypersurface singularity, answering a question of Thom.

\medskip
PPP: ``{\em Vari\'et\'e polaires II}''\footnote{This is the article \cite{T 82}.} is a continuation of those thoughts? 

\medskip
BT: No, ``{\em Vari\'et\'e polaires I}'' is about relative polar curves, ``{\em Vari\'et\'e polaires II}'' is about absolute polar varieties, essentially. 

\medskip
PPP: What came first for you?

\medskip
BT: I think they were always almost parallel in my mind. Just after the Carg\`ese paper, I thought that this cannot be an accident valid only for hypersurface singularities, that there has to be a characterisation of Whitney conditions which is completely general. So I set as a program, as a kind of dream, to try to get an algebraic characterisation of Whitney conditions. This led to ``{\em Vari\'et\'e polaires II}''. Actually, it uses a somewhat different set of ideas than the Carg\`ese paper or ``{\em Vari\'et\'e polaires I}''. It is more about induction. There is a key lemma there, which was partly inspired by a lemma of Brian\c con and Henry for surfaces\footnote{The lemma was published in \cite{BH 80}.}. Basically, the key lemma says that if we have a one-dimensional stratum in a stratified space, and a stratum which is adjacent to it and satisfies the Whitney conditions, then the polar curve of that stratum is empty. It's a very technical lemma about limits of tangents. This implies that if we take a generic hypersurface containing the small stratum, and we intersect with the big stratum, then the Whitney conditions are preserved. So, it's an inductive proof, but the spirit of the proof is somewhat different from the Carg\`ese proof.

\medskip
PPP: In both papers appears therefore this induction strategy, to cut by hyperplanes. Were you consciously inspired by Lefschetz's hyperplane section theorems? 

\medskip
BT: Absolutely not. Of course, I was aware of Lefschetz's theorems. The thing that really motivated me was the following idea. If we have a non-singular projective hypersurface and we take a generic hyperplane section, then this section is determined by the degree. If we take the cone over a non-singular projective variety, if we know the  topology of the variety, then we know the topology of all its plane sections. What was at the back of all this in my mind was that this is completely false for singular hypersurfaces. So, we have to understand better how different it is, where the difference lies. There are many ways to attack this. One of them is what I do. The main point of my proof is to show that to relate the polar varieties of a big stratum at the points of a small stratum with the Whitney conditions, you can make an induction by cutting by general planes containing the small stratum, which we can always assume to be linear.

\medskip
PPP: Let us pass now to a more recent paper about polar curves, which you wrote with Evelia Garc\'{\i}a Barroso, ``{\em Concentrations multi-\'echelles de courbure...}''\footnote{It is the article \cite{GT 99}.}. What was your motivation?

\medskip
BT:  It started with a result of Langevin\footnote{From the article \cite{L 79}.}. Very roughly speaking, the result said that the limit of the curvature of the Milnor fibers is the intersection number of the relative polar curve with the hypersurface. I found it very nice and I had discussed it with Langevin.  I learnt like this the exchange formula between the number of critical points and the curvature, which for me is a very beautiful result. Since I was interested in the polar curve, I wanted to understand better the geometry behind this curvature. I had asked myself before the following question for a plane curve singularity: {\em can I determine its topological type just by watching the dynamics of the degeneracy of the Milnor fiber to the singularity?} I didn't actually make examples, but I started to think about this in a different way. Basically, the fact is that one can localise the exchange formula near each component of the polar curve. I tried to make out a proof, which I obtained for branches. 

The geometric idea is quite simple: By the Merle-Smith theorem\footnote{From \cite{Sm 75} and \cite{M 77}.}, the branches of the polar curve have determined contact with the given branch, encoded by the polar invariants. This means that if the (general) direction determining the polar curve moves, its branches cannot move much, and the more they have contact with the given branch, the less they can move. Therefore their intersections with the Milnor fiber cannot move much. But if you have on the Milnor fiber points where the tangent moves a lot when the points moves a little, that means there is much curvature. Moreover, since in general there are several branches of the polar with the same contact with the curve, by the exchange formula all the curvature of the Milnor fiber concentrates in its intersections with several wedges thinning towards the origin. Then, one only has to verify that formulas confirm this and compute the rates at which the curvature tends to infinity in each wedge.

Then I realised that there should be an interesting theorem behind it, for arbitrary plane curve singularities. 
I suggested this as a subject for Evelia, to understand the behaviour of the polar curve in the reducible case. There was a lot of work. Finally Evelia succeeded in understanding very well what happens\footnote{This is the article \cite{G 00}.}. Then we realised that a Spanish mathematician, Eduardo Casas-Alvero, had results in the same direction proved by completely different methods, using classical italian infinitely near points\footnote{See the article \cite{C 90}.}. But I was very happy with the geometric result of Evelia, which related the Eggers diagram with the behaviour of the transforms of the polar curve in the resolution process of the curve and the Puiseux expansions of the branches of the polar curve. I thought that such a precise description of the polar curve in the reducible case was new and should be useful. To prove that it is useful, we set up to write the paper about concentration of curvature in the reducible case.  It turned out to be substantially more complicated than the branch case. In particular it showed that yes, in many cases the behaviour as $\lambda\to 0$ of the Lipschitz-Killing curvature of the Milnor fiber $f(x,y)=\lambda$ viewed as a surface in $\R^4$, the different rates at which it tends to infinity as $\lambda$ tends to zero in different wedges in $\C^2$, determine the equisingularity type of the limit curve singularity. In general, this idea that different rates of vanishing or going to infinity have both an algebraic and a geometric significance far beyond Newton polygon type results has been in my mind for many years, since I heard, in the early 1980's, a lecture by George Wanner in Geneva explaining an example of limit cycles of a quadratic vector field where if you drew three limit cycles on a sheet of paper, the fourth one has a radius of kilometres. It resonated very much with the polar invariants.

\medskip
PPP: I think that Evelia's thesis was the first use of Eggers' thesis\footnote{Eggers' thesis was published as \cite{E 82}.}. 

\medskip
BT:  As I realised, around the time of Carg\`ese, how important the polar curve and its branches were in order to relate the geometry of a Milnor fiber with that of a hyperplane section, I gave Merle as a subject for his {\em Th\`ese de troisi\`eme cycle} to compute the polar invariants of a plane branch and see if they were related to the Puiseux exponents. I already knew that the first Puiseux exponent was there. Merle did an excellent job and totally described the situation for branches\footnote{This is the article \cite{M 77}.}. Many years later I discovered that a similar work had been done by H. J. S. Smith\footnote{In \cite[Article 18]{Sm 75}.}. The subject which Brieskorn gave to Eggers for his thesis was to generalise to the reducible case the work done for branches by Merle. Eggers' work was also very nice, he introduced the {\em Eggers diagram} which is very useful\footnote{It turns out that, as shown in \cite{GGP 19}, Favre and Jonsson's {\em valuative tree} from \cite{FJ 04} may be also obtained as a projective limit of a variant of Eggers' diagrams, the so-called Eggers-Wall trees.} but as you said, at first it had no consequence one could think of because in the reducible case the polar invariants do not determine the equisingularity class. In her thesis, Evelia found how to correct this and to relate the Eggers diagram to other invariants like the resolution graph.

\medskip
\section{Bonnesen-type inequalities}  \label{sect:bonnesenwork}
\medskip

\medskip
PPP: In the middle of these works about polar varieties you proved your Bonnesen-type inequalities\footnote{See the article \cite{T 82bis}.}. Was there a relation between the two themes of research? 

\medskip
BT: No. 

\medskip
PPP: How did you come to the Bonnesen theme of research? 

\medskip
BT: This is somehow strange. When I introduced the Milnor numbers of plane sections, I proved a formula for the multiplicity of the product of two ideals in terms of their {\em mixed multiplicities}\footnote{See \cite[Section I.2.5]{T 73}.}. One of the possible definitions of mixed multiplicities -- you might call it a {\em Eudoxian approach}\footnote{This makes reference to the article \cite{GP 24}.} -- is that the multiplicity of $\mathfrak{m}^a \mathfrak{n}^b$ for $a$ and $b$ positive integers is a polynomial in $a$ and $b$, whose coefficients are the mixed multiplicities of the ideals $\mathfrak{m}$ and  $\mathfrak{n}$. I found this very striking, that there was such a simple formula. In fact, it had been anticipated by Bhattacharya\footnote{In the article \cite{B 57}. Its author, Phani Bhusan Bhattacharya, was an Indian mathematician who, after his PhD in 1954 at Delhi University, spent two years in Exeter with David Rees before going back to Delhi University. We are grateful to Prof. Jugal Verma for this information.}, but more in the classical framework of studying Hilbert functions and so on. 

At that time I was thinking in terms of the Minkowski inequality for functions, which was a mistake, because the good analog was the Brunn-Minkowski inequality in the theory of convex bodies. I was not aware of it at the time. Anyway, for ideals the analogous inequality would be  $e(\mathfrak{n}_1 \mathfrak{n}_2)^{1/d} \leq e(\mathfrak{n}_1)^{1/d} + e(\mathfrak{n}_2)^{1/d}$. I was very keen to know if this inequality was true. I saw that it would be true if there was some inequality between successive mixed multiplicities, namely,  if the sequence of mixed multiplicities is what is now called {\em log-convex}. 

\medskip
PPP: Before proving this Minkowski-type inequality for multiplicities, did you have ideas of applications for it? 

\medskip
BT: Absolutely not. 

\medskip
PPP: It was pure curiosity...

\medskip
BT: Yes, pure curiosity. 

\medskip
PPP: It came from your exploration of a dictionary. 

\medskip
BT: No, at that time the dictionary was not there. My first question about mixed multiplicities was pure curiosity. Could something like that be true? 

\medskip
PPP: But I don't understand why you asked this question. 

\medskip
BT: I don't understand either. I looked at the formula $e(\mathfrak{n}_1  \mathfrak{n}_2)=\sum_{i=0}^d{d\choose i}e(\mathfrak{n}_1^{[i]},\mathfrak{n}_2^{[d-i]})$ expressing the multiplicity of the product of two primary ideals in a local $d$-dimensional ring as a combination of mixed multiplicities. I had in mind my question about the inequalities $\frac{\mu^{(i+1)}}{\mu^{(i)}}\geq \frac{\mu^{(i)}}{\mu^{(i-1)}}$ between the mixed multiplicities of the maximal ideal and the jacobian ideal of a hypersurface with isolated singularity and I saw that a general version $\frac{e_{i+1}}{e_i}\geq\frac{e_i}{e_{i-1}}$ with $e_i=e(\mathfrak{n}_1^{[i]},\mathfrak{n}_2^{[d-i]})$\footnote{This notation designates the multiplicity of an ideal generated by $i$ general elements of $\mathfrak{n}_1$ and $d-i$ general elements of $\mathfrak{n}_2$.} implied $e(\mathfrak{n}_1 \mathfrak{n}_2)^{\frac{1}{d}}\leq e(\mathfrak{n}_1)^{\frac{1}{d}}+e(\mathfrak{n}_2)^{\frac{1}{d}}$. I verified it in the two-dimensional case.  And then I had this question in my mind. Actually, I wrote it down in the Carg\`ese paper. Later, in 1975 or something like that, I had an idea to prove it. That idea was to use the Hodge index theorem or rather its local form, the negative definiteness of the intersection matrix of the exceptional divisor of a resolution of a normal surface singularity.

It seems amazing that no one had thought to combine Bertini's theorem and the Hodge index theorem to prove inequalities between three consecutive numbers of a sequence. But that's apparently the case. The way I defined the mixed multiplicities, it is clear you can divide by functions, generators of the two ideals, so that when you have three consecutive mixed multiplicities, they are multiplicities of ideals on a two-dimensional space. That is the key. It is because I knew that mixed multiplicities were obtained by mixing generators of both ideals and taking the multiplicity of the resulting ideal, that I saw that if I cut  by everything else, I can isolate three consecutive numbers. You use Bertini to see that when you cut by sufficiently general elements, you get a normal surface, and then it is the Hodge index theorem. So, that gave me the proof. There is something funny here. It was in 1975 or something like that. Iarrobino and Eisenbud were in Paris. One day Iarrobino said --  we are  friends, we were talking all the time --  ``Eisenbud asked me a question about multiplicities''. Then he told me the question, but he had somehow deformed it a little, or perhaps I misunderstood what he said. I answered ``No, that cannot be. I know the question, and I know the answer.'' It was a consequence of what I had just done. 

\begin{figure}[h!] 
  \centering 
  \includegraphics[scale=0.7]{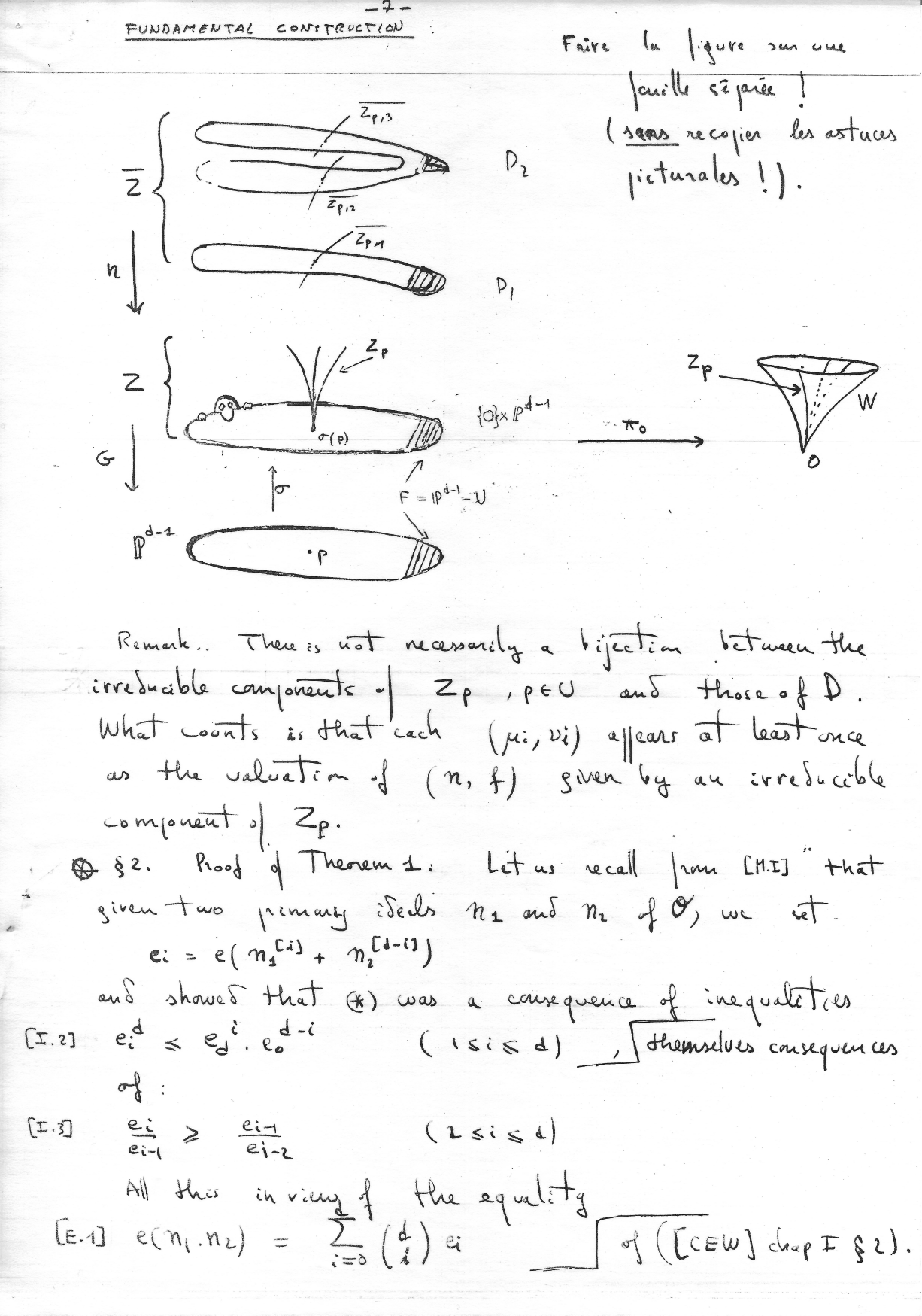} 
  \caption{Page of the manuscript of the article  \cite{T 78}}
   \label{fig:Page-manuscript}
  \end{figure} 

\medskip
PPP: That is what led to your appendix to the Eisenbud-Levine paper\footnote{It is the text \cite{T 77ter}.}? 

\medskip
BT: Yes, absolutely. Then, I asked myself about the equality case.  What happened? In 1978, after writing the appendix, I went to Nagoya invited by Matsumura. While I was there, he told me ``You know, I received a paper by Rees and Sharp, about your theorem\footnote{This was the article \cite{RS 78}.}.'' It was another proof, more general, because out of laziness, I had proved it under simplifying assumptions, like Cohen-Macaulay. Rees and Sharp gave a different proof, following Rees' beautiful ideas on multiplicities.  Then I realised that it was interesting to characterise the equality case. So, I wrote a proof of the equality case, in the Ramanujam memorial volume\footnote{This is the paper \cite{T 78}. In Figure \ref{fig:Page-manuscript} is represented a page of the manuscript of that article. Bernard explains that ``This is my geometric understanding of Rees valuations, as Monique and I rediscovered them.''}. That was the time when I realised that the analogy I had for the Minkowski inequality was wrong, that the good analogue was the Brunn-Minkowski inequality for convex bodies. I went from the Alexandrov-Fenchel-type inequalities for mixed multiplicities to Brunn-Minkowski-type inequalities. I believe that in convex geometry the same path was travelled, perhaps not in the same direction. Once Minkowski had defined mixed volumes it was perhaps a natural generalisation of the isoperimetric inequality for a plane domain, which concerns three consecutive volumes, the area of the domain, its perimeter, and the area of the unit disk.  And then I started to think that maybe such an inequality could be proved using toric geometry.

\medskip
PPP: This was your first meeting with toric geometry?

\medskip
BT: Yes and no. My first meeting was thanks to Oda, when I visited Japan. Oda was thinking about some problems linked with resolution of singularities and Newton polyhedra. He told me something about Newton polyhedra. I did not see the connection between Newton polyhedra and toric geometry, of which he was very well aware. After talking to Oda, I  thought ``this guy has some interesting ideas, I should read what he wrote''. One of the things he had written at the time was a book on toric geometry\footnote{It is the book \cite{O 78}.}. Of course, I had already a vague idea of it, because I knew Demazure's foundational paper on toric geometry\footnote{It is the paper \cite{D 70}.}. I vaguely knew that there were some varieties defined by fans. I was in a kind of haze. When I realised that the inequalities I knew were very similar to inequalities in convex geometry, I started to think that there was something there I wanted to understand. Then I understood that the linearity regions of the support function of a rational polytope defined a Demazure fan and you could attach a variety to a polytope. That is where I realised that there is a connection between toric geometry and convex bodies geometry. And then it was very easy. 

\medskip
PPP: Was it at that time that you learned the existence around Arnold of people interested by the same kind of things? And that Khovanskii had also proved what is now called the {\em Khovanskii-Teissier inequality}? 

\medskip
BT: Absolutely. Khovanskii wrote me a letter at that time, when he learnt that we had done the same thing.  He wrote that he had done this correspondence between convex bodies and toric varieties independently. He explained that there had been a talk at the Arnol'd seminar about my Minkowski-type inequalities for multiplicities.  I think that  Arnol'd was interested in this because my inequalities implied a kind of analogue of the Petrovskii-Oleinik inequalities. At the end of the talk, Khovanskii rose up and said that he could prove the same inequalities for line bundles. So he was inspired by my inequalities because he had already done the groundwork. It's very funny, we had the same idea at the same time, that made us friends. 

\medskip
PPP: It was because of this new-born interest in toric geometry that you accepted to give a Bourbaki talk on Stanley's theorem? At the beginning of that Bourbaki paper\footnote{It is the article \cite{T 80III}.}, you explained the basics of toric geometry. 

\medskip
BT: Yes. I mean, at the same time I was doing this, Stanley was using the same correspondence between polytopes and toric varieties to prove the McMullen conjecture\footnote{In the article \cite{S 80}.}. It's funny that we used the same correspondence to prove, using two aspects of Hodge theory, two different things concerning convex bodies. Again, the idea that to a lattice polytope you can associate a fan and a convex support function, and thus a projective toric variety, bloomed is several minds more or less simultaneously. At the root is Demazure's construction of toric varieties in \cite{D 70} which allows one to travel from convex bodies to algebraic geometry, the opposite of the classical direction inspired by Newton's polygon..

\medskip
\section{Valuation theory}  \label{sect:valwork}
\medskip

\medskip
PPP: Let us now pass to valuation theory. You worked a lot on the local uniformization problem\footnote{Let $X$ be a reduced germ of algebraic, analytic or formal variety and let $\nu$ be a Krull valuation of its local ring. To {\em uniformize it locally} means to find a modification of $X$ on which the center of $\nu$ is a regular point.}...

\medskip
BT: I think that my main contribution to local uniformization was to approach the problem by reduction to the case of rational valuations, instead of the classical reduction to rank one valuations, and to put forward the idea that the semigroup of a rational valuation on an excellent equicharacteristic local domain contains all the information we need to understand local uniformization. Of course,  I also made this conjecture that any  affine singularity can be resolved by a toric map after an embedding in some larger space\footnote{This conjecture appears in \cite{GT 00} and is repeated in the local uniformization version in \cite{T 03}, \cite{T 14} and \cite{T 23}.}. For me, one step towards the proof of this conjecture is the local uniformization. 

Why does {\em local uniformization} necessarily come in? If you resolve singularities, of course you uniformize valuations. So, if you start from the conviction that uniformizing a valuation is something that comes out of its semigroup, then you say OK, in order to understand resolution of singularities, I have to understand local uniformization via the semigroup. Why did this idea exist? It is because of my paper with Rebecca Goldin, which we already discussed\footnote{See Section \ref{sect:evresint}.}. It was a posteriori a big surprise for me, that for plane branches the semigroup really gives you an algorithm to resolve singularities which is blind to the characteristic. Again, I became more or less obsessed with this idea. I wanted to understand how the semigroup contained the information. As always, it's not precise at the beginning. I was working in a fog, trying to identify some solid things I can think about and on which I may try to make progress. 

 Again, my approach is to try to understand local uniformization in general instead of performing what Zariski or Abhyankar called {\em reduction to the rank one case}. In that case the value group is contained in the reals. If you look at it this way, you may have a semigroup which is simpler, because it's in the reals. But you cannot use toric geometry, you don't have any tools which might extend to positive characteristic, apart from blowing-ups of course. I wanted to concentrate on valuations for which I have tools. And my main tool is toric geometry. That means I have to study valuations which are rational -- there is no residue field extension -- but are of arbitrary rank. In that framework, my contribution is to say that the geometry of the associated graded ring is the key. Because it defines a generalised toric variety which can be resolved by a toric morphism of its ambient space independently of the characteristic, and the specialisation from the singular germ to its associated graded ring has {\em simultaneous resolution at the points picked by the valuation}. Perhaps the associated graded ring, or the corresponding toric variety, are of infinite embedding dimension, but nevertheless we get a toric variety, in the sense that we have binomial equations, perhaps infinitely many of them in infinitely many variables. So you have simple equations, binomial ones, and you can work with them. 

Of course, Abhyankar succeded to prove local uniformization for surfaces in arbitrary characteristic\footnote{In the article \cite{A 56}.}. But it's a kind of {\em tour de force}, as is the proof \cite{CP 09} of Cossart-Piltant in dimension $3$. So, I look at the graded ring associated to the valuation, at its geometry, at its presentation by binomial equations. Looking back, it seems natural. Then treat the case when the semigroup is finitely generated. That is essentially the case of Abhyankar valuations. This is very easy from this viewpoint. The way I had written the proof in my second paper on valuations\footnote{This is the paper \cite{T 14}.} is clumsy. But basically the proof is very simple. Next comes the case of non-Abhyankar valuations. This is work in progress. 

The question is natural. You have a semigroup which is not finitely generated. You want to find an embedding where a toric modification uniformizes. Well, you have to approximate your semigroup with finitely generated semigroups. How to do that? That is essentially what I achieved, except that I do not approximate my valuation by valuations with finitely generated semigroups, but by {\em semivaluations} with finitely generated semigroups. And that's it. Then you have to work with the equations, look carefully and nothing more. 

\medskip
PPP: The broad picture here is that you want to develop another general method of resolution of singularities, which uses toric geometry. 

\medskip
BT: Yes. 

\medskip
PPP: During this period of research, you made an aside with Pedro Gonz\'alez P\'erez about the Nash modification\footnote{This is the article \cite{GT 14}.}. Can you speak about your interest in Nash modifications? In some sense, it goes back to your interest in limits of tangent hyperplanes, but now you were thinking inside a toric geometric framework. 

\medskip
BT: Yes. Of course, we had this Semple-Nash conjecture in mind, that one can resolve singularities by iterating Semple-Nash modifications. That came out in part from the fact that I was unhappy with the fact that toric varieties were always normal. I said to Pedro that maybe, if we can define the non-normalised Semple-Nash modification of a toric variety, which leads us to work in the framework of non-necessarily normal toric varieties, then perhaps things become simpler. Because experimentally we know that normalisation kills a lot of information. We proved a result which is a little bit better that what Hironaka had proved using different methods, but nothing definitive.  We could not prove that the Semple-Nash modification eventually resolves\footnote{This is not always true, as was proved a little after the interview by Federico Castillo, Daniel Duarte, Maximiliano Leyton-\'Alvarez and Alvaro Liendo in the article \cite{CDLL 24}.}.  

\medskip
PPP: Do you have a strategy to get the result for any variety if you know it for toric varieties? 

\medskip
BT: Yes, going back to the specialisation to some toric variety. There are many possible specialisations to a toric variety. 

\medskip
PPP: Would it use everything else you developed for local uniformization?

\medskip
BT: Yes. Because, in order to prove that you can resolve by Semple-Nash modifications, it is enough to prove that you can uniformize valuations by Semple-Nash modifications. 

\medskip
PPP: This is a theorem? 

\medskip
BT: Yes. I am sure that Zariski knew it, and that Hironaka knows it. I mean, when we have a specific process and we want to prove that it ends with something non-singular, it suffices to show that for every valuation, eventually the center becomes non-singular. That comes from the properness of the Zariski-Riemann space over the starting variety. 

\medskip
PPP: Recently, you collaborated with Ana Bel\'en de Felipe on a problem of valuation theory\footnote{This led to the article \cite{FT 20}.}. Is this also related with your general program for local uniformization? 

\medskip
BT: It was like this. With Javier Herrera Govantes, Miguel Angel Olalla Acosta and Mark Spivakovsky, I wrote a paper\footnote{This is the paper \cite{HOST 12}.} about extending valuations to the completion. That is a very important problem on which Spivakosvky and I had thought separately. I have a very specific conjecture. Because I know how to do a lot of things for complete Noetherian equicharacteristic local rings. But if you want to prove local uniformization, you have to prove that you can reduce the excellent case to the complete case. And there appears a problem, which I identified. The conjecture is that, when you have an excellent equicharacteristic Noetherian local ring and a rational valuation, then you can extend it to a quotient of the completion in such a way that the value group is the same. Surprisingly, none seems to have thought that this is important. I had stated that in my first valuation paper\footnote{This is the paper \cite{T 03}.}, and it attracted the interest of Spivakovsky. That is why we wrote this paper together. Mark was interested in an abstract way by all possible extensions to the completion, which is a very beautiful problem in commutative algebra. This problem is still mostly open in general, but Mark solved it for valuations of rank one, already a long time ago\footnote{See \cite{HOST 14}.} and introduced new methods to attack it. I solved it for Abhyankar valuations in my second valuation paper\footnote{That is, in \cite{T 14}.}. It uses the paper that we wrote about the general case. We are still stuck with arbitrary non-Abhyankar valuations. But I think that essentially this is the only problem which remains to be solved in order to get a proof of the local uniformization. This very specific problem about completion. 

In our paper we use the fact that when we pass from a ring to its Henselization, the value group of a valuation doesn't change. We were not aware that this was a well-known result for many people in the model theory of valued fields. There existed a proof, but not a nice one in my opinion. So, I said that there should exist a natural proof. At the same time, I was looking at a paper by Lafont, who was a student of Samuel, about Henselization. I thought that there should exist a way to capture the behaviour of valuations in a Henselization using what Lafont has done. This is what Ana Bel\'en and I set up to do. 

Things turned out to be more complicated than I expected. But they were very instructive. Classically, in model theory, people were aware that if we have a field with a valuation, and you took the Henselization of the field, then the value group did not change. Of course, this doesn't really tell you much about what happens when you have a valued subring of the field, and you look at the way the Henselization changes the value semigroup. And this, in the end, is really what I am interested in. It is not so much a question about the value group of the field, but about the value semigroup of the ring. Anyway, a good first step was to look at the group, and we set out to do it using the fact that the Henselization could be constructed as the limit of very simple extensions, called Nagata extensions, defined by just one polynomial with certain properties.  The reason I liked to work on this is that if we look at what we do for this type of polynomials in Newton's algorithm of approximation of roots, we get a sequence which might not be convergent, but which is always pseudo-convergent in the sense of Ostrowski. So, we started from this beautiful idea of Ostrowski, that there is a concept of pseudo-convergence with respect to a valuation. In \cite{FT 20}, we give a new proof, which is more constructive, of the fact that if we take the Henselization of a ring, the value group does not change. The usual proof by Galois ramification is something which is a few lines, but does not tell us anything apart from the result itself. Using our construction, we reproved also a very nice theorem of Kuhlmann. When we Henselise, we know that the elements of the Henselization are approximated by elements of the smaller field, but we want to know more precisely how they are approximated. In the case of rank one it is clear and in rank greater than one, the theorem of Kuhlmann gives a very precise statement about these approximations. 

\medskip
PPP: For your program about local uniformization, you also collaborated with Dale Cutkosky. Can you explain how this collaboration enters in the program? 

\medskip
BT: In two ways. A long time ago I wanted to understand how complicated value semigroups can be. Dale had also the same kind of preoccupation, and we wrote two papers on that subject\footnote{Those are the papers \cite{CT 08} and \cite{CT 10}.}. More recently, we wrote a paper with 
Hussein Mourtada\footnote{This is the article \cite{CMT 21}.}. 
Dale's preoccupations and mine are very close, but not the same, of course. He wants to create generating sequences for valuations. If you look at the restricted world of rational valuations,  which is my world, then such sequences are formed by elements of the ring whose images in the graded ring of the valuation generate it. They are the most important elements geometrically. For example, in the case of plane branches, they are the polynomials whose roots approximate the Puiseux roots of the branch up to a fixed Puiseux exponent. My viewpoint is a bit different, as I want to start from the semigroup. But we are very close. The general problem is, {\em given a valued extension $R\subset S$ of local domains and a generating sequence for $R$, how to extend it to $S$?}

The advantage of Dale's approach is that you do not need completion. It does not depend on the problem I just mentioned, it is very straightforward. So, in the paper with Dale and Hussein we set out very special cases, but still very significant, to explain how to build generating sequences. It uses ideas coming from toric geometry, but other ideas as well. I think it is a very significant paper. Let me explain why. I should insist on this: if we look at resolutions in general, or at local uniformization, there are two philosophies. One of them is what I would call the Jung philosophy, which is to present your ring as a finite extension of a regular ring, or at least a ring whose singularities you know how to resolve. Then you try to pull up the fact that the ring downstairs is regular in order to deduce resolution for the ring upstairs. Jung actually uses embedded resolution of the discriminant. In characteristic zero it works for surfaces, except that it is rather complicated, and at the end you cannot make it work in all dimensions, for gluing up reasons.

\medskip
PPP: What are those reasons? 

\medskip
BT: When you use Jung's method in higher dimensions, you produce quasi-ordinary singularities\footnote{This is explained in \cite{L 75} and in \cite{P 11}.}. In characteristic zero, you know how to resolve quasi-ordinary singularities, but you have to glue up all those resolutions. In positive characteristic it doesn't work, because quasi-ordinary singularities can be very bad. So, that's one philosophy, using finite projections. In characteristic $p$, you look for good projections, where the algebraic extension is simple. But if you start with a given projection, in general all the trouble comes when the projection is not tame. This is where Abhyankar spent years, in order to understand this type of extension. 

Assuming that you want to uniformize a valuation, the other philosophy, which is my approach, is to use the semigroup of the valuation in order to choose an embedding which is nice,  instead of starting from a projection. If you want to relate the two methods, then look at this. Take a plane branch say $x=t^4,y=t^6+t^7$ and assume that you are in characteristic $2$. Its semigroup is generated by $4, 6$ and $15$. Then, if you take any linear projection $ax+by$ onto a line, it won't be tame. Instead, embed your curve into three-dimensional space  by $x=t^4,y=t^6+t^7, z=y^2-x^3-x^2y=-t^{15}$, where you have a $z$ axis which corresponds to $t^{15}$, and project the curve to that axis. This projection is tame\footnote{Another example is the completion $k[[x]][y]/(y^p-x^{p-1}(1+y))$ of the Artin-Schreier ring in characteristic $p$. As an extension of $k[[x]]$ it is not tame. But if I look at it as $k[[y]][x]/(x^{p-1}-y^p(1+y)^{-1})$, I have a tame extension of $k[[y]]$.}. 

\medskip
PPP: What does {\em tame} mean? 

\medskip
BT: That the characteristic does not divide the index of the small value group of the natural valuation downstairs, defined on $k[[x]]$ or $k[[y]]$, in our case $\Z$, into the big value group upstairs, which is again $\Z$. For Abhyankar valuations, my method corresponds to embedding the singularity into a larger space, where I have a finite projection onto a smooth space, which is tame. I prove that such an embedding always exists. In fact, I don't use the tame projection. 

\medskip
PPP: But then, why do you need it?

\medskip
BT: I don't need it. I say just that if I want to compare the two philosophies of resolution of singularities, I explain it like this. What I use is that certain Jacobian minors are non zero. If you look at it in terms of projections, it says that some projection is tame. But for me, the non vanishing of the jacobians corresponds to a coprimality condition of minors of a matrix of exponents which tells me that if I do such and such a toric modification of the ambiant space, the strict transform will be smooth. So, what we do in that paper is more related to the projection philosophy, because we explain how to construct generating sequences of a valued ring when we view the ring as a finite extension of a regular subring and the valuation is an extension of a valuation of the regular subring. If the extension is tame, then the result is very strong and it gives immediately the toric interpretation. Or I should say it gives a torific embedding: a {\em torific embedding} for a valuation on the ring of a singularity is an embedding of the singularity in a larger (possibly infinite-dimensional) non-singular space where it can degenerate to the associated graded ring of the valuation, which is the algebra of a generalised toric variety. Then an embedded toric resolution of a large enough part of the associated graded ring should also uniformize the valuation. The main result in my paper with Dale and Hussein is very satisfactory, it shows that the toric viewpoint is very effective. But it's not my philosophy, which is distinct from using a projection. Nevertheless, it is good to know that this kind of result exists. 

I should add that Hussein Mourtada and Bernd Schober have recently introduced\footnote{In the article  \cite{MS25}.}, for hypersurfaces in positive characteristic, a notion which they expect to play the role which quasi-ordinary singularities play in characteristic zero and that notion is relative to a projection to a non-singular space. These hypersurfaces are not quasi-ordinary but they can be ``lifted'' to quasi-ordinary singularities in characteristic zero in an equisingular way. I find it very promising.

\medskip
\section{Other mathematical works}  \label{sect:otherpapers}
\medskip

\medskip
PPP: What about your collaboration with Lipman? 

\medskip
BT :  Lipman was interested in Lipschitz saturation from a more algebraic viewpoint than could be found in my paper with Pham, and he wrote excellent papers about it. He invited me to Purdue a few times. In the mid-1970's, there appeared the paper of Brian\c con and Skoda which answered a question of Mather\footnote{This is the paper \cite{BS 74}.}: given a germ of holomorphic function on $(\C^n,0)$, can one bound the power of $f$ which is contained in its ideal $j(f)$ of partial derivatives? It was known that $f$ is integral over $j(f)$. The proof of Brian\c con and Skoda that $f^n\in j(f)$ relied on a deep analytical theorem of Skoda, while the result was clearly algebraic. I was not particularly interested in looking for an algebraic proof until one day at Harvard, Griffiths gave me a copy of his book with Harris\footnote{That is, of the book \cite{GH 78}.} and as I leafed through the chapter on residues, it struck me that there was a simple analytical proof, at least in an important special case, using the residue theorem and the interpretation of integral dependence by inequalities of moduli. This was not an algebraic proof of course but I knew one person who understood integral closure and algebraic residue theory very deeply and could transmute it into an algebraic proof, and by chance I was flying to visit Lipman in Purdue a few days later! So we wrote our paper\footnote{This is the article \cite{LiT 81}.}. It contributed to the development of a fruitful interaction between analytic residue theory and commutative algebra, as exposed in \cite{Y 05}.

\medskip
PPP: Would you like to speak also about your Comptes Rendus note concerning an arrangement of lines\footnote{This is the note \cite{T 90}.}? 

\medskip
BT: There is a mistake in that note. But it started from an interesting question. And the result is true. It's just that the way I prove it is wrong. Let me explain this. If you take any reduced complex plane curve singularity, its equisingularity type can be defined over $\mathbb{Q}$: one can find an equisingular plane curve singularity defined by an equation with rational, even integral, coefficients. This is quite striking. {\em I wondered if this was true in higher dimensions}. I know that Ignacio Luengo asked himself the same question independently around the same time. Is it true for surfaces for example? I wanted to look for a way to find a counterexample. The simplest way is to find a configuration of real lines in the plane, which cannot be defined over $\mathbb{Q}$. It is classical, I think I learnt this from the book of Gr\"unbaum\footnote{This is the book \cite{G 72}.}, that there are configurations of nine lines which are not all definable over $\mathbb{Q}$ because there exists a meeting point of four lines where the cross-ratio is not in $\mathbb{Q}$. An important fact is that this configuration is rigid, meaning that it has no deformations. The idea I had was very simple: take this arrangement of lines, complexify it, and take the cone over that. This gives you an arrangement of planes passing through the origin in $\mathbb{C}^3$. 

I claimed that this surface cannot be defined over $\mathbb{Q}$ even up to an equisingular deformation. The only thing one needs to say is that, if one has an equisingular family of surfaces like that, then one would deduce an equisingular family of tangent cones, and the rigidity implies that no member could be defined over $\mathbb{Q}$.

\medskip
PPP: What was the error? 

\medskip
BT: If you look only at the nine lines, then for Galois-theoretic reasons it can be defined over $\mathbb{Q}$. It is Jean-Benoit Bost who told me that, therefore the resulting surface singularity was definable over $\mathbb{Q}$. However, this is easy to correct by adding a well chosen tenth line to break the symmetry. In fact Adam Parusi\'nski and Lauren\c tiu P\u{a}unescu have written a nice paper\footnote{It is the paper \cite{P-P}.} correcting this and another more technical error of my note.

\medskip
\section{Interest in History of Mathematics}  \label{sect:histmath}
\medskip

\medskip
PPP: I would like to ask you now about two of your papers which are more historical, namely ``{\em Apparent contours from Monge to Todd}''\footnote{It is the article \cite{T 92}.} and  ``{\em A bouquet of bouquets for a birthday}''\footnote{It is the article \cite{T 93}.}. How do you write such a paper? Do you explore the literature, do you ask questions to other mathematicians? 

\medskip
BT: I think I used a very imperfect method. Let us talk about the paper on {\em apparent contours}. The message I wanted to transmit was that, as far as I know,  the idea of apparent contour really originated in the work of Monge  on fortifications. I thought it was interesting to follow a geometric idea from its inception in practical problems to its use in the theory of invariants of projective varieties, essentially in the complex domain. But my method is very bad, I would not recommend it. What happened is that I was interested in the history of this idea, so I read some Monge, some Poncelet, I was fascinated by the way Poncelet writes. Monge is also a very precise writer. But I must confess that I did not read the works of historians. To a large extent I followed my nose. 

I was very much impressed by two papers\footnote{The articles \cite{To 37I} and \cite{To 37II}.} of Todd. Because it is not easy at all  to find a formula which transforms equivalence classes of apparent contours or critical loci into invariants of projective geometry. Todd achieved that, he found the right formulas. I think it is a beautiful story. Of course, it is completely forgotten, now one only speaks about Chern classes. But if you think about it, everything is very natural. You can tell the same story from the viewpoint of critical points, or from the viewpoint of obstructions to extending vector fields. Because you project a projective variety onto a non-singular space, you try to lift the vector fields which are tangent to the image and then you find obstructions and it is a very natural problem to describe those obstructions. I wanted to explain the continuity thread of these ideas. But that's purely my viewpoint. As I said, it's a bad method from the viewpoint of history, because I did not look for the sources or for other papers on the same subject. Of course, I knew the papers of Chern on characteristic classes, I knew some of the work of Weil and Serre and all that. But that doesn't help. Because they are people which intend to prove results, and except perhaps in the papers of Chern, the ideas are not at all visible below the surface. While this is not at all the case in my opinion with Monge and Poncelet or with Todd. I just wanted to stress the continuity of thought. You can say it is a paper in the history of mathematics, but from a non-historian, in a very special frame of mind. It is certainly a counterexample of what a paper in the history of mathematics should be. 

\medskip
PPP: Would you say that when you read a paper in history of mathematics, this is what interests you most, to follow the development of ideas? 

\medskip
BT: Yes. In the historical development of ideas in mathematics, there are bifurcation points. What is interesting is the whole story, the continuity part and the bifurcation part, when some radically new ideas  appears, then what people do with the new ideas. In my opinion, history of mathematics should be focused in this succession of continuities and bifurcations. Of course, it is an enormous work. I believe some historians do not agree with this viewpoint, because they look for some kind of more cultural or more global view. Typically, the british school of projective geometry had this idea of invariants of projective varieties which had certain properties, that were for instance additive. And then other people coming from other schools, a purely Riemannian school, or a school of classification of fiber bundles, recast the theory. That had an enormous success. From this viewpoint I think, the key paper was that of Chern. But this is a completely different approach. Chern's paper concentrates on the differential-form viewpoint on Schubert calculus, 
whereas Todd concentrated on the geometric aspect, tangent hyperplanes, critical points of linear maps to projective spaces, and so on. Ultimately, the two got the same results. Except that in the case of Todd it was just the tangent bundle of an algebraic variety and in the case of Chern it was any fiber bundle with connections and so on. One theory fits into the other, but it is not really a good fit. There is a change of viewpoint which you have to do.  To me that was quite interesting and I tried to explain a little bit of it in that paper. 

\medskip
PPP: You mentioned several names, like Monge, Poncelet, Todd. Do you remember how you got conscious of the existence of all these actors, since the first appearance of polar curves in your work? 

\medskip
BT: I learnt at school about transformation by polars\footnote{In French it was called {\em transformation par polaires r\'eciproques}. It used a non-degenerate conic in a projective plane. With respect to it, any point has a polar line and conversely, in a way which preserves incidences between points and lines.} and that this was basically due to Poncelet. In fact, it was projective duality. So I tried to read a little of the literature of that time. I saw for instance that Bobillier wrote that Poncelet's transform was a special case of a general duality principle, and that Poncelet answered that duality was just giving a non-degenerate conic. This was very interesting, because it was conversation about conceptions, about ideas. How do you think of duality? Nowadays we don't have that type of conversation in writing, in journals. But when you read the papers of Bobillier and Poncelet, you realise that they were really discussing viewpoints. 

\medskip
PPP: Do you find part of this spirit in MathOverflow? 

\medskip
BT: Yes, I think in MathOverflow often you find people giving some type of opinion about meaning. I agree, in a way it is there. 

I was also fascinated by Poncelet's ``{\em Trait\'e des propri\'et\'es projectives des figures...}''. And in that paper I started explaining what I thought about all this. To me it was satisfying, because I put into writing something I thought interesting. 

\medskip
PPP: When did you discover the work of Todd? Because I think it is not easy to pass from Poncelet to Todd. 

\medskip
BT: The seminar which gave birth to the Lecture Notes in Mathematics 777\footnote{This is the volume \cite{DPT 80}.} came from the fact that I found Du Val's wonderful papers from the 1930s. Then I saw that there was this quadruple: Du Val, Hodge, Coxeter, Todd. They were all born between 1903 and 1908 and apparently shared a common attitude towards Mathematics. I thought that they were really great geometers. In particular, I looked at Todd's papers. I was very interested by them. They have a computational aspect: you have to find the right combinations of classes which are invariant. I was always interested in such semi-combinatorial problems. At that time I was not at all interested in Chern classes for singular varieties. I was just happy to see how Todd used the concept of {\em polar locus}. To go from Poncelet's {\em transformation par polaires r\'eciproques} to Todd's {\em polar loci} is quite a jump. You go from a non-degenerate plane conic to using the conormal space to define duality for any non-singular projective variety. Of course, in the meantime, many things had appeared in the 19th Century, with Clebsch, Noether, etc., and then the Italian geometers.

\medskip
PPP: In fact one has a version of relative polar curve in Pl\"ucker's papers. 

\medskip
BT: Yes. But that is so to speak the Poncelet tradition, except that Poncelet never published the equation of a polar curve. That is really surprising. Poncelet was apparently guided by pure thought, by what came to be called {\em synthetic geometry}. It is Pl\"ucker who published the equations. 

\medskip
PPP: Do you think that Poncelet truly worked like this? That he did not proceed like Newton, keeping the equations for himself because he wanted for a philosophical reason not to put them in the treatise? 

\medskip
BT: I have no idea. But I know that in the second edition of his treatise he showed some bitterness. He was prevented to do mathematics for a long time, because of his professional obligations. During that time, Pl\"ucker produced the formulas which bear his name using equations. Poncelet said that younger people had developed his ideas and found beautiful results, or something to that effect. 

\medskip
PPP: Would you say that what you know about history of mathematics comes more from reading than from conversations? 

\medskip
BT: Yes. 

\medskip
PPP: Let us pass now to the other paper, ``{\em A bouquet of bouquets for a birthday}'', which you dedicated to John Milnor. In it you speak about much more recent things. Do you remember how you wrote it? 

\medskip
BT: That paper also is very imperfect. I was trying to explain the heritage of Milnor's bouquet of spheres theorem\footnote{This is \cite[Theorem 6.5]{M 68}, which states that the Milnor fibers of a complex polynomial with isolated critical point at the origin has the homotopy type of a bouquet of $\mu$ spheres of middle dimension.}. For me it was an initiation to the homotopic viewpoint. I read papers by Hamm and L\^e\footnote{In particular the paper \cite{HL 91}.}, about the homotopical depth of varieties and how this corresponds in some cases to the decomposition of fibers as bouquets of spheres. I made a big mistake in that paper: I did not look at the literature from the Dutch school. For instance, at the time Siersma had also obtained a bouquet decomposition theorem in new cases\footnote{This was published in \cite{Si 95}. His results were soon generalised by Tib\u{a}r in \cite{T 96}.}. So I was trying to describe the heritage of Milnor's bouquet, but I missed a good part of it.

\medskip
\section{Interest in Philosophy of Mathematics} \label{sect:philomath}
\medskip

\medskip
PPP: Let us speak now about philosophy of mathematics. It is another topic which interests you a lot. Can you explain how you became interested in it? You wrote several papers of this nature and you coorganised at least one seminar on this subject. 

\medskip
BT: I think it all started because I need to give meaning to things. Maybe that is why I cannot really read scientific books, because very quickly the attempt to give meaning becomes too complicated for me. So I have to stop and leaf through it until I find a sentence which obviously has meaning for me. And then I start from there. That is what I did for example with Thom's book ``{\em Stabilit\'e structurelle et morphog\'en\`ese}''. In it I found many sentences with a lot of meaning. So, to come back to your question, it is basically because of that. I was interested to grasp what it means {\em to understand a proof}. 

I remember that in classes pr\'eparatoires I used to say that if we were really intelligent, then we would not have to do mathematics. It is a joke, of course. But already at that time, I was concerned by what means to understand the world, in particular through mathematics. During that period I followed lectures of Claude Levi-Strauss at the Coll\`ege de France about the isomorphisms of the structures of myths in South America and North America. I believe that unconsciously I realised that the same human brain structures trying to give meaning to the world and to find causes produce, in different environments, isomorphic systems of myths. Much later this idea came to the surface.

\medskip
PPP: So, you can say that your fundamental philosophical question is ``{\em What is meaning?}'' 

\medskip
BT: Yes, I think you can say that. Let me recall that I was attracted at first to singularity theory by the fact that there were in this field wonderful people like Zariski, Hironaka, Thom or Pham. I was very much attracted by the way they talked and by the kind of things they explained, they were at the same time very focused and very modest, they showed no desire to overwhelm anybody. I say this because when they spoke, it was easy to see that what they talked about had meaning for them. Not necessarily for me, but when they spoke about anything at all, for them it had a strong meaning. They were absolutely into it. I think that this is what attracted me, that these people knew what they were talking about in a very deep sense. So, it comes back to me trying to understand how I give meaning to things, how I understand a proof. Then I started reading a little philosophy of mathematics, Albert Lautman for example. 

\medskip
PPP: What else? 

\medskip
BT: I think that the most important one was the book called ``{\em The mathematical experience}'' by David and Hersh. Its underlying idea was the opposition between formalism and meaning. I read it many years ago and it was a big encouragement for me, that I was not the only one for which this kind of distinction is important. Then I read some Cavaill\`es, but I did not find anything in his writings which spoke to me for this question. Then Wittgenstein of course. But Wittgenstein was not at all my cup of tea. I think that his philosophy is a beautiful failure. Nevertheless,  I have a great admiration for him,  I think he was a very remarkable person and philosopher. Anyway, no philosophical text made a big impression on me. Then I started thinking by myself. At the same time I was having discussions with my friend Jean Petitot. 

Little by little I came to the conclusion that meaning is something which takes place in our unconscious mind. And so I tried to explore this idea\footnote{In the papers \cite{T 05} and \cite{T 09}.}. This is something I understand now much better. It is a philosophical problem which has been around for a long time, but which philosophy of mathematics does not really recognize. For exemple, what I write about Lautman is that he has very nice descriptions about the structure of mathematics, but he fails completely to explain why we do that, why we give meaning to this and not to that, and so on. I try to explain that his philosophy of mathematics is very coherent and satisfactory, but still it falls short for me.

I think that in a way, for what I am trying to do, a lot of philosophy is not very relevant. Because philosophy tries to describe things that are absolute, things that can be explained consciously, and I am going in a different direction. I am convinced that the important part is not conscious, that we cannot have an absolute description, and I am trying to find some interesting echoes of what goes on in the brain. And if that succeeds, I think it will already be a great success.

The impression of objectivity that we have comes from the fact that all our perceptual systems and our subconscious worlds work basically in the same way. So we have the impression that the things which we find interesting and which are also found interesting by many other people, are objective. But mathematics are objective only in the sense that a being with an isomorphic perceptual system, an isomorphic atavic experience of the world and an isomorphic brain structure would create isomorphic mathematics. If we were octopuses, we would not have the same interests or methods. So the philosophical view on mathematics as an absolute because it is logically sound is for me a serious misconception.

\medskip
PPP: Let us pass now to the seminar on cognition you coorganised. 

\medskip
BT: I coorganised it with Giuseppe Longo and Jean Petitot, at \'Ecole Normale Sup\'erieure. 

\medskip
PPP: In which period? 

\medskip
BT: In the early nineties I think. Giuseppe arrived in the computer science department, which at that time was common with the mathematics department. We discovered that we had a common interest in cognition. I knew that Jean Petitot was also very interested. The three of us had very similar motivations. 

I will speak only of mine. As I told you before, I was always interested in the understanding of mathematics. What happens when we understand a statement or a proof? I think that after doing a lot of mathematics, I can give meaning to the objects which come into play. But the meaning is something which happens in my subconscious mind. I cannot control it. I think that we have to accept this reality. And this is what I try to explain in my philosophical papers. Now, thanks to the progress of cognitive sciences, we begin to understand how this works. For example, we cannot understand how meaning is born. But we can understand, like Petitot or others did, and this was the subject of our seminar in fact, how our perceptual system works. 
 
 For instance, why do we recognize that a line is a line? I think that we have a lot of completely unconscious drives, which are adapted to our perceptual system, and which make it so that we find a line more interesting than a random curve, more meaningful if you like. To me, a visual line is an extremal state of an assembly of neurons in our perceptual system. We may call this extremal state a {\em preline} or a {\em protoline}. Then we may conceive of an object which has all the properties of a protoline according to our perceptual system, and we call that a {\em mathematical line}. It embodies not only the {\em visual line}, but also the perception we have when we walk at constant speed, what is called the {\em vestibular line}\footnote{See the article \cite{T 09}.}. 
 
 Contrary to what I did in history of mathematics, I tried to develop this idea only after a lot of discussion with competent people: Jean Petitot, Daniel Bennequin, Alain Berthoz, Giuseppe Longo, Francis Bailly. So, we were a small group of people interested in this kind of problems. Some people did much more. For instance, Petitot became a professor of cognitive science at \'Ecole Polytechnique. I continued this, I could say, as a hobby. But it is an important hobby for me. Because we can begin to grasp what it means to understand in mathematics. And also that some philosophical questions about understanding are meaningless. We tend to attribute some superhuman quality to mathematical objects, a little bit like Lautman did. I think this cannot bring us very far in our understanding. 

\medskip
PPP: Do you feel that there is a unique cerebral process for understanding, or that we have various ways of understanding?  

\medskip
BT: What startled me is that when you say that you understand a proof, it is not the outcome of a logical process. It is the sensation that you can give meaning to every step of the proof. But meaning is very hard to define. When you become a mathematician, the really hard thing is to give meaning to the objects you study.  If you cannot give meaning, then you never become a mathematician. Unfortunately, school teaching is not geared to try to explain to students the meaning of objects. The teaching concentrates in giving them the rules according to which you can manipulate objects. But for many people, if they do not have meaning, then they do not care about the rules. 

Let us come back to the example of the real line. The first thing to do is making the students understand that a number is like a point of a line, that adding two numbers is just like translation. If you don't have a picture like that in your mind, then you will not understand anything about analysis for example. And this is meaning. Because as I said, the vestibular line, which is produced when we walk at constant speed and direction, is also a special state of an assembly of neurons. There is something which I call the {\em Poincar\'e-Berthoz isomorphism}, which sort of identifies the vestibular line and the visual line. And that is, I think, the real meaning of the mathematical line. It has steps, it has a direction, like the vestibular line, but at the same time it has continuity and no orientation, like the visual line. If you identify them, then you get an orientation of the visual line. To me, if you can convince a student that the mathematical line has these two aspects, then you make enormous progress in many elementary problems of mathematics.  

Similarly, we have a protomathematical perception of parallelism. I insist on the fact that many of these things are just not conscious. Maybe this is a bit rough, but I see them as special states of assemblies of neurons. And to answer your question, I think that there is only one process of understanding, for each individual. 

\medskip
PPP: Do you think that nowadays we are able to explain better the Aha! moment of understanding, like the famous one described by Poincar\'e\footnote{In the third chapter of his 1908  book ``{\em Science et m\'ethode}''.}? He wrote that he prepared for that moment for a long time, that everything was in his mind, and that at that moment a kind of chemical reaction took place. But he cannot say more. Are we able to describe better this kind of ``chemical reaction''? 

\medskip
BT: First of all, I don't like to think about a chemical reaction, because this is too reductionist. But I think that what we do when we think about a mathematical problem,  is that we have various objects to which we give meaning. We have sufficiently worked with them to make them familiar objects. Then we ask questions.\par To ask a question about an object to which you have given a meaning is like creating a tension.  I think the Aha! moment is a moment of relaxation. The energy goes to a minimum. All this tension which was there because we did not understand, because something which had meaning for us did not fit with our organizing pulsions and our previous experiences, either as a primate or as a scientist, suddenly, for some reason which of course I cannot describe, relaxes, it goes to a minimum. Many times, it is a brutal, catastrophic event if you like. I want to emphasize that this is not a reductionist or mechanistic viewpoint. Because the physics, chemistry, electricity of what goes on in the brain is far beyond our comprehension. We can never write formulas, we can only make rough models of what goes on. The complexity of the system is huge. So we can make some kinds of analogies, like I just did, but the details escape us. Physicists like to convince us that they are quantum phenomena, inside synaptic connections for example. That may be, but I rather like to think that they are systems of huge complexity, of which we cannot hope to have any correct description. But it is like clouds or like a haze, we can make a rough outline of the shape, and get a rough understanding of what goes on inside. And that to me would already be a fantastic progress. But we are very far from that.

\medskip
\section{Main results, notions and conjectures}  \label{sect:notionsconj}
\medskip

\medskip
PPP: Which of your results are you particularly fond of? 

\medskip
BT: Interesting question... I already mentioned the {\em idealistic Bertini Theorem}\footnote{As it appears in \cite[Proposition 2.7]{T 73}. It was discussed in Section \ref{sect:detres}.}. Then perhaps the introduction of the {\em polar invariants} attached to the components of the relative polar curve of a hypersurface singularity\footnote{See footnote \ref{polinv}.}. And what I called the {\em Jacobian Newton polygon}\footnote{See footnote \ref{JNpolyg}.}. I liked very much the fact that this is constant in an equisingular deformation\footnote{That is, \cite[Th\'eor\`eme 6]{T 77}.}, because it is somehow surprising. I think it has interesting consequences. Of course there are the results about the intersection multiplicities of the relative polar curve with the hypersurface\footnote{In \cite[Proposition 1.2]{T 73} and what follows.}. Because also this led to many interesting consequences. I am not giving them in any chronological order, but as they come to my mind. Then there is this result with Rebecca Goldin, more generally the fact that the specialisation of a plane branch to its associated monomial curve is equisingular, in the sense of resolution of singularities\footnote{This is \cite[Section 6.3]{T 03}.}. That is a result I like very much, because I think it is inspiring. 

Then, there is something I thought about for a long time and then I succeeded proving, the Pl\"ucker formula for projective varieties in terms of the minimal Whitney stratification\footnote{This appears in \cite[Chapter 8]{GiT 18}.}. It illustrates the usefulness of the existence of a minimal Whitney stratification\footnote{This is explained in \cite[Chapter VI, Section 3]{T 82}.}. I like this result because I think it is very natural. It is true in the complex domain and not true in the real domain\footnote{See \cite{BT 79}.}, showing that there is still something to understand in the real domain. But the fact that there is a formula for the degree of the dual of a projective variety which depends only on the topology of the strata and the topology along the strata  is something which has been on the back of my mind for many years. In fact, I wrote a sketch of proof in a paper for the International Congress of Mathematics in 1983\footnote{This is the article \cite{T 84}.}, but I wrote down a complete proof only many years later. I did not work on it for a long time, but when I finally wrote it out, I was quite happy. It is a good example of a result in projective geometry which is a consequence of a much more general result in singularity theory applied to the case where the singularity is the vertex of the cone over a projective variety. It is also a formula which does not involve the intersection numbers of cohomology classes. It only involve topological Euler characteristics associated to general linear sections of the projective variety. The only things that move are the linear spaces. It is more geometric than cohomological formulas, more to my taste. I wonder how many other invariants of projective varieties one can describe as invariants of general singularities viewed in the special case of cones.

More recently, there is the proof of local uniformisation for Abhyankar valuations, which I think is quite neat. And even more recently, there is the general proof, if it works out, and I think it does. 

\medskip
PPP: This is your present theme of research? 

\medskip
BT: Yes. It is a relatively complicated proof, but I like it.

\medskip
PPP: This was about results, but you also defined various notions. Could you try to list the main ones, from your viewpoint? 

\medskip
BT: First of all, the definition of {\em local polar variety}, which I did with L\^e. For me, it started when I encountered the relative polar curve. Only after some reflection came out the fact that one should also define what I called {\em absolute polar varieties}. I think that is a useful definition\footnote{See footnote \ref{fn:absversrel}.}. There is also the notion of {\em graded ring of a valuation}. It is a generaliation of constructions of Samuel and Rees, which turned out to be a very significant object.

\medskip
PPP: This goes back to the appendix to Zariski's book...

\medskip
BT: In a way, yes. It was done there for the case of a branch inspired by a construction we had made with Monique at Harvard, which we called $\overline{\rm gr}\ \mathcal O$. Then I extended it to general valuations in the 1990's when I was at the \'Ecole Normale. 

\medskip
PPP: And in full generality? 

\medskip
BT: In my first paper on valuations, which was called ``{\em Valuations, deformations and toric geometry}''\footnote{It is the article \cite{T 03}.}, where I really studied that carefully. So, polar varieties, graded rings ...

\medskip
PPP: Also the sequence $\mu^*$? 

\medskip
BT: Yes, one should perhaps add it to the list. In fact, it is not really the definition of $\mu^*$, it is this idea that you should look not only at the germ of a space, but also at the same time at its generic linear sections. It is less a definition than a state of mind. 

\medskip
PPP: What about this ``{\em question of B.Teissier}'', as in the title of an article of Tevelev\footnote{It is the article \cite{Tev 14}.}? I think that it is an important idea, which is present in many of your papers. 

\medskip
BT: Yes, that is the center of most of my recent work. That is a question which interests me very strongly  and Tevelev's result is very encouraging, it is a ``proof of concept": if embedded resolutions exist at all, then toric embedded resolutions in my sense exist and they are even ``cofinal" among embedded resolutions in the sense that any embedded resolution can be realised as a diagram of strict transforms in a toric embedded resolution. 

\medskip
PPP: How did you imagine it? Were you consciously influenced by Kouchnirenko's notion of {\em Newton non-degenerate hypersurface}? 

\medskip
BT: Absolutely not. I knew of Kouchnirenko's paper, which had been presented in our Seminar\footnote{Kouchnirenko's paper \cite{K 76} was presented by Merle in the chapter \cite{M 80} of \cite{DPT 80}.}. I found the idea very beautiful because it added geometry to the Newton polyhedron, but after understanding how it worked, I lost interest. I was influenced by the result for branches proved with Rebecca Goldin. It is really the idea on which I have been working maybe for twenty years now. Of course I am aware of the fact that it means that any singularity becomes ``non-degenerate" in some sense after an appropriate re-embedding. As I said before, I knew from the beginning that it was going to be long and complicated. But it gave me a lot of satisfactions. One of the problems which was on my mind for some time was that if you have a valuation, you know you have local uniformization, but you would like to describe a sequence of birational operations which produce regular local rings which fill up the valuation ring. As soon as you experiment a little, you see that in dimensions $>2$ this is going to be very hard if you use blowing up. But if you use toric geometry, it is almost easy. It follows from a generalisation of the Jacobi-Perron algorithm. 

\medskip
PPP: Can we speak a little more about the context of introduction of each of these notions you spoke about? Let us take first the polar invariants. What was the guiding problem? 

\medskip
BT: As I told you before, I was inspired by Zariski's conjecture on multiplicity. I translated Zariski's conjecture into a simpler conjecture, which is that if you have a family of hypersurface singularities with constant Milnor number, then it has constant multiplicity. Once I knew that the Milnor number was a topological invariant, it was a very natural step to do. Then this becomes a conjecture which you can try to treat with equations or whatever. In order to do that, I thought I could prove a more general result, which is that if $\mu$ is constant, then the Milnor number of a generic hyperplane section is also constant. Then we could go down recursively and we get the constancy of the multiplicity. This turned out not to work, by the counterexample of Brian\c{c}on and Speder. But, before knowing that, my problem was to relate the Milnor number of the hypersurface singularity with the Milnor number of a generic hyperplane section. What stands between the two is exactly the intersection number of the polar curve with the Milnor fiber. That is how I introduced the polar curve, because I wanted to understand the difference between the two Milnor numbers. This gave the formula which eventually became the L\^e-Greuel-Teissier formula, which says that the sum of Milnor numbers of the hypersurface singularity and a generic hyperplane section of it is the intersection multiplicity of the generic polar curve with the hypersurface. 

Then, from there, inspired by Zariski's work on discriminants and equisingularity\footnote{As in \cite{Z 63}, say.}, it was more or less natural to ask why just look at the discriminant of a projection to a non-singular space of the same dimension? Why not look at all projections to all linear spaces of all dimensions? It is because the idea of cutting by hyperplanes worked out well, that I thought that something interesting was going to happen also by projecting to linear spaces of lower dimensions. Something interesting indeed happened, because that gave the notion of absolute local polar variety. And that was useful for Whitney conditions. 

\medskip
PPP: Both viewpoints are related, because sections are fibers of those projections. Therefore, everything got unified.

\medskip
BT: Yes, exactly. That was also a source of satisfaction, to understand that everything was in the same framework. 

\medskip
PPP: What about the minimal Whitney stratification? Was it your project to prove its existence, or did it come as a byproduct of something else? 

\medskip
BT: No, it was not my project. Once I had found a numerical characterisation of Whitney conditions, it became obvious that one has to give some kind of description of the entire Whitney stratification of a space by equimultiplicity conditions. I did that, and it turned out that this result gave the existence of a minimal Whitney stratification. 

\medskip
PPP: Let us pass now to your work about valuations. You said that your starting point was your paper with Rebecca Goldin. How did you arrive at that problem? 

\medskip
BT: I told you before that one night before giving the talk I did not know what I was going to speak about. Suddenly, I do not know why, I thought that this old appendix to Zariski's book could give me a way to resolve the singularities of a plane branch. 

\medskip
PPP: So, it came as a surprise, you suddenly had the idea. 

\medskip
BT: Yes, absolutely. Then I thought ``Let's try''. I started making a computation for a branch with two characteristic pairs. I saw that it worked out perfectly. The next day my talk was to redo part of the appendix to Zariski's book. Rebecca Goldin had come to Paris for one year. She is a very determined person. There was an agreement  between the \'Ecole Normale and Harvard that students of Harvard could come to the \'Ecole Normale for one year and conversely. This agreement had died for ten years or more. But Rebecca was determined to come to Paris, so she bothered Harvard's administration till they had to admit that the agreement had never been terminated. So, she did come to Paris. Then, she came to the Singularities seminar. She seemed interested, therefore I offered her to write this down together. She keeps doing many interesting things. She is part of a group of people called STATS who dissect the statistics of the American government and the media, to see which ones are meaningful and which ones are not, this kind of things. 

To come back to the approach to resolution of singularities, for me it presents an interesting viewpoint: we are taught that the right way to look at varieties is the  ``intrinsic" one, independent of its embedding in an affine or projective space. We talk of a very ample line bundle instead of a projective embedding. Everything has to be on the variety, essentially bundles or sheaves on the variety. This goes far, leading to derived geometry for example, or to classifying spaces, recently to topoi. In my approach, you understand the embedded resolution of a plane branch, say, only after re-embedding it, in a suitable way, in a higher dimensional space where it can be specialized in an equiresolvable way to the monomial curve with the same semigroup, whose toric resolution is blind to the characteristic. Of course I could rephrase all this in terms of a generalised Rees algebra to make it look like everything is intrinsic, but that is not what speaks to my imagination. This ``equiresolvable'' or more generally ``equisingular in some sense`` specialisation of a ring to an associated graded ring corresponding to a filtration is the phenomenon I begin to understand in general, for local uniformization.

\medskip
\section{Working methods}  \label{sect:workmeth}
\medskip

\medskip
PPP: Can you describe the process which leads to your research articles? For instance, do you wait to have all the details of the proofs before starting to write them? 

\medskip
BT: No, I don't think that is a good method. Because for example, as I told you when I wrote the first paper on the polar invariants, I had a test, which was the Thom-Sebastiani sum. Its analytic type depends on the unit you put in front of the functions whose sum you take. So, this was a very elementary test for equisingularity invariants. The analytic type varies, the equisingularity type does not and we have something we can compute with. So, I had the idea of computing the Milnor number of the generic hyperplane section of a Thom-Sebastiani sum, thinking that there has to be some formula, because everything depended only on the equisingularity types. It turned out that the formula was not in terms of the Milnor numbers of the components, but in terms of the polar invariants of the components. But when I started writing the paper, to answer the question, I had no idea what the result would be. My problem was to compute the Milnor number of a generic hyperplane section of a Thom-Sebastiani sum. OK, I start,  I make computations, I try to understand, and so on. And the result popped out. I had no idea of it before writing the paper. 

\medskip
PPP: Therefore, you make discoveries during the process of writing. In general, do you reorganise the actual process a lot to reach the final version of the paper? Do you move sections around? 

\medskip
BT: Yes. I write very much by successive approximations. 

\medskip
PPP: I wonder how much is this process universal. Have you discussed this with other mathematicians? 

\medskip
BT: Not everybody does like this. I am sure that many people write notes and then organise them once they know exactly what they want to put into the paper. In general I don't work like that, because I like to follow the thread. 

\medskip
PPP: And while you work on a problem, do you consult a lot of literature? Or do you instead isolate yourself from literature? 

\medskip
BT: Neither. I don't read lots of literature, but at the same time, if I need to clarify some notion... for instance, if you take this volume of Bourbaki...

\medskip
PPP: {\em Algebra 4 to 7}...

\medskip
BT: ... it is here because I am thinking about {\em key polynomials} and somehow I think they are related with the primitive element theorem. This theorem  is true for a field extension only when there are finitely many intermediate fields. I think that not enough attention has been given to this fact. This is a condition I want to understand better. As far as I know, this volume of Bourbaki contains the only text which treats the primitive element theorem in this way, although it may be found perhaps in work of Helmut Hasse or Emil Artin which I do not know. In most texts you see only that it is true in characteristic zero. Here it says that it is true in any characteristic, provided that there are finitely many intermediate fields. 

\medskip
PPP: Is this question related to your program about resolution of singularities? 

\medskip
BT: Yes, indirectly yes. 

\begin{figure}[h!] 
  \centering 
  \includegraphics[scale=0.7]{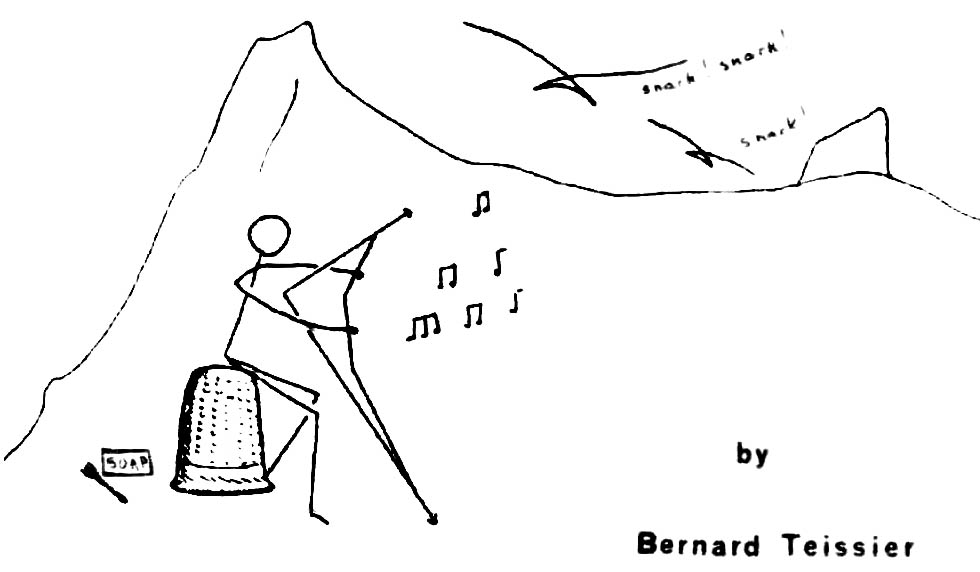} 
  \caption{The cover drawing of ``{\em The hunting of invariants in the geometry of discriminants}''}
  \label{fig:DrawingHIGD}
  \end{figure}

 \begin{figure}[h!] 
  \centering 
  \includegraphics[scale=0.5]{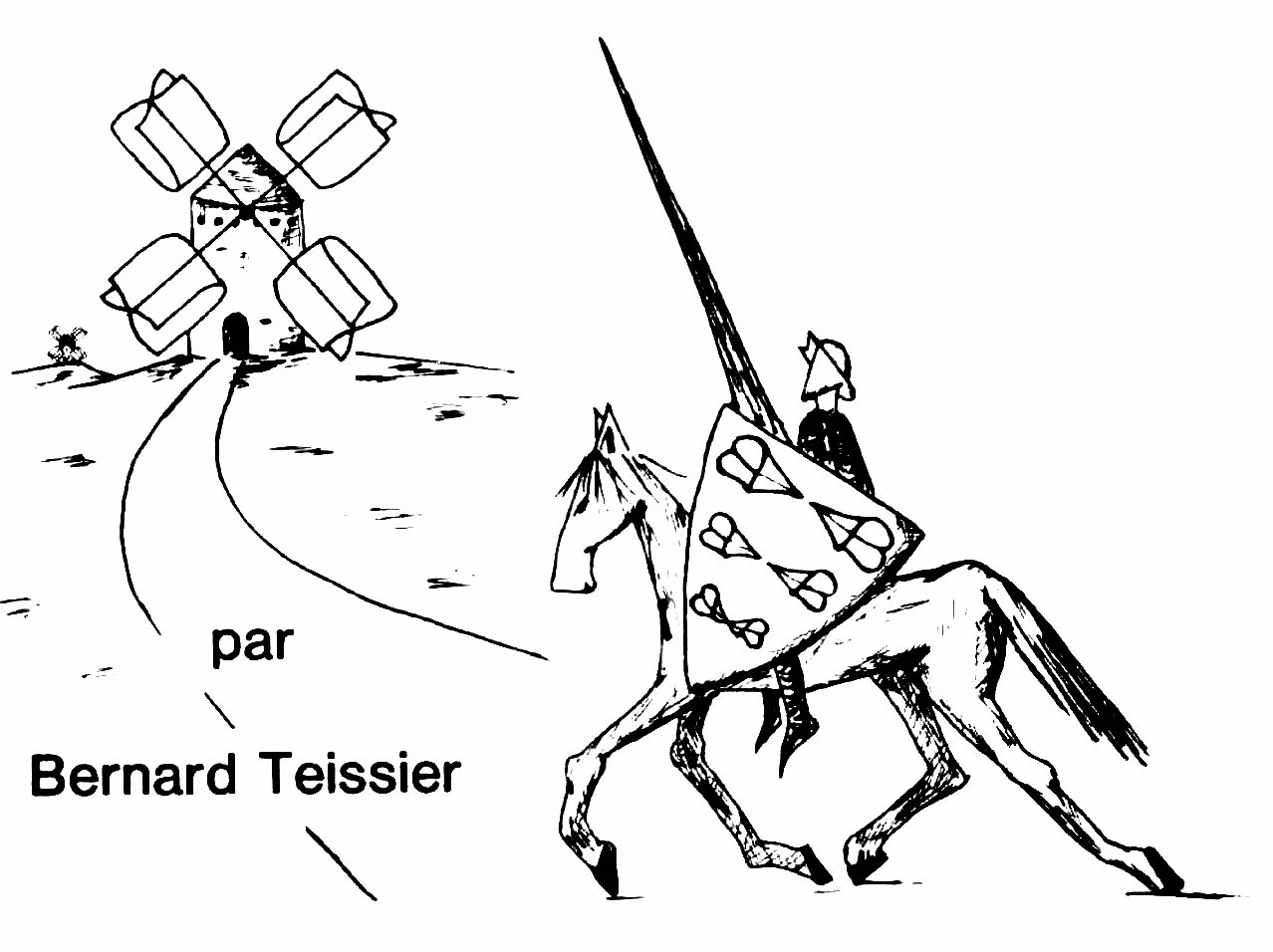} 
  \caption{The cover drawing of ``{\em Vari\'et\'es polaires II}''}
  \label{fig:DrawingVPII}
  \end{figure}

\medskip
PPP: I imagine that your writing practice changed with the appearance of LaTex? 

\medskip
BT: Oh yes. 

\medskip
PPP: Can you describe those changes? 

\medskip
BT: That is an interesting question. Before, I wrote a manuscript and I gave it to a wonderful secretary, Marie-Jo L\'ecuyer, from the Centre de Math\'ematiques, who not only could type very well, but could also detect errors.

\medskip
PPP: Impressive! Was she mathematician? 

\medskip
BT: No, not at all. She said ``Here you say this, there you say that, but it should be the same notation.'' She had a very coherent mind. She produced beautiful typescripts. At the time the Centre de Math\'ematiques produced booklets which were sent out as preprints. I liked to illustrate their covers. For example, the cover of  ``{\em The hunting of Invariants}" is a drawing I made inspired by a hike during the 1976 Oslo Conference to a mountain called Falketind\footnote{See Figure \ref{fig:DrawingHIGD}.}. And the cover of ``{\em Vari\'et\'es polaires 2}" reflects the fact that it is the notes of a course given in Madrid in 1980\footnote{See Figure \ref{fig:DrawingVPII}.}. In this case, since I could not draw a horse, I made a sketch and Marie-Jo made the final drawing. That lasted until I left for \'Ecole Normale. And just about the same time arrived LaTex, at the end of the 1980s. From then on, I type with LaTex. Before that, I had even produced a font which allowed a number of mathematical notations. 

\medskip
PPP: Have you used it in some papers? 

\medskip
BT: Yes, but only preliminary redactions which were superseded or abandoned later, like the notes of ``Stone lectures" on algebraization which I gave at Northeastern in 1988.

\medskip
PPP: One could imagine that with LaTex, people began to write much less carefully. Because you can change anything at any moment, not like with a typed text, when one needs to type again everything. Do you feel that this happened to you? 

\medskip
BT: You see, I had and still have a lot of respect for this secretary. As much as possible, I did not want to burden her with fifteen different versions of the same page. I tried to produce relatively clean manuscripts. Probably you are right: I thought more before writing, in order not to produce dirty manuscripts. But with LaTex of course, you can start typing and keep correcting as you go so the first version can be rather approximative.  For this, LaTex is very useful. And also to change names, notations. LaTex also gives you a lot of freedom with notations. I think that one has to be very careful in choosing notations. That is an important thing. 

For example, if we talk about polar invariants or the Jacobian Newton polygon, and it also answers one of your earlier questions, in this paper about polar varieties I produce a Jacobian Newton polygon for a hypersurface and a generic hyperplane section. I proved that it is constant under Whitney equisingular deformation. But then, because all generic hyperplane sections are Whitney equisingular,  you can do the same for the generic hyperplane section, you get a new polygon, and so on. When you put all these polygons together, do you get a convex polygon? That is, is there an inequality  between the largest polar invariant of the hypersurface and the smallest polar invariant of its generic hyperplane section? That was a problem which I had in my mind, and I decided I would not publish the paper before I had an answer to that. The answer is not obvious but once I had  invented a notation for describing Newton polygons I easily found a counterexample to this convexity.

\medskip
PPP: The generalised fractional notation\footnote{See footnote \ref{fn:notnewtonpol}.}? 

\medskip
BT: Yes. And because I had that notation, I could find a counterexample. I think I would have never found it otherwise. Because it is so difficult to imagine all these things at the same time. And so, I decided to publish the paper. The counterexample is towards the end\footnote{See \cite[Remark 5.3]{T 77}.}. 

\medskip
PPP: Are there other notations you introduced? 

\medskip
BT: I don't think so. Well, perhaps the notation $E_YX$ (or $B_YX$) for the result of the blowing up of $Y$ in $X$, extended to $E_YZ$ for the blowing up of the inverse image of $Y$ in $Z$ by a morphism $Z\to X$. It appears in \cite{T 82}.

\medskip
PPP: How do you prepare a talk? For instance, do you write detailed notes, or just a short plan? 

\medskip
BT: Usually, I write detailed notes, which then I forget. I mean, I write the detailed notes, I take them with me, and then I never use them. 

\medskip
PPP: This work helps you to organise your mind, to remember so well what you want to tell, that you don't need them any more. 

\medskip
BT: Exactly.

\medskip
PPP: But you could not give the talk directly without writing those notes. 

\medskip
BT: If you do not do that, then even for very simple computations you can get lost on the blackboard. Alain Connes has a very nice sentence: ``You don't prepare the lecture, you prepare the lecturer.'' 

\medskip
PPP: A talk is seen as a way to present your results. But have you used talks just to clarify your ideas? 

\medskip
BT: I am not sure I did it on purpose. But the preparation of the talk forces you to clarify your thoughts, that is inevitable. If you think about the best way to present something, then you improve your understanding. 

\medskip
PPP: Do you use talks just to see what is understandable or not, according to the reactions of people? 

\medskip
BT: No. Probably I should do that more. But I try beforehand to adjust a talk to the audience. 

\medskip
PPP: If you would give advice for giving a talk, what would be your main principles? For me, what you just said would be the most important principle: ``Think to whom you will speak.'' Would you also put it as number one?

\medskip
BT: Yes, I guess so. 

\medskip
PPP: And then? 

\medskip
BT: Then, ``Limit the number of ideas.'' 

\medskip
PPP: It is indeed very difficult to digest even three new ideas. 

\medskip
BT: Yes, that's it. You should concentrate and see what are the ideas. And that is quite difficult. I remember giving a talk in Cuernavaca about the space of orderings on $\mathbb{Z}^r$. It has a natural topology. A basis of open sets is indexed by pairs $(a,b)$ of distinct elements of $\mathbb{Z}^r$. The basic open set corresponding to such a pair is formed by the orders $\prec$ such that $a \prec b$. Then, if you look at it with this topology, Sikora proved\footnote{In the article \cite{S 04}.} that for $r\geq 2$ it is homeomorphic to the Cantor set. I found a proof of this fact based on toric geometry\footnote{See \cite{T 18}.}. But that hinged on a characterisation of the case when a projective limit of finite topological spaces is homeomorphic to the Cantor set. It is known that the Cantor set is a projective limit of finite topological spaces. Then you can ask yourself the converse. There is a simple condition. This was interesting to me. In the talk I had to teach people about toric geometry, because the finite sets turned out to be the sets of orbits in a projective system of toric varieties. So, I had a kind of quandary: do I insist on the characterisation of the projective limits which are Cantor sets or do I spend a long time speaking about toric geometry? Finally I decided to insist on toric geometry, because I thought that was the most useful thing for the audience. If they understand what a projective system of toric varieties is, then the rest would be easy. So, it is that kind of decision. 

\medskip
PPP: Let us discuss now about examples, which are very useful in talks. Those found in your papers are in general illustrations you look for after the theory is well understood, or do you put the examples which allowed you to make progress? 

\medskip
BT: I put the examples which allowed me to understand something.

\medskip
\section{Reading style}  \label{sect:readstyle}
\medskip

\medskip
PPP: You told me that while you are working on a paper, you don't read a lot. But except in that situation, how would you qualify your reading style of mathematical books and articles?  

\medskip
BT: I don't want to be presumptuous, but I remember that once I had a discussion with Thom about how he read books. He told me something which went straight to my heart: ``I don't read books. I leaf through them until I find a sentence which I find interesting.'' Then he looked at what is around the sentence and maybe at the whole book. I think that in mathematics this description also applies  to me. Of course, there is the title. But if I don't find a sentence which catches my imagination in one way or another, like this statement of the primitive element theorem, then I do not spend more time with it. 

\medskip
PPP: This is interesting, because some time ago, Norbert A'Campo took the volume in honor of his eightieth birthday\footnote{It is the book \cite{P 23}.} and said that it was a big problem for him, because each time he started reading the introduction to a paper in it, he began to dream about what he just read and he never passed the introductions... 

\medskip
BT: He is lucky, because he starts reading the introduction. I need to find some sentence somewhere which speaks to me. 

\medskip
PPP: Are there mathematical books which, after discovering some interesting sentence, you read linearly from beginning to end? 

\medskip
BT: No, never. 

\medskip
PPP: Not even Hilbert and Cohn-Vossen's ``{\em Geometry and the imagination}'' for instance? 

\medskip
BT: I did not read it linearly. I believe I read everything in it, but not linearly. 

\medskip
PPP: Even a book by Thom?

\medskip
BT: Indeed, maybe there is an exception, I think that I read ``{\em Pr\'edire n'est pas expliquer}'' linearly. 

\medskip
PPP: That book is more about philosophy of mathematics. 

\medskip
BT: Yes, indeed.

\medskip
PPP: What are your preferred mathematical books ? 

\medskip
BT: I very much like to browse in mathematical books and papers.  I do that in libraries, take a book out, open it, see if I can find something meaningful for me. It can be a book on anything. This has advantages, because I think it gives some kind of superficial culture. It has a big drawback in that you very seldom get to the deep ideas of the subject or to what made people look at this or that kind of problems. I am thinking particularly of combinatorics. I don't really understand where do many of the topics come from, because each subject has its history... 

\medskip
PPP: You browse also through old books and papers, do you go back for instance till the fifteenth century? 

\medskip
BT: No, I don't go that far back. The furthest I went with some serious reading is Monge, Poncelet and maybe Legendre too. 

\medskip
PPP:  If we speak about mathematical papers, for your research, what would be the few ones which were really important for you, in which you learnt a lot? 

\medskip
BT: Of course, Zariski's papers on many subjects: saturation, resolution, local uniformization... I learnt a lot from them. Then I would quote Laufer's book on normal two-dimensional singularities\footnote{This is the book \cite{L 72}.}. I read most of it. Also the book of Mark Kac I mentioned before. Papers of Hironaka obviously had quite an influence on me. This helped me several years ago when I wrote with L\^e an exposition of Hironaka's mathematical work\footnote{This is the article \cite{LT 08}.}. 

\medskip
PPP: Have you read completely his great paper on resolution of singularities of algebraic varieties in characteristic zero? 

\medskip
BT: In a way yes, in a way no. I did not read it line by line, but I know its  structure, I know what each part says. You can say I know it, but you cannot say I read it. 

\medskip
PPP: Do you know part of it because he told you about it? Did Hironaka like to tell mathematical stories? 

\medskip
BT: He liked to explain mathematics. He never explained how he came to do this or that. But when he explained to us the proof of resolution in the complex-analytic case, he explained what the difference was with the algebraic case. Using that, plus of course reading the paper itself, I got a coherent picture. I tried to explain that in the foreword to the book which was finally published\footnote{This is the book \cite{AHV 18}.}. Why is it interesting to read that book? Because the 1964 paper was not in final form, the proof could have been simplified. But the way to understand how it could be simpler is to read the book on the complex analytic case.

\medskip
\section{What is Singularity Theory?}  \label{sect:descrsing}
\medskip

\medskip
PPP: Passing now to a more global viewpoint, how would you describe {\em Singularity Theory} to mathematicians outside the field? What are the main interests of this theory? And in which sense can we say it is a single theory? 

\medskip
BT: Basically, it is a geometric theory. It is a theory of shapes. These shapes are constrained, because they have to be algebraic or analytic. They cannot be fractals. When people ask me what I do, I sometimes say that I do geometry outside fractal geometry. 

\medskip
PPP: Well, it is not completely true, because you said you meet Cantor sets. And you told me that you are fascinated by the fractal shapes of valuation spaces. 

\medskip
BT: Yes, that is true, but that is somehow on the boundary of the subject. I would say that if it is a theory of shapes, it is also a theory of projections, of the points where the projection of a shape does not behave like at all nearby points. And basically that's it, I think. It studies also bifurcations of shapes in families. For instance, look at any solid and at all the linear projections, in all directions. The parameter is the direction of projection and the discriminant is the set of directions where the image changes as a shape. Already there you have a fantastic amount of singularity theory going on. But there is also a lot of commutative algebra and homological algebra closely related to singularity theory which have developed in the last decades. 

\medskip
PPP: And speaking now more to people inside singularity theory, how would you describe the evolution of this theory since you entered the field? 

\medskip
BT: It has widened considerably. For example, extending the ideas of singularity theory to foliations and dynamical systems. This is an important aspect of the evolution, it did not exist fifty years ago. The geometric ideas have extended to positive characteristic much more than fifty years ago. The idea of discriminant of equations of hypersurfaces, in particular as a source of invariants existed, for plane curves, quite explicitly in the 19th century. Its use was extended to higher dimensions by Jung and then Zariski and Abhyankar but the general idea, this mixture of bifurcations \`a la Thom, plane sections, that is something that did not exist. I think that is an important development: we have many more discriminants now than we used to have. There are many more shapes to study, auxiliary shapes we could call them, than fifty years ago. And also the connection between analysis, for instance Hodge theory, and singularity theory, has very much developed. Also, and  I would put it under the same heading as analysis, the ideas connected to motivic integration. Using a linear combination of spaces as an object of the singularities world is quite a development of the last decades. The development of bilipschitz geometry gives us a key to understand much, much better than before the relationship between the actual shape of a surface singularity and the geometry of its resolution, its hyperplane sections and its polar curves.

Resolution of singularities in characteristic zero has also evolved much. There are simpler, more functorial versions, versions which are not based on blowing-up non singular centers but on {\em weighted blowing-ups}. In positive characteristic, as I mentioned, there has been the {\em tour de force} of Cossart and Piltant in dimension three.\par Why do we want to resolve singularities? We could say that it is because we want to do analysis, and we cannot do analysis on a singular space, but we can do analysis on a resolution. But the first statement is not entirely true, we can now do some analysis on singular spaces. Also we understand much better what it means to do analysis on resolutions than fifty years ago.\par While the Minimal Model Program was developing, another type of approach to birational or bimeromorphic geometry was born, motivated in part by problems coming from complex dynamical systems. It develops geometry and in particular intersection theory, on projective systems of birational maps and the Zariski-Riemann manifold, which is a space of valuations with similarities with the Berkovich analytization. Now it is really a topic, I think for instance of the works of Favre and Jonsson\footnote{In \cite{FJ 04}.}, or Gignac and Ruggiero\footnote{In \cite{GR 21}.}. In fact, thanks to singularity theory, an algebraic geometry of infinite dimensional spaces is developing, whether it is spaces of valuations, the spaces corresponding to the associated graded rings of valuations, the spaces of arcs and their algebra, the spaces whose coordinates are the coefficients of power series to geometrize approximation theorems. I have been working in part of this and I am fascinated by the whole.  
In the preface to the first volume of the ``{\em Handbook of Geometry and Topology of Singularities}''\footnote{This is the volume \cite{H 20}.} I wrote that Singularity theory sits within Mathematics much as Mathematics sits within Science. It interacts in an important way with every field of Mathematics, from Number Theory to PDE's to Combinatorics, Differential Geometry, Mathematical Physics, Algebraic Statistics and Logic, to name a few.

\medskip
\section{Extra-mathematical aspects}  \label{sect:non-math}
\medskip

\medskip
PPP: Let us discuss now some extra-mathematical aspects of your life. You met your wife Maryvonne very young. You had twins and you are happy grandparents. Was it difficult to harmonize the personal and professional lives?  

\begin{figure}[h!] 
  \centering 
  \includegraphics[scale=0.5]{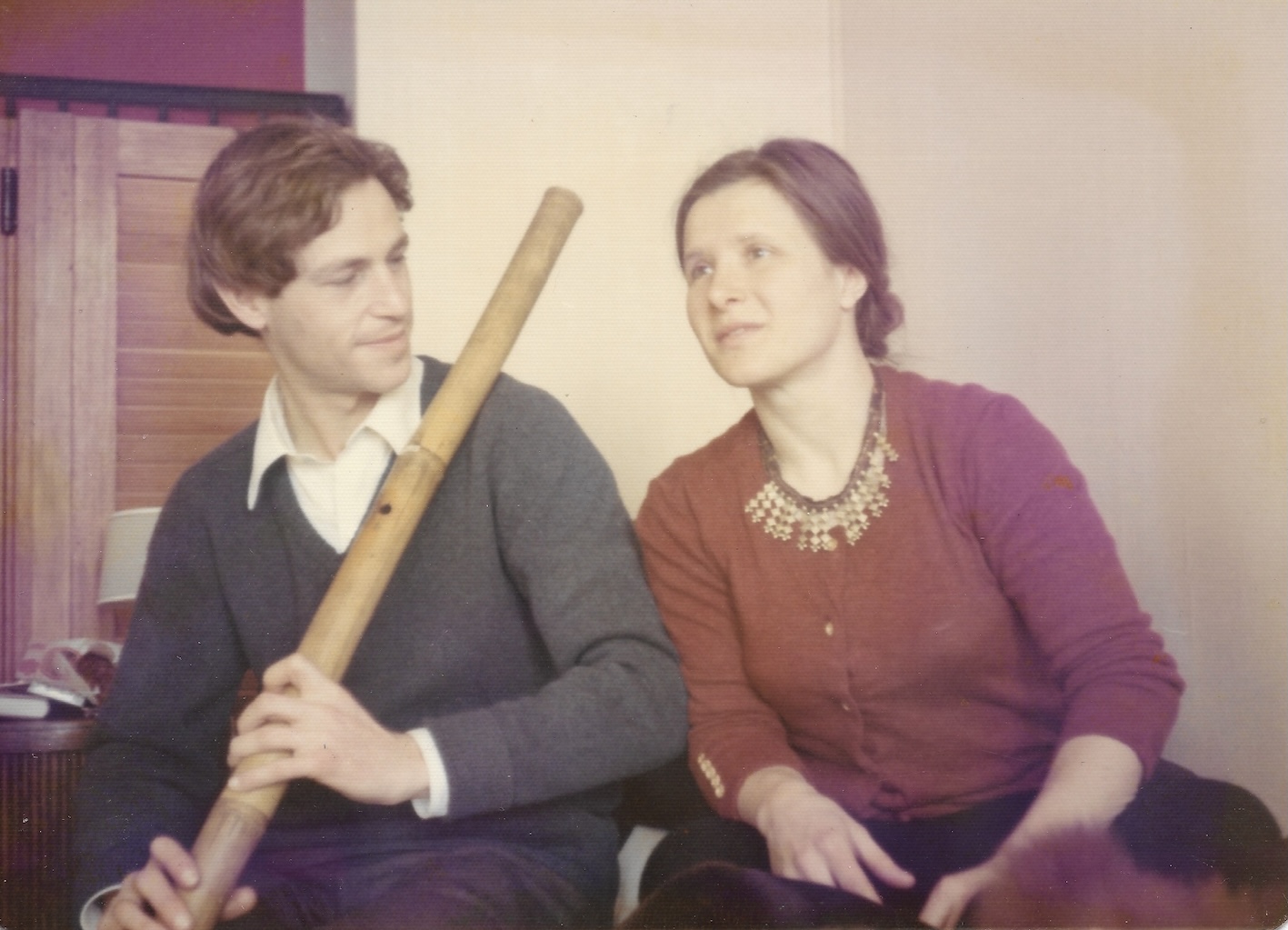} 
  \caption{Bernard and Maryvonne in January 1973}
   \label{fig:Bernard-Maryvonne73}
  \end{figure}

\begin{figure}[h!] 
  \centering 
  \includegraphics[scale=0.75]{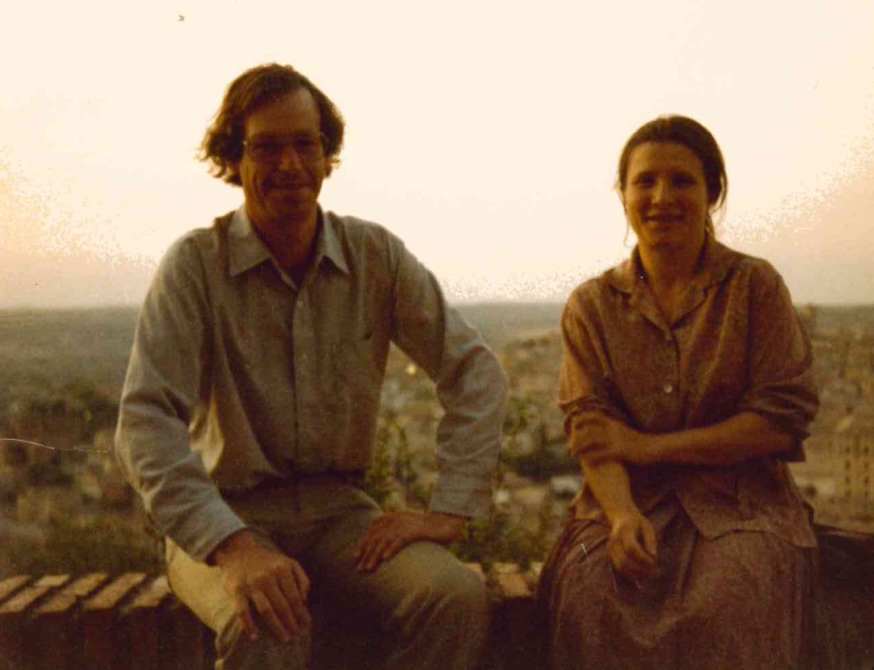} 
  \caption{Bernard and Maryvonne around 1980}
   \label{fig:Bernard-Maryvonne80}
  \end{figure}

\begin{figure}[h!] 
  \centering 
  \includegraphics[scale=0.75]{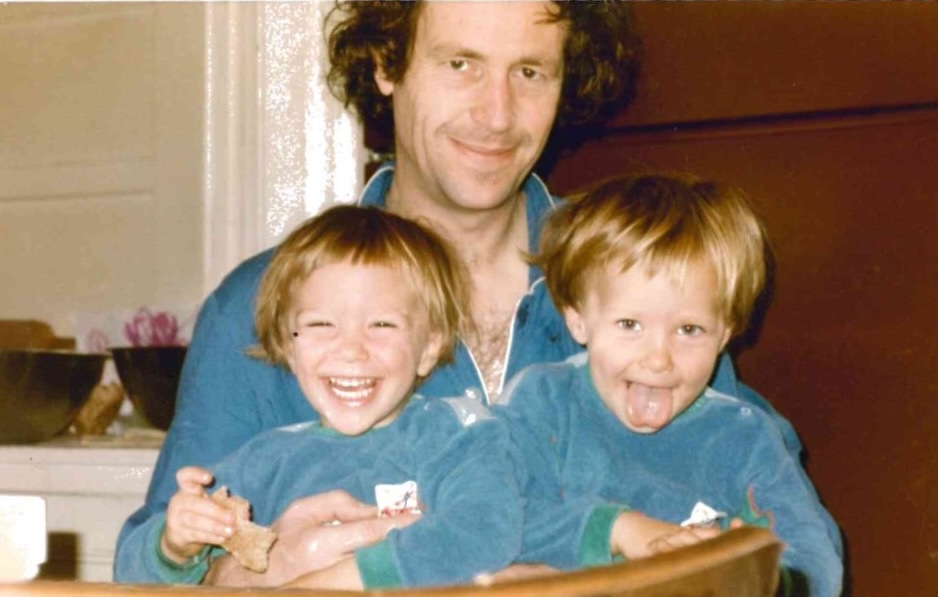} 
  \caption{Bernard with his twins, Jean and Anne, in 1986}
   \label{fig:Bernard-Anne-Jean}
  \end{figure}

\medskip
BT: I am ashamed to say this, but I think it was easy for me and difficult for Maryvonne. She made a lot of sacrifice, more than me. Because we travelled a lot, she could not fully realise her potential as a mathematician, which was quite important. I think I helped sufficiently with the children when they were small. But most of the time, most of the burden was on Maryvonne's shoulders. I am not especially happy with this, but that's how things worked out. 

\medskip
PPP: I think that your children don't have the feeling that you neglected them. Both are researchers. In other fields than mathematics, but they saw you were happy doing that and they chose to do the same kind of activity. 

\medskip
BT: Perhaps that's my best contribution to their education, in addition to summers in our mountain hut. Maryvonne contributed in so many other ways. She was the one who had the idea to enroll them in the \'Ecole des enfants du spectacle, because she had a friend, Marilyse Beffa, whose children had been there. She took care of their musical education, their violin lessons, their gymnastics lessons. These activities were necessary to be admitted in that school, which is a public school for children who have a lot of extra activities, be it acting, music, skating, playing chess, circus, gymnastics, etc. They were very happy there I think and were with kids who all had some intense activity. When they had to leave that school they went to the ``\'elite" Lyc\'ee Louis Le Grand, in different classes. In the middle of the first term we were summoned by both main teachers separately and they had the same grievance about both children: ``They laugh too much". While our children were in that Lyc\'ee, Maryvonne was very active in the parents' association. Later she directed Jean towards  the \'Ecole Sup\'erieure de Physique et Chimie Industrielles, an engineering school where he really thrived. 

\medskip
PPP: Can you tell me how you met Maryvonne? 

\medskip
BT: Oh, we were both in a course given by Michel Mend\`es-France, about a mixture of topology, combinatorics and analysis. 

\medskip
PPP: When was this? 

\medskip
BT: In the fall of 1966. 

\medskip
PPP: Where did this course take place? 

\medskip
BT: At the Institut Henri Poincar\'e. Then Michel Mend\`es-France organised a hitchhiking trip to Afghanistan. Maryvonne and I took part in it, and that's how we really met.

\medskip
PPP: I think that you went together at Harvard and that she worked on her thesis there. 

\medskip
BT: Yes, indeed. 

\medskip
PPP: Can you tell a little about her mathematical work? 

\medskip
BT: She was a student of Choquet and she worked on combinatorial topology, ultrafilters, subjects which are, I think, still quite interesting. I remember that in the eighties we met a set-theorist, Akihiro Kanamori. We started discussing mathematics and at some point Maryvonne said something and Aki said ``do you know the work of Daguenet?'' Daguenet is her maiden name. 

\medskip
PPP: Very nice! 

\medskip
BT: I regret I couldn't do more for her mathematical career. But then after a while our children took a lot of our energy. 
Until her retirement she taught Mathematics for Computer Science at Paris 7 University. She liked that because she had a varied audience, not only computer science students but also mathematicians, physicists and chemists, and also because she was free to choose what she taught. 

\medskip
PPP: Let us discuss about your extra-mathematical interests. I am thinking for instance about mountain, scientific essays and literature. I remember that you told me that at Harvard you discovered Stephen Jay Gould, and that you liked him very much. 

\medskip
BT: I think that it was in the 1980s. I had already read some of his books, which treat of evolution in a way which is both very concrete and conceptual. But there I saw that he was giving lectures in the Science Center. This was my recreation, I think it was on Wednesdays at eleven. I went down three floors and Stephen Jay Gould would be there telling stories. 

\medskip
PPP: Did you discuss with him? 

\medskip
BT: Not really. I just enjoyed listening to him telling stories to undergraduates about evolution. 

\medskip
PPP: What non-mathematical books mattered particularly for you? I remember for instance that you like to quote from time to time from Maimonides' ``{\em The guide for the perplexed}''. Why did you like it? 

\medskip
BT: As I told you before, I browse through a book till I find a sentence which interests me. I was doing exactly that with the book of Maimonides, in a library. I think that the sentence which got my attention said that time was not an attribute of God. Then I started reading more carefully. I discovered a very sharp mind at work on the problem of logic versus religion.  I enjoyed that very much. Another author in that line which I enjoyed is Fontenelle. He too had a very sharp mind, he dissected religious literature and tradition in a very interesting way. 

\medskip
PPP: What about the literature you like? The writers which come to my mind because you told me several times about them are Hermann Hesse and Umberto Eco. 

\medskip
BT: Yes, especially Hesse because of his ``{\em Glasperlenspiel}'', which I read in French and which struck me very much. The novel presents a kind of historical tapestry woven around a game, the glass beads game, in which the characters are all interested and try to excel, but which is never made explicit, you never know what the game is. For me, it is a description of what culture meant for Hesse, including the fact that it cannot be the object of a definition.  But before that, in another vein completely, I was very fond of Alexandre Dumas father. I think that the character of Athos in ``{\em Les Trois Mousquetaires}'' had a lot of influence on me because he is the most morally complex character, an example of rectitude without being righteous. I read that in my early teenage years, and it stayed with me. I must add that I also read a lot of science fiction.

\begin{figure}[h!] 
  \centering 
  \includegraphics[scale=0.35]{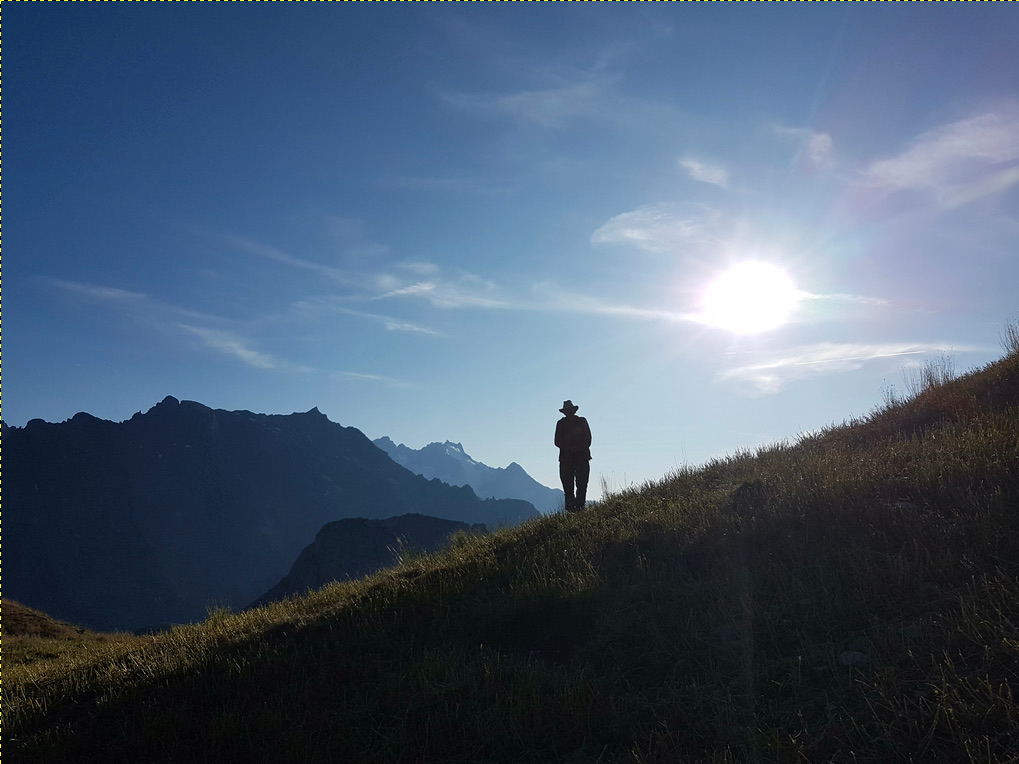} 
  \caption{Bernard in the Alps in 2019}
   \label{fig:Alpes}
  \end{figure}

\medskip
PPP: Let us speak now about your passion for mountains. 

\medskip
BT: I like being in the mountains, I like their landscapes, the freedom they bring, and also the fact that they present some danger, that you have to be constantly aware of the weather, of the time, of the terrain. 

\medskip
PPP: Do you think about mathematics while you hike? 

\medskip
BT: Sometimes. Nietzsche said -- Nietzsche also had some influence on me -- that the only thoughts that have some value are those which come when you are walking. Of course, he tended to be a bit extreme, but I think there is some truth in that. I did a lot of hiking alone, when I was a teenager or a bit later. I think that hiking puts you in a special frame of mind. I cannot say that I ever had an idea while hiking. But meditation while hiking is something I enjoy very much. And of course I like skiing. And the mountain people in general. This is true whether you are in the Alps or in Tenerife, or in the Himalayas. These people have some special qualities. They are hard-working, they think of their posterity, they show great solidarity, they can also trust you immediately but you must respect the enormous work that they and their ancestors have done to make the mountain livable.

\medskip
PPP: You hiked in the Himalayas?

\medskip
BT: To call that hiking is quite an exaggeration. But in 1967 Maryvonne and I went from Afghanistan to the valley of Hunza in Northern Pakistan. I could not really hike because I had twisted my ankle. But we met people there. 

Mountain is important for me also because of its silence and vistas. It is a special kind of atmosphere, of fusion with the world, which I enjoy very much.  Two or three years ago I went up to my hut, which is very close to where the  picture you show above was taken, and that year there was a drought. Usually, you can hear the stream but not that time. I was there at the end of the afternoon, there was no wind, it was absolute silence and a view like in the picture. There are moments like this in the mountains. That is I think a wonderful gift, which you cannot find anywhere else. 

I am very grateful to Maryvonne for sharing the effort of rebuilding, over many summers, an abandoned shepherds' hut 550 meters  in altitude above the road, when there were no helilifts but only our backs and donkeys lent to us by shepherds. Of course we also shared the pleasures of living there in the summer for more than 50 years and teaching mountain life and values to our children.

\medskip
PPP: I have the impression that you also like aphorisms. I remember one from my PhD years. I complained that I did not understand mathematics. You quoted a great artist, I don't remember who, saying at seventy or eighty that he began to understand what it meant to paint. 

\medskip
BT: That's Matisse. 

\medskip
PPP: Ah, thank you. So, now that you are almost eighty, do you feel you begin to understand what is mathematics? 

\medskip
BT: Unfortunately, I don't think I can say that. But I still try.

\medskip
PPP: Thank you very much Bernard for this wonderful interview!

\bigskip

{\bf Acknowledgments.} This research was funded, in whole or in part, by l'Agence Nationale de la
Recherche (ANR), project SINTROP (ANR-22-CE40-0014) and Labex CEMPI (ANR-11-LABX-0007-
01).   I am grateful to Evelia Rosa Garc\'{\i}a Barroso, Pedro Daniel Gonz\'alez P\'erez and Paul-Emmanuel Timotei for their remarks on a previous version of this interview.

\vfill 
\noindent
\textbf{\small{Authors' addresses:}}
\smallskip
\

\noindent
\small{P.\ Popescu-Pampu,
  Univ.~Lille, CNRS, UMR 8524 - Laboratoire Paul Painlev{\'e}, F-59000 Lille, France.
  \\
\noindent \emph{Email address:} \url{patrick.popescu-pampu@univ-lille.fr}}
\vspace{2ex}

\medskip
\end{document}